\documentclass[12pt]{amsart}
\usepackage{amsmath,amssymb,amscd,array, mathrsfs }

\usepackage{amsmath,amscd,amssymb,amsthm,array}
\usepackage{color}

\usepackage{amsmath,amssymb,mathrsfs,amsthm, tikz-cd,mathrsfs}
\usepackage{comment}
\def\url#1{\expandafter\s

\tring\csname #1\endcsname}

\def\mmat #1,#2,#3,#4,{\text{\small\arraycolsep=3pt $
\begin{pmatrix}#1&#2\\#3&#4\end{pmatrix}$}}

\usepackage{hyperref}
\usepackage{capt-of}

\usepackage{multirow}
\usepackage[all]{xy}

\usepackage{comments}
\newComments\SBe{Said}{blue}
\newComments\SBo{Sofiane}{blue}
\newComments\AM{Nacer}{blue}
\newComments\DL{DL}{red}
\newComments\QEh{QEh}{blue}

\def\mmat #1,#2,#3,#4,{\text{\small\arraycolsep=3pt $
\begin{pmatrix}#1&#2\\#3&#4\end{pmatrix}$}}

\usepackage{lscape}
\usepackage{tikz-cd}
\usepackage{enumerate}

\usepackage{DLdef1}
\usepackage{multicol}

\def\mmat #1,#2,#3,#4,{\text{\small\arraycolsep=3pt $
\begin{pmatrix}#1&#2\\#3&#4\end{pmatrix}$}}

\renewcommand {\ssbegin}[2][*]
 {\refstepcounter{subsection}%
\if#1*
\addcontentsline{toc}{subsection}{\thesubsection.\hspace*{1pc} #2}%
\else
\addcontentsline{toc}{subsection}{\thesubsection.\hskip 1pc #2. #1}%
\fi
 \def \secno {\gdef \secno {}{\ssecfont
\thesubsection.\hskip 2ex}%
 }%
 \begin{#2}}

\renewcommand {\sssbegin}[2][*]
 {\refstepcounter{subsubsection}
\if#1*
\addcontentsline{toc}{subsubsection}{\thesubsubsection.\hskip 1pc #2}%
\else
\addcontentsline{toc}{subsubsection}{\thesubsubsection.\hskip 1pc #2. #1}
\fi
 \def \secno {\gdef \secno {}{\ssecfont \thesubsubsection.\hskip 2ex}%
 }%
 \begin{#2}}

\renewcommand {\parbegin}[2][*]
 {\refstepcounter{paragraph}
\if#1*
\addcontentsline{toc}{paragraph}{\theparagraph.\hskip 1pc #2}%
\else
\addcontentsline{toc}{paragraph}{\theparagraph.\hskip 1pc #2. #1}
\fi
 \def \secno {\gdef \secno {}{\ssecfont \theparagraph.\hskip 2ex}%
 }%
 \begin{#2}}

\rmnameii{etr}{etr}
\rmnameii{evv}{ev}
\DeclareMathOperator{\h}{\mathcal{H}}

\DeclareMathOperator{\K}{\mathbb{K}}
\newcommand{\Z}{\mathbb{Z}}

\newcommand {\A}{{\cal{A}}}

\newcommand{\g}{\mathfrak{g}}
\newcommand{\Ll}{{\mathrm{L}}}
\newcommand{\Rr}{{\mathrm{R}}}

\newcommand{\al}{\alpha}
\newcommand{\ga}{\gamma}

\newcommand{\la}{\lambda}

\newcommand{\prs}{\langle\;,\;\rangle}
\def\br{[\;,\;]}
\newcommand{\esp}{\quad\mbox{and}\quad}
\begin{document}

\title[Pseudo-Euclidean Novikov superalgebras ]{Pseudo-Euclidean Novikov Superalgebras: Structure and Properties}

\author{Sa\"id Benayadi}
\address {Universit\'e  de Lorraine, Laboratoire IECL, CNRS-UMR 7502, UFR MIM, 3 rue Augustin Fresnel, BP 45112, 57073 Metz Cedex 03, France}
\email{said.benayadi@univ-lorraine.fr}

\author{Sofiane Bouarroudj}

\address{Division of Science and Mathematics, New York University Abu Dhabi, P.O. Box 129188, Abu Dhabi, United Arab Emirates.}
\email{sofiane.bouarroudj@nyu.edu}

\author{Hamza El Ouali}
\address{Universit\'e Cadi-Ayyad,
	Facult\'e des sciences et techniques,
	B.P. 549, Marrakech, Maroc.}
\email{eloualihamza11@gmail.com}

\thanks{S. Bo. was supported by the grant NYUAD-065}



\begin{abstract}
A pseudo-Euclidean Novikov superalgebra $\A$ is a Novikov superalgebra endowed with a  non-degenerate symmetric bilinear form $\prs$ such that all left multiplication operators are $\prs$-antisymmetric. In this case, the associated Lie superalgebra $(\A^{-},\prs)$ is a flat pseudo-Euclidean Lie superalgebra. In this paper, we investigate the structure of pseudo-Euclidean Novikov superalgebras. In particular, we introduce a distinguished subclass, called \emph{Milnor  superalgebras}, and prove that any pseudo-Euclidean Novikov superalgebra whose two-sided ideal is non-degenerate belongs to this class. We provide a method for constructing pseudo-Euclidean Novikov superalgebras.

We also introduce a double extension procedure for pseudo-Euclidean Novikov superalgebras and show that every such superalgebra with a degenerate two-sided ideal can be obtained via this method. Furthermore, we establish that any pseudo-Euclidean Novikov superalgebra is either a Milnor superalgebra or can be obtained by a sequence of double extensions starting from a Milnor superalgebra. As an application, we provide a complete classification of pseudo-Euclidean Novikov superalgebras of total dimension at most four over the field of real numbers.

\end{abstract}


\maketitle

\thispagestyle{empty}
\setcounter{tocdepth}{2}
\section{Introduction}
\subsection{Flat pseudo-Euclidean Lie (super)algebras}
A pseudo-Riemannian Lie group is a Lie group \(G\) endowed with a left-invariant pseudo-Riemannian metric \(\mu\). Its Lie algebra \(\mathfrak{g} = T_eG\), equipped with the bilinear form \(\prs = \mu_e\), is called a pseudo-Riemannian (or pseudo-Euclidean) Lie algebra. It is well known that there exists a unique torsion-free connection compatible with the metric, called the \emph{Levi-Civita connection}. This connection induces a bilinear product on \(\mathfrak{g}\), called the \emph{Levi-Civita product}, defined by the Koszul formula:
\[
2\langle u \bullet v, w \rangle
= \langle [u,v], w \rangle
- \langle [v,w], u \rangle
+ \langle [w,u], v \rangle,
\qquad \text{ for all } \, u,v,w \in \mathfrak{g}.
\]
This product satisfies the identities
\[
u \bullet v - v \bullet u = [u,v], \qquad
\langle u \bullet v, w \rangle + \langle v, u \bullet w \rangle=0,
\qquad \text{ for all } \, u,v,w \in \mathfrak{g}.
\]
In the framework of pseudo-Euclidean non-associative (super)algebras, and in particular Leibniz (super)algebras, this construction has been extended in \cite{BBE}. See also \cite{BoO, BoO2} for the antisymmetric case. 

The curvature tensor of the Levi-Civita connection at the identity element $e$ is given by
\[
\mathcal{R}(u,v) := \Ll_{[u,v]} - [\Ll_u, \Ll_v],
\]
where \(\Ll_u\) denotes the left multiplication operator defined by
$
\Ll_u(v) = u \bullet v,$ for all $u,v \in \mathfrak{g}.$

A pseudo-Riemannian Lie group \((G,\mu)\) is said to be \emph{flat} if its curvature $\mathcal{R}$ vanishes identically. In this case, the associated Lie algebra \((\mathfrak{g}, \br,  \prs)\) is called a \emph{flat pseudo-Riemannian \textup{(}or pseudo-Euclidean\textup{)} Lie algebra}. The vanishing of the curvature is equivalent to the fact that \((\mathfrak{g}, \bullet)\) is a left-symmetric algebra. Consequently, a flat pseudo-Riemannian Lie algebra can be viewed as a left-symmetric algebra endowed with a symmetric non-degenerate bilinear form for which all left multiplication operators are antisymmetric with respect to this form.

In the case where the metric is Riemannian (or Euclidean), that is, positive definite (or negative definite), J.~Milnor \cite{M} proved that a Euclidean Lie algebra $(\mathfrak{g}, \br, \prs)$ is flat if and only if $\mathfrak{g}$ decomposes orthogonally as $
\mathfrak{g} = \mathfrak{b} \oplus \mathfrak{u},
$
where $\mathfrak{u}$ is an abelian ideal, $\mathfrak{b}$ is an abelian subalgebra, and $\mathrm{ad}_u$ is $\prs$-antisymmetric, for all $u\in \mathfrak{b}$.

In recent years, several works have been devoted to the study of flat pseudo-Euclidean Lie algebras with arbitrary signatures. In the Lorentzian case, that is, for signature $(1,n-1)$, it was shown in  \cite{DI} that any flat Lorentzian Lie algebra is necessarily solvable. In \cite{Medina2}, the authors introduced a construction method, called the \emph{double extension}, for flat pseudo-Euclidean Lie algebras, which is an analogue of the double extension procedure for quadratic Lie algebras developed by A.~Medina and P.~Revoy \cite{MR}. They proved that flat Lorentzian nilpotent Lie algebras can be obtained via this method and that they are at most $3$-step nilpotent. Furthermore, M.~Gu\'ediri \cite{guediri} classified $2$-step nilpotent flat Lorentzian Lie algebras and showed that such algebras are trivial extensions of the $3$-dimensional Heisenberg Lie algebra $\mathfrak{h}_3$, where the restriction of the metric to $\mathfrak{h}_3$ is Lorentzian,  and $\mathfrak{h}_3$ has a  degenerate center. Later, M.~Ait Ben Haddou, M.~Boucetta and H.~Lebzioui \cite{ABL} further investigated flat Lorentzian Lie algebras. They proved that flat Lorentzian Lie algebras with degenerate center can be obtained by the double extension method. In addition, it was shown in \cite{BL1} that non-unimodular flat Lorentzian  Lie algebras can also be constructed using this approach. A complete classification of flat Lorentzian nilpotent Lie algebras was later obtained in \cite{Bajo}.

For signature $(2,n-2)$, M.~Boucetta and H.~Lebzioui \cite{BL2} studied flat pseudo-Euclidean Lie algebras and showed that every nilpotent flat pseudo-Euclidean Lie algebra of this signature arises as a double extension of a nilpotent flat Lorentzian Lie algebra. 
Moreover, in \cite{LH}, the author proved that any unimodular flat pseudo-Euclidean Lie algebra of signature $(2,n-2)$ with degenerate center arises as a double extension of a unimodular flat Lorentzian Lie algebra.

More recently, in \cite{BBE}, the authors introduced the notion of the Levi-Civita product associated with pseudo-Euclidean non-associative (super)algebras, in particular left Leibniz (super)algebras, together with a corresponding notion of curvature. They investigated flat pseudo-Euclidean left Leibniz (super)algebras and developed a double extension procedure in this setting. In particular, they proved that any flat pseudo-Euclidean left Leibniz (super)algebra can be obtained by a sequence of double extensions starting from a flat pseudo-Euclidean Lie (super)algebra. Consequently, the classification of flat pseudo-Euclidean left Leibniz (super)algebras reduces to the classification of flat pseudo-Euclidean Lie (super)algebras.
\subsection{Novikov (super)algebras} Novikov algebras constitute an important subclass of left-symmetric algebras. A Novikov algebra is a left-symmetric algebra $(\A,\bullet)$ satisfying $$
(u \bullet v)\bullet w = (u \bullet w)\bullet v,\quad \text{ for all $u,v,w \in \A.$}
$$
It is well known that the commutator
$
[u,v] = u \bullet v - v \bullet u,$ for all $u,v \in \A,$ defines a Lie algebra structure on $\A$, denoted by $\A^{-}$.

In \cite{LH2}, the author introduced the notion of a pseudo-Euclidean Novikov algebra $(\A,\bullet,\prs)$, that is, a Novikov algebra endowed with a non-degenerate symmetric bilinear form such that all left multiplication operators are $\prs$-antisymmetric. It follows that $(\A,\bullet,\prs)$ is a pseudo-Euclidean Novikov algebra if and only if $(\A^{-},\prs)$ is a flat pseudo-Euclidean Lie algebra and
$
[\Rr_u, \Rr_v] = 0,$ for all $ u,v \in \A,
$
where $\Rr_u$ is the right multiplication operator, and $\bullet$ coincides with the Levi-Civita product associated with $(\A^{-},\prs)$.

In \cite{LH2}, the author also studied pseudo-Euclidean Novikov algebras in the Lorentzian case. In particular, using Milnor’s theorem, it was shown that the Levi-Civita product of a flat Euclidean Lie algebra is Novikov. Moreover, it was proved that a pseudo-Euclidean Novikov algebra $(\A,\bullet,\prs)$ such that the derived ideal $[\A^{-},\A^{-}]$ is  non-degenerate satisfies the orthogonal decomposition
$
\A^{-} = [\A^{-},\A^{-}]^{\perp} \oplus [\A^{-},\A^{-}],
$
where both components are abelian, $[\A^{-},\A^{-}]$ is non-degenerate, and $\mathrm{ad}_u$ is $\prs$-antisymmetric for all $u \in [\A^{-},\A^{-}]^{\perp}$. In particular, such algebras are $2$-step solvable. In the degenerate case, it was shown that these algebras can be obtained via the double extension process from flat Euclidean Lie algebras.

In \cite{BEL}, the authors generalized these results to arbitrary signature. They proved that pseudo-Euclidean Novikov algebras are $2$-step solvable and established analogous structural decompositions when the derived ideal is non-degenerate. Furthermore, they introduced a double extension by blocks for pseudo-Euclidean Novikov algebras and showed that, in the degenerate case, such algebras are obtained by a  double extension starting from pseudo-Euclidean Novikov  algebras such that the derived ideal $[\A^{-},\A^{-}]$ is non-degenerate .

In a related direction, H.~Lebzioui \cite{LH1} studied the class $\mathcal{C}$ of flat pseudo-Euclidean Lie algebras $(\mathfrak{g},\br,\prs)$ such that the derived ideal $[\mathfrak{g},\mathfrak{g}]$ is contained in the left annihilator of the Levi-Civita product. In fact, Novikov algebras are subclasses of $\mathcal C$. This condition implies, in particular, that these Lie algebras are $2$-step solvable. Moreover, this is equivalent to $(\A,\bullet,\prs)$ being a pseudo-Euclidean left-symmetric $\Ll$-algebra, where $\A$ is the underlying vector space of $\mathfrak{g}$. It was shown that every Lorentzian algebra in $\mathcal{C}$ with degenerate derived ideal can be obtained by a flat double extension. In the nilpotent case, the derived ideal is necessarily degenerate, and the algebra is obtained by successive flat double extensions starting from an abelian Euclidean Lie algebra. In the non-degenerate case, a complete characterization was obtained for any signature. The study of this class was further developed in \cite{BO}, where the authors provided an inductive description of all elements of $\mathcal{C}$ in arbitrary signature by introducing double extensions by algebras of dimension $1$ and $2$. They also showed that pseudo-Euclidean Novikov algebras form a subclass of pseudo-Euclidean left-symmetric $\Ll$-algebras. In addition, in \cite{BO1}, the authors investigated Lie-admissible algebras whose product is a biderivation of the associated Lie algebra, and studied in particular pseudo-Euclidean Novikov algebras satisfying this condition.

Zelmanov \cite{Z} studied Novikov algebras endowed with a positive definite symmetric bilinear form for which the right multiplications are symmetric. He proved that such algebras are necessarily associative and provided their classification. The Lorentzian case was later investigated by Gu\'ediri \cite{G}, who obtained a complete classification in that setting.

In \cite{AB, JZ}, the authors studied Novikov superalgebras $(\mathcal{A}, \bullet)$ (see the definition of Novikov superalgebras in the Background section) endowed with an even  nondegenerate symmetric bilinear form $\prs$ which is invariant, that is,
\[
\langle u \bullet v, w \rangle  = \langle u, v \bullet w \rangle , \quad \text{ for all } u,v,w \in \mathcal{A}.
\]
Such structures are called \emph{symmetric Novikov superalgebras} or \emph{quadratic Novikov superalgebras}. They showed that such superalgebras are associative satisfying an additional condition, and that the associated Lie superalgebra $(\mathcal{A}^-, \prs)$ is a quadratic $2$-step nilpotent Lie superalgebra. Several examples are also provided. Moreover, the authors in \cite{AB} introduced the notion of double extensions in order to obtain inductive descriptions of symmetric Novikov superalgebras.

In this paper, we study Novikov superalgebras endowed with a non-degenerate symmetric bilinear form for which left multiplications are antisymmetric.  

\subsection{On the results}
Section \ref{Myback} collects the preliminary material required for our work. For completeness and to establish notation, we include a review of key concepts, some of which may be familiar to specialists. We introduce the notions of $\Ll$- and $\Rr$-superalgebras. If both conditions are satisfied, we define the notion of an $\Ll\Rr$-superalgebra. We also show that any $\Ll\Rr$-superalgebra is Lie admissible. Moreover, we show that the corresponding  Lie superalgebra is $2$-step solvable, see Proposition \ref{My LR}.

In section \ref{Mypseudo}, we introduce the concept of pseudo-Euclidean Novikov superalgebra, see Definition \ref{My Nov}. These are Novikov superalgebras endowed with a non-degenerate symmetric bilinear form $\prs$ such that all left multiplication operators are $\prs$-antisymmetric.  We show that $(\A,\bullet, \prs)$ is a pseudo-Euclidean Novikov superalgebra if and only if $(\A^{-},\prs)$ is a flat pseudo-Euclidean Lie superalgebra with the property that 
\[
[\Rr_u, \Rr_v]=0, \qquad \text{ for all } \, u,v \in \A,
\]
where $\Rr_u$ denotes the right multiplication operator. In this case, the product $\bullet$ coincides with the Levi-Civita product associated with $(\A^{-},\prs)$. This equivalence allows us to reduce the study of pseudo-Euclidean Novikov superalgebras to that of flat pseudo-Euclidean Lie superalgebras together with commuting right multiplication operators.

First, we establish a key structural result (see Lemma~\ref{Lemma11}): if $(\A,\bullet)$ is a superalgebra endowed with a non-degenerate  symmetric bilinear form such that all left multiplication operators are antisymmetric and the right multiplication operators commute, then
\[
(u \bullet v)\bullet w = 0, \qquad \text{ for all }\, u,v,w \in \A.
\]
We then prove that a pseudo-Euclidean superalgebra is Novikov if and only if it is an $\Ll\Rr$-superalgebra if and only if it is a left Leibniz $\Ll$-superalgebra (see Proposition~\ref{Pr1}).

Furthermore, by applying Engel's theorem, we show that a pseudo-Euclidean Novikov superalgebra $(\A,\bullet)$ is nilpotent if and only if all left multiplication operators $\Ll_u$ are nilpotent for every $u \in \A$, which is also equivalent to the nilpotency of the Lie superalgebra $\A^{-}$ (see Proposition~\ref{Nil}).

While studying the structure theory of pseudo-Euclidean Novikov superalgebras, we are led to introduce a distinguished subclass, called \emph{Milnor}  superalgebras (see Definition \ref{defmilnor} and  Proposition~\ref{milnor}). In particular, we prove that any pseudo-Euclidean Novikov superalgebra whose two-sided ideal is non-degenerate is necessarily a Milnor superalgebra (see Theorem~\ref{ ideal non degenere}). Furthermore, we show that if $(\A,\bullet,\prs)$ is a pseudo-Euclidean Novikov superalgebra such that $\A \bullet \A$ is non-degenerate, then $\A^{-}$ is nilpotent if and only if the product $\bullet$ is trivial, see Corollary \ref{nilpotent}. 

We  introduce the notion of representation of pseudo-Euclidean Novikov superalgebras  and define the associated product. More precisely, given such a superalgebra $(\A, \bullet, \prs)$, we define a new bilinear product $\star$ on $\A$ by means of the relation
$
\langle u \bullet v, w \rangle = \langle u, v \star w \rangle,$ for all $ u,v,w \in \A,$ and we show that this product endows $\A$ with a structure of a $2$-step nilpotent Lie superalgebra. Furthermore, we establish a characterization of pseudo-Euclidean Novikov superalgebras in terms of bimodule isomorphisms. More precisely, we prove that the existence of a suitable bilinear form is equivalent to the isomorphism between certain naturally associated $\A$-bimodules (see Theorem \ref{Thm7.11}). Finally, motivated by this representation and inspired by the $T^*$-extension for quadratic algebras introduced in \cite{Bordemann}, we provide a method for constructing pseudo-Euclidean Novikov superalgebras (see Theorems \ref{constriction} and \ref{constriction1}).

 In Section \ref{d-ext}, we introduce the notion of \emph{double extension} for pseudo-Euclidean Novikov superalgebras, see Theorem \ref{double-ex1}, Theorem \ref{double-ex2}, Theorem \ref{double-ex3}, and Theorem \ref{double-ex4}. This construction is inspired by the classical double extension, originally introduced by Medina and Revoy \cite{MR} for quadratic Lie algebras.

 In section \ref{con-ex}, we show that every pseudo-Euclidean Novikov superalgebra whose two-sided ideal is degenerate can be obtained by a double extension (see Theorems~\ref{central} and \ref{central1}). 
More precisely, every non-trivial even (resp. odd) pseudo-Euclidean nilpotent Novikov superalgebra can be obtained by a finite sequence of double extensions of  even (resp. odd)  pseudo-Euclidean nilpotent Novikov superalgebras, starting from the trivial superalgebra (see Corollary~\ref{start-abel} and Corollary~\ref{start-abel2}). 

Finally, we prove that any pseudo-Euclidean Novikov superalgebra is either a Milnor superalgebra or can be obtained by a sequence of double extensions from a Milnor superalgebra.

Section~\ref{classi} is devoted to the classification of all pseudo-Euclidean Novikov superalgebras of total dimension at most four. The classification is carried out over the field of real numbers.


\section{Backgrounds}\label{Myback}

Let $\mathbb{K}$ be an arbitrary field of characteristic zero. The group of integers modulo $2$ is denoted by $\Z_{2}$. 

 Let $V=V_{\bar 0}\oplus V_{\bar 1}$ be a superspace  defined over ${\mathbb K}$. The parity of a homogeneous element $v\in V_{\bar{i}}$ is denoted by $|v|:=\bar{i}$. The element $v$ is called \textit{even} if $v\in V_{\bar 0}$ and \textit{odd} if $v\in V_{\bar 1}$. The {\it superdimension} of $V$ is $\mathrm{sdim}(V)=a+b\epsilon$, where $\epsilon^2=1$, and $a=\mathrm{dim}(V_{\bar 0})$, $b=\mathrm{dim}(V_{\bar 1})$. Usually, $\mathrm{sdim}(V)$ is shorthanded as $a\mid b$. The total dimension of $V$, which is $a+b$, is denoted by  $\mathrm{dim}(V)$.  
 
 Throughout the text, all elements are supposed to be homogeneous unless otherwise stated. 
 
 A linear map $\varphi:V\rightarrow W$ between superspaces is called \textit{even} if $\varphi(V_{\bar{i}})\subset W_{\bar{i}}$ and \textit{odd} if $\varphi(V_{\bar{i}})\subset W_{\overline{i}+\overline{1}}$.

Let $V=V_{\bar 0}\oplus V_{\bar 1}$ be a vector superspace. We denote by $\Pi$ the \textit{change of parity functor} $\Pi: V\mapsto \Pi (V)$, where $\Pi(V)$ is another copy of $V$ such that $ \Pi(V)_{\bar 0}:=V_{\bar 1};~~\Pi(V)_{\bar 1}:=V_{\bar 0}$.  Elements of $\Pi(V)$ shall be denoted by $\Pi(v),\text{ for all } v\in V$. In fact, $|\Pi(v)|=|v|+{\bar 1}$ for every homogeneous $v\in V$.

 Let $U$ and $V$ be two vector superspaces. Given a morphism $G: U \to V$, we define the morphism
$
G^\Pi : \Pi(U) \to \Pi(V)
$
defined by $G^\Pi(\Pi(u))=\Pi(G(u))$. It is straightforward to verify that this construction makes the assignment
$
U \mapsto \Pi(U)
$ into a functor satisfying $\Pi^2 = \mathrm{id}$.

Moreover, given a morphism $G: U \to V$, one can construct morphisms $
\Pi(G) : U \to \Pi(V), $ defined by $ \Pi(G)(u) = \Pi(G(u)),
$ and
$
G\Pi : \Pi(U) \to V,$ defined by  $(G\Pi)(\Pi(u)) = G(u).
$
Clearly, we have the factorization $
G^\Pi = \Pi(G\Pi)$, see
\cite{Youri}.

Let \(V\) be a vector superspace. For homogeneous elements \(x,y\in V\), we denote by \(x^*,y^*\in V^*\) their duals . For every \(x^*\otimes y^*\in V^*\otimes V^*\) and  \(u\otimes v\in V\otimes V\), we define the pairing
\[
\langle x^*\otimes y^*, u \otimes v\rangle :=(-1)^{|y||u|}x^*(u)\,y^*(v).
\]
The symmetric product (graded symmetrisation) is given by
\[
x^*\odot y^*=x^*\otimes y^* + (-1)^{|x||y|}\,y^*\otimes x^*.
\]
\subsection{Bilinear forms on a superspace}
Let $V$ and $W$ be two superspaces. Let $\beta \in \text{Bil}(V,W)$ be a homogeneous bilinear form. Recall that the Gram matrix $B=(B_{ij})$ associated with $\beta $ is given by the following formula:
\begin{equation*}\label{martBil}
B_{ij}=(-1)^{|\beta ||v_i|}\beta (v_{i}, w_{j})\text{~~for the basis vectors $v_{i}\in V$ and $w_{j}\in W$.}
\end{equation*}
This definition allows us to identify a~bilinear form of $\mathrm{Bil} (V, W)$ with an element of $\mathrm{Hom}(V, W^*)$. 
Consider the \textit{upsetting} of bilinear forms
$u:\text{Bil} (V, W)\rightarrow \text{Bil}(W, V)$ given by the formula 
\begin{equation*}\label{susyB}
u(\beta )(w, v)=(-1)^{|v||w|}\beta (v,w)\text{~~for any $v \in V$ and $w\in W$.}
\end{equation*}
In terms of the Gram matrix $B$ of $\beta $, we have 
\begin{equation*}
u(B)=
\left ( \begin{array}{cc}  R^{t} & (-1)^{|B|}T^{t} \\ (-1)^{|B|}S^{t} & -U^{t} \end{array} \right ),
\text{ for $B=\left (  \begin{array}{cc}R & S \\ T & U \end{array} \right )$}.
\end{equation*}
Suppose now that $V=W$. Following \cite{L}, we say that: 
\begin{itemize} \item[(i)] the form
$\beta $ is  \textit{symmetric} if  and only if $u(B)=B$;

\item[(ii)] the form $\beta $ is \textit{anti-symmetric} if and only if $u(B)=-B$.

\item[(iii)] the form $\beta $ is non-degenerate if $u\in V$ satisfies $\beta (v,u)=0,$ for any $v\in V$, then $u = 0$,
\item[(iv)] the form $\beta $ is even if $\beta (V_{\bar{0}}, V_{\bar{1}})=\beta (V_{\bar{1}},V_{\bar{0}})=\{0\}$,
 \item[(v)] the form $\beta $ is odd if  $\beta (V_{\al}, V_{\al}) = \{0\}$, for all $\al\in\mathbb{Z}_2,$ 

\end{itemize}

From now on, we will denote a bilinear form with $\prs$ instead. 

A \emph{pseudo-Euclidean} vector superspace $(V,\prs)$ is a vector superspace $V$ equipped with a non-degenerate symmetric bilinear form. It is called \emph{Euclidean} if the form $(V,\prs)$ is  positive definite.

Let $(V,\prs)$ be a pseudo-Euclidean vector superspace. Let $f$ be a homogeneous endomorphism of $V$ of parity $\al \in \{\bar{0},\bar{1}\}$.
\begin{itemize}
  \item[(i)] The adjoint of $f$ with respect to $\prs$, denoted by $f^*$, is the unique homogeneous endomorphism of parity $\al$ satisfying 
  \[
  \langle f(u), v\rangle = (-1)^{\al|u|} \langle u, f^*(v)\rangle, \quad \text{for all } u,v\in V.
  \]
\item[(ii)] The endomorphism $f$ is said to be:
  \begin{itemize}
    \item Symmetric (or self-adjoint) with respect to $\prs$ if $f^* = f$, i.e.,
    \[
    \langle f(u), v\rangle = (-1)^{\al|u|} \langle u, f(v)\rangle, \quad \text{for all } u,v\in V.
    \]

    \item Anti-symmetric  (or anti-self-adjoint) with respect to $\prs$ if $f^* = -f$, i.e.,
    \[
     \langle f(u), v\rangle = -(-1)^{\al|u|} \langle u, f(v)\rangle, \quad \text{for all } u,v\in V.
    \]
  \end{itemize}
\end{itemize}
Instead of saying that $f$ is symmetric (resp. anti-symmetric) with respect to $\prs$, we adopt the terminology that $f$ is $\prs$-symmetric (resp.  $\prs$-anti-symmetric).


\ssbegin{Lemma}
\label{evendimension}
Let $V = V_{\bar{0}}\oplus V_{\bar{1}}$ be a finite-dimensional  vector superspace.
\begin{enumerate}
    \item[$(i)$] $V$ admits an even pseudo-Euclidean bilinear form $\prs$ if and only if $\dim V_{\bar{1}}$ is even.
    \item[$(ii)$] $V$ admits an odd pseudo-Euclidean  bilinear form $\prs$ if and only if $\dim V_{\bar{0}} = \dim V_{\bar{1}}$.
\end{enumerate}
\end{Lemma}

\begin{proof}
(i) $(\Rightarrow)$ The result follows from the fact that the restriction $\prs|_{V_{\bar 1}}$ is
non-degenerate and anti-symmetric. It is well known that, in the non-super setting, such a bilinear form can exist only on vector spaces of even dimension.

$(\Leftarrow)$ Conversely, suppose $\dim V_{\bar{1}}$ is even.  Choose a non-degenerate symmetric bilinear form on $V_{\bar{0}}$ and a non-degenerate skew-symmetric bilinear form on $V_{\bar{1}}$ (which exists since $\dim V_{\bar{1}}$ is even). Their orthogonal direct sum defines an even non-degenerate symmetric bilinear form on $V$.

(ii) $(\Rightarrow)$ Let $0 \neq u \in V_{\bar 1}$. Since $\prs$ is odd and non-degenerate, there exists $v \in V_{\bar 0}$ such that $\langle u,v\rangle =1$. Set $E:=\operatorname{Span}\{u,v\}$.
We claim $E^\perp \cap E=\{0\}$. Indeed, let $w \in E^\perp \cap E$.
Put $w=\alpha u+\beta v$ for some $\alpha,\beta \in \mathbb{K}$. We have
\[
0=\langle u,w\rangle=\langle u,\alpha u+\beta v\rangle=\beta,
\]
and similarly,
\[
0=\langle v,w\rangle =\langle v,\alpha u+\beta v\rangle=\alpha.
\]
Thus $w=0$, proving the claim.

Now, let $x\in V$. We can write
\[
x=\bigl(x-\langle x,u\rangle v-\langle x,v\rangle u\bigr)+\bigl(\langle x,u\rangle v+\langle x,v\rangle u\bigr)\in E^\perp \oplus E.
\]
Hence,
\(
V=E\oplus E^\perp,
\)
and therefore
\(
\operatorname{sdim}(V)=\operatorname{sdim}(E^\perp)+(1|1).
\)
The restriction $\prs|_{E^\perp}$ is again odd and non-degenerate, so the same argument
can be applied inductively to $E^\perp$. Since $V$ is finite-dimensional, the result
follows.

$(\Leftarrow)$ Conversely, suppose $\dim V_{\bar{0}} = \dim V_{\bar{1}} = n$. Let $\{e_1, \dots, e_n\}$ and $\{f_1, \dots, f_n\}$ be bases of $V_{\bar{0}}$ and $V_{\bar{1}}$, respectively. 
Defining a bilinear form $\prs$ by
\[
\langle e_i, f_j \rangle = \langle f_j, e_i \rangle = \delta_{ij} \quad \text{and} \quad \langle e_i, e_j \rangle = \langle f_i, f_j \rangle = 0 \quad \text{for all } 1 \le i, j \le n,
\]

This yields an odd pseudo-Euclidean structure on $V$.
\end{proof}

\subsection{Novikov  superalgebras}
A superalgebra $(\A, \bullet)$  over a field $\mathbb{K}$ consists of a superspace $\A=\A_{\bar 0}\oplus \A_{\bar 1}$ and a binary operation satisfying: $\A_\alpha \bullet \A_\beta \subseteq \A_{\alpha+\beta}, \text{ for } \alpha, \beta \in \mathbb{Z}_2$. 

 For every $u \in \A$, we denote the left and right multiplication operators by $\Ll_u$ and $\Rr_u$, respectively. More precisely, $\Ll_u$ and $\Rr_u$ are defined as  follows: 
 \[
 \Ll^\bullet_u(v):=u \bullet v \text{  and }\Rr^\bullet_u(v):=(-1)^{|u||v|}v\bullet u, \text{ for all $u, v \in \A$}.
 \]
Let us define the following normalizers
$$ N_\ell(\A, \bullet) := \left\{ u \in \A,\;  \mathrm{L}^\bullet_u = 0 \right\}, \;  N_r(\A, \bullet) := \left\{ u \in \A,\;  \mathrm{R}_u^\bullet = 0 \right\}, \; \text{and} $$
\[
N(\A, \bullet):= N_\ell(\A, \bullet)\cap  N_r(\A, \bullet).
\]
We refer to these as the {\it left normalizer}, {\it right normalizer}, and (two-sided) {\it normalizer}, respectively.

Let $(\A, \bullet)$ be a non-associative superalgebra. On the underlying vector superspace $\A$, we
define the corresponding commutator 
$$[u, v]_\bullet :=u\bullet v -(-1)^{|u||v|} v\bullet u,\; \text{ for all } u,v\in\A.$$ 
The superalgebra $(\A, [\; , \; ]_\bullet)$ will be denoted by  $\A^-$. The superalgebra $(\A, \bullet)$ is called Lie-admissible
superalgebra, if $\A^-$ is a Lie superalgebra.  In this case,  $\A^-$ is called the sub-adjacent Lie superalgebra of $(\A, \bullet)$ and the latter is called a compatible (Lie-admissible) superalgebra structure on
the Lie superalgebra $\A^-$. 

Similarly, on the underlying vector superspace $\cal A$, we define the corresponding anti-commutator:
$$[u, v]_\bullet^+ :=u\bullet v +(-1)^{|u||v|} v\bullet u,\; \text{ for all } u,v\in\A.$$ The superalgebra $({\cal A}, [\; , \;]_\bullet^+)$ will be denoted by ${\cal A}^+$.

A superalgebra $(\A, \bullet)$ is called left-symmetric (resp. right-symmetric) if and only if 
 \[
(u,v,w) = (-1)^{|u||v|} \, (v,u,w), \quad (\text{resp. } (w, u,v) = (-1)^{|u||v|} \, (w,v,u)), \text{ for all $u,v,w\in {\mathcal A}$,}
 \]
 where 
 $(u,v,w) :=(u \bullet v) \bullet w - u \bullet (v \bullet w)$  is the associator.  It is well-known that such superalgebras are Lie-admissible.

Let $(\mathcal{A}, \bullet)$ be a left-symmetric (resp. right-symmetric) superalgebra. Then the left (resp. right) multiplication operators satisfy the identity:
\begin{equation}
\mathrm{L}^\bullet_{[u,v]_\bullet} = \left[\mathrm{L}^\bullet_u, \mathrm{L}^\bullet_v\right]
\quad (\text{resp. } \mathrm{R}^\bullet_{[u,v]_\bullet} = -\left[\mathrm{R}^\bullet_u, \mathrm{R}^\bullet_v\right]), \quad \text{for all $u, v \in \A$.}
\label{mul}\end{equation}

It is obvious that $(\mathcal{A}, \bullet)$ is a left-symmetric superalgebra if and only if $\left(\mathcal{A}, \bullet_{\text {opp }}\right)$ is a right-symmetric superalgebra where  $u \bullet_{\text {opp }} v:=-(-1)^{|u||v|}v \bullet u$. 

Following \cite{Xu}, a superalgebra $(\A, \bullet)$ is called a Novikov superalgebra if it is a left-symmetric  superalgebra
and the following identity is satisfied:
$$(w\bullet u)\bullet v=(-1)^{|u||v|} (w\bullet v)\bullet u, \quad \text{for all $u,v,w\in {\cal A}$.}
$$

Let $(\A, \bullet)$ be a superalgebra. Then $(\A, \bullet)$ is a  Novikov  superalgebra if and only if 
\begin{align}\label{Nvleft}
&\mathrm{L}^\bullet_{[u,v]_\bullet} = \left[\mathrm{L}^\bullet_u, \mathrm{L}^\bullet_v\right],\esp  \Rr^\bullet_u\circ \Rr^\bullet_v=(-1)^{|u||v|} \Rr^\bullet_v\circ \Rr^\bullet_u,\quad \text{for all $u, v \in \A$.}
\end{align}
\sssbegin{Remark}
Novikov algebras admit two superizations:  Novikov superalgebras, and Balinsky-Novikov superalgebras. We will study Balinsky–Novikov superalgebras in a forthcoming paper. Novikov and Balinsky–Novikov structures on four-dimensional real Lie superalgebras have been investigated in \cite{BR}.
\end{Remark}

Let $(\mathcal{A},\bullet)$ be a Novikov superalgebra. For $n\geq 1$, we define the following sequences of subsuperspaces of $\mathcal{A}$:
\begin{align*}
\mathcal{A}^0 & =   \mathcal{A}, &  
\mathcal{A}^n & =  \mathcal{A}^{\,n-1}\bullet \mathcal{A}^{\,n-1}; & 
\mathcal{A}^{(0)} &= \mathcal{A}, & 
\mathcal{A}^{(n)} &= \mathcal{A}\bullet \mathcal{A}^{(n-1)};\\
\mathcal{A}^{\langle 0\rangle} &= \mathcal{A}, & 
\mathcal{A}^{\langle n\rangle} &= \mathcal{A}^{\langle n-1\rangle}\bullet \mathcal{A}; &
\mathcal{A}_1 & = \mathcal{A}, & 
\mathcal{A}_n &  = \sum_{i=1}^{n-1} \mathcal{A}_i \bullet \mathcal{A}_{n-i}, \quad n\ge 2.
\end{align*}
\begin{enumerate}
\item[(i)] If there exists $n\in\mathbb{N}$ such that $\mathcal{A}^n=\{0\}$, then $\mathcal{A}$ is said to be \emph{solvable}. 

\item[(ii)] If there exists $n\in\mathbb{N}$ such that $\mathcal{A}^{(n)}=\{0\}$, then $\mathcal{A}$ is called \emph{left nilpotent}.

\item[(iii)] If there exists $n\in\mathbb{N}$ such that $\mathcal{A}^{\langle n\rangle}=\{0\}$, then $\mathcal{A}$ is called \emph{right nilpotent}.

\item[(iv)] If there exists $n\in\mathbb{N}$ such that $\mathcal{A}_{n}=\{0\}$, then $\mathcal{A}$ is called \emph{ nilpotent}.
\end{enumerate}

A superalgebra $(\A, \bullet)$ is called an  \emph{$\Ll$-superalgebra} (resp. \emph{$\Rr$-superalgebra}) if and only if 
 \begin{equation}
u\bullet (v\bullet w)
= (-1)^{|u||v|} \, v\bullet (u\bullet w) , \quad (\text{resp. } (w\bullet u)\bullet v
= (-1)^{|u||v|} \, (w\bullet v)\bullet u)
 \label{LR-eq}\end{equation}
for all $u,v,w\in {\mathcal A}$.

A superalgebra $(\A,\bullet)$ is called an \emph{$\Ll\Rr$-superalgebra} if it is both an
$\Ll$- and an $\Rr$-superalgebra.

The following proposition can be viewed as a superization of a result by Burde, Deschamps, and Dekimpe \cite{BDD}.
\sssbegin{Proposition}\label{My LR}
If $(\A, \bullet)$ is an $\Ll\Rr$-superalgebra, then $\A^{-}$ is a Lie superalgebra. Moreover, $\A^{-}$ is $2$-step solvable.
\end{Proposition}
\begin{proof}
By a direct computation using the identities of $\Ll\Rr$-superalgebras, one can show that $\A^{-}$ is a Lie superalgebra.

We now show that $\A^{-}$ is $2$-step solvable. Let $u,v,w,z\in \A$. Using the identities $\Ll$ and $\Rr$, we obtain
\begin{align*}
(u\bullet v)\bullet (w\bullet z) 
&= (-1)^{|v|(|w|+|z|)}(u\bullet (w\bullet z))\bullet v \\
&= (-1)^{|v|(|w|+|z|)+|u||w|}(w\bullet (u\bullet z))\bullet v \\
&= (-1)^{|v||w|+|u||w|+|v||u|}(w\bullet v)\bullet (u\bullet z) \\
&= (-1)^{|v||w|}u\bullet ((w\bullet v)\bullet z) \\
&= (-1)^{|v||w|+|v||z|}u\bullet ((w\bullet z)\bullet v) \\
&= (-1)^{(|u|+|v|)(|w|+|z|)}(w\bullet z)\bullet (u\bullet v).
\end{align*}
Thus, $(u\bullet v)$ and $(w\bullet z)$ supercommute in $(\A, \bullet)$. Now, 
\begin{align*}
[[u,v], [w,z]]&= [u,v]\bullet [w, z]-(-1)^{(|u|+|v|)(|w|+|z|)} [w,z]\bullet [u,v]\\&= (u\bullet v)\bullet (w\bullet z)- (-1)^{|u||v|}(v\bullet u)\bullet (w\bullet z) -(-1)^{|w||z|} (u\bullet v)\bullet (z\bullet w)\\ &\quad+ (-1)^{|u||v|+ |w||z|}(v\bullet u)\bullet (z\bullet w)-(-1)^{(|u|+|v|)(|w|+|z|)} (w\bullet z)\bullet (u\bullet v)\\& \quad +(-1)^{(|u|+|v|)(|w|+|z|)} ((-1)^{|z||w|}(z\bullet w)\bullet (u\bullet v) +(-1)^{(|u||v|} (w\bullet z)\bullet (v\bullet u)) \\&\quad-(-1)^{(|u|+|v|)(|w|+|z|)+|z||w|+|u||v|} (z\bullet w)\bullet (v\bullet u)\\&=0.
\end{align*}
Thus, $\A^-$ is $2$-step solvable.
\end{proof}

Following \cite{AAO}, a superalgebra $(\A,\bullet)$ is called a \emph{left-Leibniz superalgebra} if
\begin{equation}
(u\bullet v)\bullet w
= u\bullet (v\bullet w)
- (-1)^{|u||v|} \, v\bullet (u\bullet w),
\quad \text{for all } u,v,w\in\A.
\label{Leibniz}\end{equation}

 \subsection{The Levi-Civita product}

 In \cite{BBE}, the notion of a Levi-Civita product for pseudo-Euclidean non-associative superalgebras was introduced. In this work, we only require the Levi-Civita product associated with pseudo-Euclidean Lie algebras, which we recall in the following Proposition.
 
\sssbegin{Proposition}[\cite{BBE}]\label{LVCT} Let $(\g,\br,\prs)$ be a pseudo-Euclidean Lie superalgebra. Then there exists a unique product $\bullet$ on $\g$ satisfying
\begin{align}
[u,v] &=[u,v]_\bullet, \label{torsion1}\\
\langle u\bullet v , w\rangle
&= -(-1)^{|u||v|}\langle v , u\bullet w\rangle. \label{compatibl1}
\end{align}
More precisely, the product $\bullet$ is given by
\begin{equation}\label{lv-lie}
2\langle u\bullet v , w\rangle
= \langle [u,v] , w\rangle
- (-1)^{|u||v|+|u||w|}\langle [v,w] , u\rangle
+ (-1)^{|v||w|+|u||w|}\langle [w,u] , v\rangle.
\end{equation}
The product $\bullet$ is called the \emph{Levi-Civita product}
associated with $(\g,\br,\prs)$.	
\end{Proposition}
According to the relations \eqref{compatibl1} and \eqref{torsion1}, for any $u \in \g$, the left multiplication $\Ll_u^\bullet$ is $\prs$-anti-symmetric, and the adjoint operator satisfies
$
\ad_u = \Ll_u^\bullet - \Rr_u^\bullet,
$
where $\ad_u: \g \to \g$ is defined by $\ad_u(v) = [u,v]$ for all $v\in \g$.

Following \cite{BBE}, the curvature operator on $\g$ is defined by
\[
\mathcal{R}(u,v) : = \Ll^\bullet_{[u,v]} - [\Ll^\bullet_u, \Ll^\bullet_v], \qquad \text{ for all } u,v \in \g,
\]
Moreover, if this curvature operator vanishes identically, i.e.
\[
\mathcal{R}(u,v) = 0, \quad \text{ for all } \, u,v \in \g,
\]
then $(\g,\br, \prs)$ is called a \emph{flat pseudo-Euclidean Lie superalgebra}.

It is easy to see that the curvature $\mathcal{R}$ of a pseudo-Euclidean Lie superalgebra $(\g, \br, \prs)$  vanishes identically if and only if $(\g, \bullet)$ is a left-symmetric superalgebra.

\sssbegin{Remark}
The vanishing of the curvature of $(\g,\br,\prs)$ is equivalent to the fact that 
$\g$, endowed with its Levi-Civita product, is a left-symmetric superalgebra. 
Thus, a flat pseudo-Euclidean Lie superalgebra can be viewed as a left-symmetric 
superalgebra equipped with a non-degenerate symmetric bilinear form with respect to which  
the left multiplications are $\prs$-anti-symmetric.

Indeed, assume that $(\g,\br,\prs)$ is a flat Lie superalgebra and let $\bullet$ 
denote its Levi-Civita product. The vanishing of the curvature of 
$(\g,\br,\prs)$ implies that $(\g,\bullet)$ is a left-symmetric superalgebra. 
Moreover, the Levi-Civita product satisfies
\[
\langle u\bullet v, w\rangle 
+ (-1)^{|u||v|}\langle v, u\bullet w\rangle =0, \text{ for all $u,v,w\in\g$,}
\]
which means that the left multiplications are 
$\prs$-anti-symmetric. Hence $(\g,\bullet)$ is a left-symmetric superalgebra 
equipped with a non-degenerate symmetric  bilinear form with respect to which the left 
multiplications are $\prs$-antisymmetric.

Conversely, assume that $(\g,\bullet)$ is a left-symmetric superalgebra 
equipped with a non-degenerate symmetric bilinear form $\prs$ such that
\[
\langle u\bullet v, w\rangle 
+ (-1)^{|u||v|}\langle v, u\bullet w\rangle =0,
\qquad \text{ for all } u,v,w\in\g .
\]
It is well known that $\g^{-}$ is a Lie superalgebra with Lie bracket
\[
[u,v]_\bullet=u\bullet v-(-1)^{|u||v|}v\bullet u .
\]
Moreover, the above identity shows that $\bullet$ satisfies the defining 
property of the Levi-Civita product. By uniqueness, $\bullet$ coincides with the Levi-Civita product associated with 
$(\g^{-},\prs)$. Consequently, $(\g^{-},\prs)$ is a flat pseudo-Euclidean 
Lie superalgebra.

\end{Remark}

\section{Pseudo-Euclidean Novikov superalgebras}\label{Mypseudo}

The main goal of this paper is to introduce the notion of pseudo-Euclidean Novikov superalgebras and to investigate their structures.

The problem we study can be approached in two equivalent ways.

\medskip
\noindent \textbf{Approach 1.}
Let $(\mathfrak g,[\, ,\, ],\langle\, ,\, \rangle)$ be a pseudo-Euclidean Lie superalgebra and let 
$\bullet$ denote its Levi-Civita product. We investigate Lie superalgebras for which 
$(\mathfrak g,[\, ,\, ],\langle\, ,\, \rangle)$ is flat and
\[
[\Rr_u,\Rr_v]=0,\qquad \text{ for all } u,v\in\mathfrak g.
\]

The flatness of $(\mathfrak g,[\, ,\, ],\langle\, ,\, \rangle)$ implies that $(\mathfrak g,\bullet)$ 
is a left-symmetric  superalgebra. Moreover, the additional condition 
$[\Rr_u,\Rr_v]=0$ implies that $(\mathfrak g,\bullet)$ is a Novikov superalgebra.

\medskip

\noindent \textbf{Approach 2 (reciprocal).}
Let $(\mathfrak g,\star)$ be a Novikov superalgebra equipped with a non-degenerate symmetric bilinear form 
$\langle\, ,\, \rangle$ such that
\[
\langle \Ll_u^\star v,w\rangle + (-1)^{|u||v|}\langle v, \Ll_u^\star w\rangle = 0,
\qquad \text{ for all } u,v,w\in\mathfrak g .
\]
Then $(\mathfrak g^{-},\langle\, ,\, \rangle)$ is a pseudo-Euclidean Lie superalgebra with Lie bracket
\[
[u,v]:=u\star v-(-1)^{|u||v|}v\star u .
\]
Moreover, by the uniqueness of the Levi-Civita product, $\star$ coincides with the Levi-Civita product 
associated with $(\mathfrak g^{-},\langle\, ,\, \rangle)$.

\medskip

\noindent \textbf{Conclusion.} The study of pseudo-Euclidean Lie superalgebras whose Levi-Civita structure is a Novikov superalgebra is therefore equivalent to the study of Novikov superalgebras equipped with a non-degenerate symmetric 
bilinear form satisfying
\[
\langle \Ll_u v , w\rangle + (-1)^{|u||v|}\langle v , \Ll_u w\rangle =0,
\qquad \text{ for all } u,v,w\in\mathfrak g ,
\]

\subsection{Pseudo-Euclidean Novikov superalgebras}
The following definition is a superization of a definition given by Lebzioui \cite{LH2}.
\sssbegin{Definition}\label{My Nov}
A \emph{pseudo-Euclidean Novikov superalgebra} is a triple
$(\mathcal A, \bullet, \prs)$, where
$(\mathcal A, \bullet)$ is a Novikov superalgebra  and
$\prs$ is a non-degenerate symmetric bilinear form on
$\mathcal A$ satisfying 
\[
\langle u \bullet v, w \rangle
= -(-1)^{|u||v|}\langle v, u \bullet w \rangle,
\quad \text{for all } u, v, w \in \mathcal A.
\]
This means that $(\Ll_u^\bullet)^*=-\Ll_u^\bullet$ for all $u\in \A.$
\end{Definition}
In the case where the form $\prs$ is even (resp. odd), we refer to these superalgebras as \emph{even} (resp. \emph{odd}) pseudo-Euclidean Novikov superalgebras.
\sssbegin{Remark}
Let $(\A,\bullet,\prs)$ be a pseudo-Euclidean Novikov
superalgebra. Then $(\A^-,\prs)$, where
\(
\A^-=(\A,\br_\bullet)
\)
denotes the Lie superalgebra associated with $\bullet$, is a
flat pseudo-Euclidean Lie superalgebra. Moreover, the product $\bullet$
coincides with the Levi-Civita product associated with
$(\A^-,\prs)$.

Indeed, since $(\A,\bullet)$ is a Novikov superalgebra,
$\A^-=(\A,\br_\bullet)$ is a Lie superalgebra. Furthermore, the
pseudo-Euclidean condition implies that
\[
\langle u\bullet v, w\rangle
= -(-1)^{|u||v|}\langle v, u\bullet w\rangle,
\quad \text{for all } u,v,w\in\A.
\]
Hence, the product $\bullet$ satisfies the torsion-free condition
\eqref{torsion1} and the metric compatibility condition
\eqref{compatibl1}. By the uniqueness of the Levi-Civita product,
$\bullet$ is precisely the Levi-Civita product associated with
$(\A^-,\prs)$.

Since $(\A,\bullet)$ is Novikov, it follows that $(\A^-,\prs)$
is a flat pseudo-Euclidean Lie superalgebra and satisfies
\[
[\Rr_u^\bullet,\Rr_v^\bullet]=0,
\quad \text{for all } u,v\in\A.
\]

Conversely, let $(\g,\br,\prs)$ be a flat
pseudo-Euclidean Lie superalgebra whose Levi-Civita product $\bullet$
satisfies
\[
[\Rr_u^\bullet,\Rr_v^\bullet]=0,
\quad \text{for all } u,v\in\g.
\]
Then $(\g,\bullet)$ is a Novikov superalgebra. Moreover, since the left
multiplication operators of the Levi-Civita product are $\prs$-antisymmetric, the triple
$(\g,\bullet,\prs)$ is a pseudo-Euclidean Novikov
superalgebra.
\end{Remark}

We now present a method for constructing pseudo-Euclidean Novikov superalgebras. As we shall see, this class is quite rich and abundant, which further motivates and justifies a systematic study of these structures.

Let $(\mathcal{B}, \bullet_\mathcal{B}, \prs_\mathcal{B})$ be a pseudo-Euclidean Novikov superalgebra, and let $(\h , \cdot, \Omega)$ be  an associative supercommutative   superalgebra, endowed with an even  nondegenerate invariant symmetric bilinear form $\Omega$, i.e.,
\[
\Omega(a\cdot b, c) = \Omega(a, b\cdot c), \quad  \text{ for all } \quad a,b,c \in \h.
\]
Over the superspace $\A = \A_{\bar 0}\oplus \A_{\bar 1}$, where $\A_{\bar 0}=\mathcal{B}_{\bar 0}\otimes \h_{\bar 0} \oplus \mathcal{B}_{\bar 1}\otimes \h_{\bar 1}$  and $\A_{\bar 1}=\mathcal{B}_{\bar 0}\otimes \h_{\bar 1} \oplus \mathcal{B}_{\bar 1}\otimes \h_{\bar 0}$, we define  the following structures:
\begin{equation}\label{Om-g}
\begin{array}{ccl}
(u\otimes a)\bullet (v \otimes b)&:=& (-1)^{|a||v|}(u\bullet_\mathcal{B} v)\otimes (a \cdot b), \\[2mm]
\langle u \otimes a, v \otimes b\rangle &:=& (-1)^{|a||v|}\langle u, v\rangle_\mathcal{B} \, \Omega(a,b), 
\end{array}
\end{equation}
for all $u,v \in \mathcal{B}$ and $a,b \in \h$.   
\sssbegin{Proposition} The superalgebra defined by \textup{Eq.} \eqref{Om-g} is a pseudo-Euclidean Novikov superalgebra.
\end{Proposition}

\begin{proof}
It is easy to see that $(\mathcal A,\bullet)$ is a Novikov superalgebra. 
It is also clear that $\prs$ is a non-degenerate symmetric bilinear form. 
Let $u,v,w\in {\mathcal B}$ and $a,b,c\in \mathcal H$. We have
\[
\begin{aligned}
\langle (u\otimes a)\bullet (v\otimes b),\, w\otimes c\rangle
&= (-1)^{|a||v|}\,
\langle (u\bullet_{\mathcal B} v)\otimes (a\cdot b),\, w\otimes c\rangle \\
&= (-1)^{|a||v|+(|a|+|b|)|w|}
\langle u\bullet_{\mathcal B} v,\, w\rangle_{\mathcal B}\,
\Omega(a\cdot b,c) \\
&= -(-1)^{|a||v|+(|a|+|b|)|w|+|u||v|+|a||b|}
\langle v,\, u\bullet_{\mathcal B} w\rangle_{\mathcal B}\,
\Omega( b, a\cdot c) \\
&= -(-1)^{|a||v|+|a||w|+|u||v|+|a||b|+|b||u|}
\langle v\otimes b,\, (u\bullet_\mathcal{B} w)\otimes(a\cdot c)\rangle \\
&= -(-1)^{|a||v|+|u||v|+|a||b|+|b||u|}
\langle v\otimes b,\, (u \otimes a)\bullet (w \otimes c)\rangle
\\&= -(-1)^{|u \otimes a||v\otimes b|}
\langle v\otimes b,\, (u \otimes a)\bullet (w \otimes c)\rangle.
\end{aligned}
\]
Thus, $(\A, \bullet, \prs)$ is a pseudo-Euclidean Novikov superalgebra.
\end{proof}

The following lemma is the key result underlying the whole development of this work.
\sssbegin{Lemma}\label{Lemma11}
Let $(\mathcal A,\bullet)$ be a superalgebra endowed with a non-degenerate symmetric bilinear form $\prs$ such that the left multiplication operator $\Ll^\bullet$ is $\prs$-antisymmetric. If
\(
[\Rr_u^\bullet,\Rr_v^\bullet]=0,\)
 for all $u,v\in \mathcal A,$
then
\[
(u\bullet v)\bullet w = 0,
\quad \text{for all } u,v,w\in \mathcal A.
\]
Equivalently,
\[
\Ll_{u\bullet v}^\bullet = 0
\quad\text{for all } u,v\in \mathcal A,
\quad\text{or equivalently}\quad
\Rr_u^\bullet \circ \Rr_v^\bullet = 0
\quad\text{for all } u,v\in \mathcal A.
\]
\end{Lemma}
\begin{proof}
Assume that
\(
(u\bullet v)\bullet w = (-1)^{|v||w|}(u\bullet w)\bullet v,\)
for all  $u,v,w\in\A.$

Let $u,v,w,z\in\A$. Using the fact that the left multiplication operators
are $\prs$-antisymmetric, we obtain
\begin{align*}
\langle (u\bullet v)\bullet w, z\rangle &=-(-1)^{(|u|+|v|)|w|}\langle w, (u\bullet v)\bullet z\rangle\\ &= -(-1)^{(|u|+|v|)|w|+|v||z|}\langle w, (u\bullet z)\bullet v\rangle\\  &=   -(-1)^{(|u|+|v|)|w|+|v||z|+(|u|+|v|+|z|)|w|}\langle  (u\bullet z)\bullet v, w\rangle \\ &=    (-1)^{|v||z|+|z||w|+(|u|+|z|)|v|}\langle  v,  (u\bullet z)\bullet w\rangle \\ &=    (-1)^{|z||w|+|u||v|+|z||w|}\langle  v,  (u\bullet w)\bullet z\rangle \\ &= -   (-1)^{|u||v|+(|u|+|w|)|v|}\langle (u\bullet w)\bullet v,   z\rangle \\ &= -   (-1)^{|u||v|+(|u|+|w|)|v|+|w||v|}\langle (u\bullet v)\bullet w,   z\rangle \\ &= -   \langle (u\bullet v)\bullet w,   z\rangle.
\end{align*}
Since $\prs$ is non-degenerate, it follows that $(u\bullet v)\bullet w=0$, for all $u,v,w\in \A$.
\end{proof}
\sssbegin{Proposition}\label{Pr1}
Let $(\A,\bullet, \prs)$ be a pseudo-Euclidean superalgebra  such that the
left multiplication operator $\Ll^\bullet$ is  $\prs$-antisymmetric. Then the following assertions are
equivalent:
\begin{enumerate}
\item[$(i)$]
$(\A,\bullet)$ is a Novikov superalgebra.

\item[$(ii)$]
$(\A,\bullet)$ is an $\Ll\Rr$-superalgebra.

\item[$(iii)$]
$(\A,\bullet)$ is a left Leibniz $\Ll$-superalgebra.

\item[$(iv)$]
For all $u,v\in\A$,
\begin{equation}\label{kly0}
\Ll_{u\bullet v}^\bullet = 0
\quad \text{and} \quad
[\Ll_u^\bullet,\Ll_v^\bullet]=0.
\end{equation}
\end{enumerate}
\end{Proposition}
\begin{proof}
 It is an immediate consequence of Lemma \ref{Lemma11} and the relations \eqref{Nvleft}, \eqref{LR-eq} and \eqref{Leibniz}.
\end{proof} 

\sssbegin{Corollary}
Let $(\mathcal A, \bullet, \prs)$ be a pseudo-Euclidean Novikov superalgebra. Then $(\A, \bullet)$ is right-symmetric superalgebra if and only if $\A^-$ is a 2-step nilpotent Lie superalgebra.
\end{Corollary}
\begin{proof}
Let $u,v\in \A$. According to Proposition~\ref{Pr1}, we have
\[
\Ll_{u\bullet v} = \Rr_u \circ \Rr_v = \Rr_v \circ \Rr_u = 0.
\]
Thus, $(\A, \bullet)$ is a right-symmetric superalgebra if and only if 
$\Rr_{[u,v]}=0$.

On the other hand, for the Lie superalgebra $\A^-$, we have
\[
\ad_{[u,v]} = \Ll_{[u,v]} - \Rr_{[u,v]} = - \Rr_{[u,v]}.
\]

Therefore, $\A^-$ is $2$-step nilpotent if and only if $\Rr_{[u,v]}=0$, 
which is equivalent to $(\A, \bullet)$ being right-symmetric.
\end{proof}


\sssbegin{Proposition} \label{Nil}
Let $(\mathcal A, \bullet, \prs)$ be a pseudo-Euclidean Novikov superalgebra.
Then $(\mathcal A,\bullet)$ is nilpotent if and only if the associated Lie superalgebra
$\mathcal A^{-}$ is nilpotent.
Moreover, $(\mathcal A,\bullet)$ is nilpotent if and only if the left multiplication operator $\Ll_u^\bullet$ is nilpotent, for any $u\in\mathcal A$.
\end{Proposition}

\begin{proof}
Assume first that $(\mathcal A,\bullet)$ is nilpotent.
Then it is clear that the Lie superalgebra $\mathcal A^{-}$ is also nilpotent.

Conversely, assume that $\mathcal A^{-}$ is nilpotent.
Since $(\mathcal A,\bullet)$ is a Novikov superalgebra, Eq~\eqref{kly0} implies that $(u\bullet v)\bullet w=0$ for all $u,v,w\in \mathcal A$. It follows that 
\[
\Rr_u^\bullet \circ \Rr_v^\bullet
=0, \;  \Rr_u^\bullet \circ \Ll_v^\bullet
=0, \;  \Ll_{u\bullet v}^\bullet = 0,
\quad \text{for all } u,v\in\mathcal A.
\]
Hence, in order to prove that $(\mathcal A,\bullet)$ is nilpotent, it is sufficient to show that
the left multiplication operators are nilpotent, that is,
\[
\Ll_{u_1}^\bullet \circ \Ll_{u_2}^\bullet \circ \cdots \circ \Ll_{u_n}^\bullet = 0,
\quad \text{for all } u_i\in\mathcal A,
\]
for $n$ large enough.

Recall that the adjoint operator is given by
\[
\ad_u = \Ll_u^\bullet - \Rr_u^\bullet.
\]
A direct computation yields
\[
\ad_u \circ \ad_v
= (\Ll_u^\bullet - \Rr_u^\bullet)\circ (\Ll_v^\bullet - \Rr_v^\bullet)
= \Ll_u^\bullet \circ \Ll_v^\bullet - \Ll_u^\bullet \circ \Rr_v^\bullet,
\]
and
\[
\ad_u \circ \ad_v \circ \ad_w
= \Ll_u^\bullet \circ \Ll_v^\bullet \circ \Ll_w^\bullet
- \Ll_u^\bullet \circ \Ll_v^\bullet \circ \Rr_w^\bullet.
\]
By induction, for all $n\geq 1$ and for all $u_1,\dots,u_n\in\mathcal A$, we obtain
\begin{equation}\label{nil}
\ad_{u_1}\circ \cdots \circ \ad_{u_n}
= \Ll_{u_1}^\bullet \circ \cdots \circ \Ll_{u_n}^\bullet
- \Ll_{u_1}^\bullet \circ \cdots \circ \Ll_{u_{n-1}}^\bullet\circ \Rr_{u_n}^\bullet .
\end{equation}

Since $\mathcal A^{-}$ is nilpotent, there exists $n_0\in\mathbb N$ such that
\[
\ad_{u_1}\circ \cdots \circ \ad_{u_{n_0}} = 0,
\quad \text{for all } u_i\in\mathcal A.
\]
Using~\eqref{nil}, we deduce that
\[
\Ll_{u_1}^\bullet \circ \cdots \circ \Ll_{u_{n_0}}^\bullet
= \Ll_{u_1}^\bullet \circ \cdots \circ \Rr_{u_{n_0}}^\bullet.
\]
It follows, for all $u_{n_{0}+1}\in \A$, we have
\[
\Ll_{u_1}^\bullet \circ \cdots \circ \Ll_{u_{n_0}}^\bullet\circ \Rr_{u_{n_0+1}}^\bullet
= \Ll_{u_1}^\bullet \circ \cdots \circ \Rr_{u_{n_0}}^\bullet\circ \Rr_{u_{n_0+1}}^\bullet=0.
\]
From \eqref{nil}, we obtain that
\begin{equation}
\ad_{u_1}\circ \cdots \circ \ad_{u_{n_0}}\circ   \ad_{u_{n_0+1}}
= \Ll_{u_1}^\bullet \circ \cdots \circ \Ll_{u_{n_0}}^\bullet \circ   \Ll_{u_{n_0+1}}
- \Ll_{u_1}^\bullet \circ \cdots \circ \Ll_{u_{n_0}}^\bullet \circ   \Rr^\bullet_{u_{n_0+1}} .
\end{equation}
Thus, we deduce that $$\Ll_{u_1}^\bullet \circ \cdots \circ \Ll_{u_{n_0}}^\bullet \circ   \Ll^\bullet_{u_{n_0+1}}=0, \text{ for all $u_i\in \A$}
$$
which shows that $(\mathcal A,\bullet)$ is nilpotent.

Finally, if $\Ll_u^\bullet$ is nilpotent for all $u\in\mathcal A$, then~\eqref{nil} implies that
$\ad_u$ is nilpotent for all $u\in\mathcal A$.
By Engel’s theorem, the Lie superalgebra $\mathcal A^{-}$ is nilpotent, and therefore
$(\mathcal A,\bullet)$ is nilpotent as well.
The converse implication is obvious.
\end{proof}

Let $(\A,\bullet,\prs)$ be a pseudo-Euclidean Novikov
superalgebra. Denote by
\[
[\A,\A]:=\{[u,v]_\bullet \mid u,v\in\A\}
\]
the derived ideal of the associated Lie superalgebra
\(
\A^-,
\)
and by
\[
\A\bullet\A:=\{u\bullet v \mid u,v\in\A\}
\]
the two-sided ideal of the superalgebra $(\A,\bullet)$. Then we have
\begin{equation}\label{ortogonal}
(\A\bullet\A)^\perp = N_r(\A,\bullet),
\esp 
[\A,\A]^\bot
= \{u\in\A \mid (\Rr_u^\bullet)^*=\Rr_u^\bullet\}.
\end{equation}

\sssbegin{Proposition}\label{gg abelian}
	Let $(\A, \bullet,  \prs)$ be a pseudo-Euclidean Novikov superalgebra. Then 
	\begin{enumerate}
		\item[$(i)$]  The restriction of $\bullet$ to $\A\bullet\A$ is trivial; in particular, $[\A,\A]$ is an abelian ideal of $\A^-$.
		\item[$(ii)$]  $\left[\A, \A\right]$ is a graded two-sided ideal of $(\A,\bullet)$ and $\left[\A, \A\right]^{\perp}$  is a left ideal.
		\item[$(iii)$]  If $[\A, \A]$ is nondegenerate or $\A\bullet\A$ is nondegenerate then $ [\A, \A]=\A\bullet\A$ and  $ [\A, \A]^{\perp}$ is an abelian Lie subsuperalgebra of $\A^-$.
	\end{enumerate}

\end{Proposition}
\begin{proof}

Part (i). This is an immediate consequence of \eqref{kly0} and the fact that  $[\A, \A] \subset\A\bullet \A $. 

Part (ii). Let $u,v\in [\A,\A]$ and $w\in \A$.  By virtue of \eqref{kly0},  $[u,v]\bullet w=0$  and $w\bullet [u,v]=[w,[u,v]]$ which shows that $[\A,\A]$ is a graded two-sided ideal of $(\A,\bullet)$. On the other hand, for any $u\in\A$, $v\in [\A, \A]$ and $w\in [\A, \A]^\bot$, we have  
$$\langle u\bullet w, v\rangle=-(-1)^{|u||w|}\langle w, u\bullet v\rangle=0,$$
because $[\A, \A]$ is two-sided ideal of $(\A, \bullet)$. It follows that $u\bullet w\in[\A, \A]^\bot$  and hence $[\A, \A]^\bot$ is a left ideal of $(\A, \bullet)$. 
	
Part (iii). The inclusion $(\A\bullet\A)^\perp\subset [\A, \A]^\perp$ is obvious. 

From \eqref{ortogonal}, for any $u \in (\A\bullet \A)^\bot$ we have $\Rr^\bullet_u=0$ and therefore the restriction of $\bullet$ to $(\A\bullet\A)^\perp$ is trivial. On the other hand, for any $u\in[\A, \A]^\perp$, and by virtue of \eqref{kly0}, we have $\Rr^\bullet_u([\A,\A])=\Rr^\bullet_u(\A\bullet\A)=\{0\}$.  Since $\Rr^\bullet_u$ is symmetric, for any  $u\in[\A, \A]^\bot$, we get $\Rr^\bullet_u([\A,\A]^\perp)\subset [\A,\A]^\perp$ and $\Rr^\bullet_u((\A\bullet\A)^\perp)\subset (\A\bullet\A)^\perp$. Having these remarks in mind let us show our assertion.
		
Suppose that $[\A,\A]$ is nondegenerate, i.e, $\A=[\A,\A]\oplus[\A,\A]^\perp$.  

Since $\langle[u,v],w\rangle=0$, for any $u,v,w\in[\A,\A]^\perp$, and the fact that $\Rr^\bullet_v$ symmetric and $\Ll^\bullet_v$ is anti-symmetric, it follows that 
		\[ \langle u\bullet v,w\rangle =(-1)^{|u||v|}\langle v\bullet u, w\rangle=-\langle u,v\bullet w\rangle=-(-1)^{|v||w|}\langle u,w\bullet v\rangle
		=-\langle u\bullet v,w\rangle \]
		this implies that $\Rr^\bullet_u([\A,\A]^\perp)=0$ and therefore $\Rr^\bullet_u=0$. Thus $(\A\bullet\A)^\perp=[\A, \A]^\perp$ by \eqref{ortogonal}, which implies that  $\A\bullet\A=[\A, \A].$ The restriction of $\bullet$ to $(\A\bullet\A)^\perp$ is trivial and then $[\A, \A]^\perp$ is an abelian Lie subsuperalgebra.

		Suppose now that $\A\bullet\A$ is nondegenerate, i.e, $\A=\A\bullet\A\oplus(\A\bullet\A)^\perp$.   For any $v,w\in (\A\bullet\A)^\perp$ and $u\in[\A,\A]^\perp$, we have obviously $\langle v\bullet u,w\rangle=0$ and hence $\Rr^\bullet_u=0$. Arguing as in the first case, we obtain the desired result.
		\qedhere
   
\end{proof}

We now introduce an important class of pseudo-Euclidean Novikov superalgebras. In \cite[Theorem 1.5]{M}, Milnor showed that the Lie algebra $\mathfrak g$ of a flat Riemannian (or Euclidean) Lie group decomposes as a semidirect product of an abelian ideal $\mathfrak u$ and an abelian subalgebra $\mathcal B$. Furthermore, for any $u \in \mathcal B$, the adjoint map $\operatorname{ad}_u$ is antisymmetric. 

With respect to this  decomposition  $\mathfrak{g}=\mathfrak{u}\oplus\mathcal{B}$, the
Levi-Civita product $\bullet$ satisfies
$$
\Ll^\bullet_u =
\begin{cases}
0, & \text{if } u\in \mathfrak{u}, \\[2mm]
\ad_u, & \text{if } u\in \mathcal{B}.
\end{cases}
$$
Here $\ad_u=\Ll^\bullet_u-\Rr^\bullet_u$, for all $u\in\g$. This is
equivalent to the fact that $\g=\mathfrak{u}\oplus \mathfrak{u}^\bot$, with 
\[
\g\bullet\g \subseteq \mathfrak{u}
\quad \text{and} \quad
\Ll^\bullet_u = 0,
\ \text{for all } u\in\mathfrak{u}.
\]

We enlarge the context and introduce the following definition.
	
	\sssbegin{Definition}\label{defmilnor} A Milnor (super)algebra  is a pseudo-Euclidean (super)algebra $(\A,\bullet,\prs)$ such that $(\A,\bullet)$ is Lie-admissible superalgebra and  there exists  a \emph{nondegenerate} two-sided ideal  $I$ of $(\A,\bullet)$ such that $\A\bullet\A\subset I$ and  $\Ll^\bullet_u=0$, for any $u\in I$.
		
	\end{Definition}
     The following proposition justifies this terminology and demonstrates that Milnor  (super)algebras form a subclass of pseudo-Euclidean Novikov (super)algebras.

\sssbegin{Proposition}\label{milnor} Let $(\A,\bullet,\prs)$ be a Milnor superalgebra by means of a graded ideal $I$. Then:
		\begin{enumerate}
			\item[$(i)$] Both $I$ and $I^\perp$ are  trivial subsuperalgebras of $(\A,\bullet)$. Thus $I$ is an abelian graded ideal of $\A^-$ and $I^\perp$ is an abelian Lie subsuperalgebra. 
			\item[$(ii)$] $\A=I\oplus I^\perp$ and, for any $u\in \A$,
			\[ \Ll^\bullet_u=\begin{cases}
				0& \mbox{if}\quad u\in I,\\
				\ad_u & \mbox{if}\quad u\in I^\perp.
			\end{cases} \]
		\item[$(iii)$]  The center $Z(\A^-)$ of $\A^-$ satisfies $Z(\A^-)=\{u\in\A, \, \Ll^\bullet_u=\Rr^\bullet_u=0\}$.
		
		\end{enumerate} 
	In particular, $(\A,\bullet,\prs)$ is a pseudo-Euclidean Novikov superalgebra.
		\end{Proposition}

\begin{proof} Part (i).  It is a consequence of the fact that $\Ll^\bullet_u=0$, for any $u\in I$, and the fact that $I^\perp\subset(\A\bullet\A)^\perp$ which implies $\Rr^\bullet_u=0$, for any $u\in I^\perp$,  by virtue of \eqref{ortogonal}.

Part (ii). The only thing to point out is that $I^\perp\subset(\A\bullet\A)^\perp$ and so, for any $u\in I^\perp$, $\Rr^\bullet_u=0$ and $\ad_u=\Ll^\bullet_u$. 

Part (iii). Let $v=v_1+v_2\in Z(\A^-)$ where 
		$v_1\in I$ and $v_2\in I^\perp$. Since both $I$ and $I^\perp$ are abelian Lie subsuperalgebras, we deduce that $v_1,v_2\in Z(\A^-)$. 
From $0 = \ad_{v_1} = \Ll^\bullet_{v_1} - \Rr^\bullet_{v_1}$ and 
$0 = \ad_{v_2} = \Ll^\bullet_{v_2} - \Rr^\bullet_{v_2}$, we deduce that 
$\Ll^\bullet_{v_1} = \Rr^\bullet_{v_1}$ and $\Ll^\bullet_{v_2} = \Rr^\bullet_{v_2}$. 

Moreover, since $v_1 \in I$, we have 
$\Ll^\bullet_{v_1} = \Rr^\bullet_{v_1} = 0$, and since $v_2 \in I^\perp$, we also have 
$\Ll^\bullet_{v_2} = \Rr^\bullet_{v_2} = 0$. 

Therefore, we obtain
\[
\Ll^\bullet_v = \Ll^\bullet_{v_1} + \Ll^\bullet_{v_2} = 0, \qquad
\Rr^\bullet_v = \Rr^\bullet_{v_1} + \Rr^\bullet_{v_2} = 0.
\]

The fact that $(\A,\bullet,\prs)$ is pseudo-Euclidean Novikov superalgebra is an immediate consequence of $(ii)$ and the fact that $\br_\bullet$ is a Lie bracket.  Indeed, since 
$\A \bullet \A \subseteq I$, we have $\Ll^\bullet_{u\bullet v} =0$, for all $u,v \in \A$, and 
$[\Ll^\bullet_u, \Ll^\bullet_w] \overset{(ii)}{=} 0$ for all $u \in I$ and $w \in \A$. Moreover, since $I^\bot$ is an abelian Lie subsuperalgebra,   for $w,z \in I^\perp$, we have 
\[
0 = \ad_{[w,z]} = [\ad_w,\ad_z] = [\Ll^\bullet_w-\Rr^\bullet_w, \Ll^\bullet_z-\Rr^\bullet_z] \overset{(ii)}{=} [\Ll^\bullet_w, \Ll^\bullet_z].
\]
Hence, $(\A, \bullet)$ satisfies the identities \eqref{kly0}, and thus 
$(\A, \bullet)$ is a Novikov superalgebra.
	\end{proof}

\sssbegin{Theorem} \label{ ideal non degenere}
	Let $(\A, \bullet, \prs)$ be a  pseudo-Euclidean Novikov superalgebra such that $\A\bullet\A$ is  nondegenerate. Then $(\A,\bullet,\prs)$ is a Milnor  superalgebra by means of $\A\bullet\A=[\A,\A]$ and the conclusions of Proposition \ref{milnor} hold.
\end{Theorem}
\begin{proof} 
Since $\A\bullet\A$ is nondegenerate, by Proposition~\ref{gg abelian} we have $
\A\bullet\A=[\A,\A].
$ Then, by Proposition~\ref{Pr1}, we obtain $\Ll^\bullet_u=0$ for all $u\in \A\bullet\A$. 
Hence $(\A,\bullet,\prs)$ is a Milnor superalgebra with 
$\A\bullet\A=[\A,\A]$, and the conclusions of Proposition~\ref{milnor} hold.
\end{proof}

\sssbegin{Corollary}\label{nilpotent}
Let $(\A, \bullet, \prs)$ be a   pseudo-Euclidean Novikov superalgebra such that $\A\bullet\A$ is non-degenerate. Then, $\A^-$ is nilpotent if and only if  $\bullet$ is trivial. 
\end{Corollary}

\begin{proof}
Assume that $(\A, \bullet)$ is non-trivial. According to Theorem~\ref{  ideal non degenere}, 
$(\A,\bullet,\prs)$ is a Milnor pseudo-Euclidean superalgebra by means of  $\A\bullet\A=[\A,\A]$. Since $\A^-$ is nilpotent, we have
$
Z(\A^-)\cap [\A,\A]\neq 0.$ By part (iii) of Proposition~\ref{milnor} and Eq.~\eqref{ortogonal}, we obtain
$
Z(\A^-)\subseteq [\A,\A]^{\bot}.
$
Hence \[
\{0\}\neq Z(\A^-)\cap [\A,\A]\subseteq [\A,\A]^{\bot}\cap [\A,\A],
\]
which is a contradiction since $[\A,\A]$ is nondegenerate. 

Conversely, the statement is obvious.
\end{proof}

\subsection{Representations of pseudo-Euclidean Novikov superalgebras}
Here, we show, as a consequence of the results obtained so far, that every pseudo-Euclidean Novikov superalgebra naturally carries the structure of a left Leibniz $\Ll$-superalgebra. We then introduce the notion of representation of a left Leibniz $\Ll$-superalgebra and define the associated product related to pseudo-Euclidean Novikov superalgebras. Furthermore, we establish a characterization of pseudo-Euclidean Novikov superalgebras in terms of bimodule isomorphisms. More precisely, we show that the existence of a suitable bilinear form is equivalent to the isomorphism between certain naturally associated $\A$-bimodules. Motivated by this representation, we construct a new class of pseudo-Euclidean Novikov superalgebras, analogous to the $T^*$-extension for quadratic algebras introduced in \cite{Bordemann}. 
\sssbegin{Definition}\label{Def3.1}
Let $(\A,\cdot)$ be a non-associative superalgebra satisfying some polynomial identities $(P_1),\dots,(P_n)$. Let $V$ be a vector superspace, and let 
$
r,l : \A \longrightarrow \mathrm{End}(V)
$
be two linear maps. On the direct sum $\A \oplus V$, define a product $\ast$ by
\begin{equation}\label{eq7}
(u+x)\ast (v+y) := u \cdot v + l(u)(y) + (-1)^{|u||v|}r(v)(x), \text{ for all $u,v \in \A$ and $x,y  \in V$.}
\end{equation}

We say that $(V,r,l)$ is an $\A$-bimodule if $(\A \oplus V,\ast)$ is a superalgebra satisfying the polynomial identities $(P_1),\dots,(P_n)$. In this case, the pair $(r,l)$ is called a representation of $\A$ in $V$ relative to the identities  $(P_1),\dots,(P_n)$.
\end{Definition}

Based on the above definition, we now specify the notion of representation in the case of left-Leibniz $\Ll$-superalgebras.

\sssbegin{Definition}\label{Def3.2}
Let $(\A, \bullet)$ be a left-Leibniz $\Ll$-superalgebra, $V$ be a vector superspace, and  
$
r,l : \A \longrightarrow \mathrm{End}(V)
$
be two linear maps. The pair $(r,l)$ is called a representation of $\A$ in $V$ if $(\A \oplus V,\ast)$, where $\ast$ is defined by Eq. \eqref{eq7}, is a left-Leibniz $\Ll$-superalgebra.
\end{Definition}
\sssbegin{Proposition}
Let $(\A, \bullet)$ be a left-symmetric superalgebra, $V$ be a vector superspace, and  
$
r,l : \A \longrightarrow \mathrm{End}(V)
$
be two linear maps. Then $(r,l)$ is a representation of $\A$ in $V$ if and only if, for all $u,v \in \A$, the following identities hold: 
\begin{equation}\label{eq3.13}
\begin{aligned}
\relax l([u, v]_\bullet)&= [l(u),l(v)], \quad r(u\bullet v)-(-1)^{|u||v|}r(v)\circ r(u)=[l(u),r(v)] 
\end{aligned}
\end{equation}

\end{Proposition}

\sssbegin{Proposition}\label{Prop3.3}
Let $(\A, \bullet)$ be a left-Leibniz $\Ll$-superalgebra, $V$ be a vector superspace, and  
$
r,l : \A \longrightarrow \mathrm{End}(V)
$
be two linear maps. Then $(r,l)$ is a representation of $\A$ in $V$ if and only if, for all $u,v \in \A$, the following identities hold:
\begin{equation}\label{eq3.13b}
\begin{aligned}
\relax [l(u),l(v)] &= 0, & 
l(u)r(v) &= r(u \bullet v), & 
l(u \bullet v) &= 0, & 
r(u)r(v) &= 0, & 
r(u)l(v) &= 0.
\end{aligned}
\end{equation}
\end{Proposition}

\sssbegin{Proposition}\label{Prop3.4}
Let $(\A, \bullet)$ be a left-Leibniz $\Ll$-superalgebra, and let 
$
\Ll^\bullet,\Rr^\bullet : \A \longrightarrow \mathrm{End}(\A)
$
be the left and right multiplication operators, respectively. Then $(\Rr^\bullet, \Ll^\bullet)$ is a representation of $\A$ in $\A$, called the \emph{adjoint representation} \textup{(}or \emph{regular representation}\textup{)} of $\A$.
\end{Proposition}

\begin{proof}
The result follows directly from Definition~\ref{Def3.2} and Proposition~\ref{Prop3.3}.
\end{proof}
\sssbegin{Remark}
For a $(\A, \bullet, \prs)$ a pseudo-Euclidean Novikov superalgebra, Proposition~\ref{Pr1} implies that $(\A, \bullet)$ is a left Leibniz $\Ll$-superalgebra. Consequently, any representation of a pseudo-Euclidean Novikov superalgebra naturally induces a representation of the associated left Leibniz $\Ll$-superalgebra.
\end{Remark}

\sssbegin{Remark}
Let $(\A, \bullet, \prs)$ be an even pseudo-Euclidean Novikov superalgebra. Then:

\begin{itemize}
\item[(i)] The even part $(\A_{\bar{0}}, \bullet, \prs|_{\A_{\bar{0}}})$ is a pseudo-Euclidean Novikov algebra.

\item[(ii)] The odd part $(\A_{\bar{1}}, \prs|_{\A_{\bar{1}}})$ is a symplectic vector space.

\item[(iii)] The maps (for all  $u\in \A_{\bar 0}$ and $x\in \A_{\bar 1}$)
\[
r : \A_{\bar{0}} \to \End(\A_{\bar{1}}), \quad r(u)(x) = \Rr_u^\bullet(x),
\]
and
\[
l : \A_{\bar{0}} \to \End(\A_{\bar{1}}), \quad l(u)(x) = \Ll_u^\bullet(x),
\]
define a representation of $\A_{\bar{0}}$ on $\A_{\bar{1}}$, where each $l(u)$ is $\prs|_{\A_{\bar{1}}}$-antisymmetric.

\item[(iv)] Moreover, for all $x,y \in \A_{\bar{1}}$ and $u \in \A_{\bar{0}}$, we have
\[
\langle x \bullet y, u \rangle_{\A_{\bar{0}}}
= \langle y, x \bullet u \rangle_{\A_{\bar{1}}}
= \langle y, r(u)(x) \rangle_{\A_{\bar{1}}}.
\]
\end{itemize}

Conversely, let $(\mathcal{B}, \bullet_{\mathcal{B}}, \prs_{\mathcal{B}})$ be a pseudo-Euclidean Novikov algebra, $(V, \prs_V)$ be a symplectic vector space, and let $(l,r)$ be a representation of $\mathcal{B}$ on $V$ such that $l(u)$ is $\prs_V$-antisymmetric for all $u \in \mathcal{B}$. 

Define $\A := \mathcal{B} \oplus V$ with $\A_{\bar{0}} = \mathcal{B}$ and $\A_{\bar{1}} = V$. The product $\bullet$ and the bilinear form $\prs$ are given by
\[
(u+x)\bullet (v+y)
= u \bullet_{\mathcal{B}} v + l(u)(y) + r(v)(x) + x \circ y, \; \text{for all $u,v \in \mathcal{B}$ and $x,y \in V$,}
\]
and
\[
\prs = \prs_{\mathcal{B}} + \prs_V,
\]
where the product on $V$ is defined by
\[
\langle x \circ y, u \rangle_{\mathcal{B}} = \langle y, r(u)(x) \rangle_V, \text{ for all $u\in \mathcal{B}$ and $x,y \in V$.}
\]

Assume that the following conditions hold:
\begin{align*}
\langle x \circ t,\, y \circ z \rangle_{\A_{\bar{0}}}
&= - \langle y \circ t,\, x \circ z \rangle_{\A_{\bar{0}}},\;\;
\langle r(u)(x),\, r(v)(y) \rangle_{V}
= - \langle r(u)(y),\, r(v)(x) \rangle_{V}, \\
l(x \circ y) &= 0, \esp
(x \circ y)\bullet_{\mathcal{B}} u = 0,
\end{align*}
for all $x,y,z,t \in V$ and $u \in \mathcal{B}$.

Then $(\A, \bullet, \prs)$ is an even pseudo-Euclidean Novikov superalgebra.

\medskip

The above construction is analogous to that introduced for even quadratic Lie superalgebras (see, for instance, \cite{Ko}).  Indeed, in the Lie case, the invariance of the bilinear form is expressed by
\[
\langle [u,v], w \rangle  = -(-1)^{|u||v|} \langle v, [u,w] \rangle.
\]
\end{Remark}

Let $(\A, \bullet, \prs)$ be a pseudo-Euclidean left-symmetric  superalgebra. Consider the bilinear map $\star: \A \times \A \longrightarrow \A$ defined by
\begin{equation}
\langle u\bullet v, w\rangle=\langle u, v \star w\rangle, \quad \text{ for all } u,v,w \in \A .
\label{produce associe}\end{equation}
\sssbegin{Proposition}\label{Pr associe}
     The product $\star$ defined in Eq. \ref{produce associe} satisfies the following identities:
\begin{enumerate}
\item[$(a)$]
$u \star v=-(-1)^{|u||v|}v \star u$,\; $\text{ for all } u,v\in\A$,
\item[$(b)$] $(u\bullet v)\star w-u\star(v\star w)=u\bullet(v\star w)-(-1)^{|u||v|} v\star (u\bullet w),$\; $\text{ for all } u,v, w\in\A$
\end{enumerate}

Moreover, if $(\A, \bullet, \prs)$ is a pseudo-Euclidean Novikov superalgebra, then the  product $\star$ defined in Eq.~\ref{Pr associe} satisfies:
     \begin{enumerate}
\item[$(c)$] $u \star(v \star w)=0$, \;$\text{ for all } u,v,w\in\A$,
\item[$(d)$]$u \bullet(v \star w)=0$, \;$\text{ for all } u,v,w\in\A$,
\item[$(e)$]  $(u\bullet v)\star w=-(-1)^{|u||v|}v\star (u \bullet w), \text{ for all } u,v,w \in \A$.
\end{enumerate}
\end{Proposition}

\begin{proof}
Let $u,v,w,z \in \A$.
     \begin{enumerate}
         \item[$(a)$] We have\begin{align*}
      \langle u, v\star w\rangle=&      \langle u \bullet v, w\rangle= -(-1)^{|u||v|}\langle v, u\bullet w\rangle= -(-1)^{|v||w|}\langle u\bullet w, v\rangle=-(-1)^{|v||w|}\langle u, w\star v\rangle.
            \end{align*}
Thus, the product $\star$ is antisymmetric.
\item[$(b)$]  Since $(\A, \bullet)$ is left-symmetric, then we have 
\begin{align*}
0=&\langle  (u\bullet v)\bullet w-(-1)^{|u||v|}   (v\bullet u)\bullet w-u\bullet (v\bullet w)+(-1)^{|u||v| }v\bullet (u\bullet w), z\rangle\\&= \langle u, v\star (w\star z)\rangle +\langle u, v\bullet (w\star z)\rangle -\langle u, (v\bullet w)\star z\rangle -(-1)^{|v||w|}\langle u, w\star (v\bullet z)\rangle\\&=\langle u, v\star (w\star z)+ v\bullet (w\star z)-(v\bullet w)\star z-(-1)^{|v||w|} w\star (v\bullet z)\rangle
\end{align*}
Since $\prs$ is nondegenerate, then \[
v\star (w\star z)+ v\bullet (w\star z)=(v\bullet w)\star z+(-1)^{|v||w|} w\star (v\bullet z).
\]
\end{enumerate}
Now assume that $(\A, \bullet, \prs)$ is a pseudo-Euclidean Novikov superalgebra.
 \begin{enumerate}

\item[$(c)$] According to Proposition \ref{Pr1}, we have
\begin{align*}
0\overset{\eqref{kly0}}{=}\langle (u\bullet v)\bullet w, z\rangle=\langle u\bullet v, w\star z\rangle= \langle u , v\star( w\star z)\rangle.
\end{align*}
Since $\prs$ is nondegenerate, then  $v\star( w\star z)=0$.
\item[$(d)$] According to Proposition \ref{Pr1}, we have
\begin{align*}
0\overset{\eqref{kly0}}{=}\langle (u\bullet v)\bullet w, z\rangle=\langle u\bullet v, w\star z\rangle= -(-1)^{|u||v|}\langle v , u\bullet( w\star z)\rangle.
\end{align*}
Since $\prs$ is nondegenerate, then  $u\bullet( w\star z)=0$.
\item[$(e)$] This part follows from (b), (c) and (d).\end{enumerate}\end{proof}


\sssbegin{Proposition}
Let $(\A, \bullet, \prs)$ be a pseudo-Euclidean Novikov superalgebra. Then the associated Lie superalgebra $(\A, \star)$ is $2$-step nilpotent.

Conversely, let $(\g, \star, \prs)$ be a $2$-step nilpotent Lie superalgebra endowed with a non-degenerate symmetric bilinear form. Define a bilinear product $\bullet$ on $\g$ by
\[
\langle u \bullet v, w \rangle = \langle u, v \star w \rangle,
\qquad \text{ for all } u,v,w \in \g.
\]
If identity $(e)$ of Proposition~\ref{Pr associe} is satisfied, then $(\g, \bullet, \prs)$ is a pseudo-Euclidean Novikov superalgebra.
\end{Proposition}

\begin{proof}
The first assertion follows directly from parts (a) and (c) of Proposition~\ref{Pr associe}. 

Conversely, assuming that $({\mathfrak g}, \star)$ is a $2$-step nilpotent Lie superalgebra and that condition $(e)$ holds, one verifies that the Novikov identities are satisfied. The proof is straightforward and follows the same arguments as in Proposition~\ref{Pr associe}.
\end{proof}

\sssbegin{Remark}
It is worth noting that similar constructions have been studied in the non-super setting. Let $\g$ be a $2$-step nilpotent Lie algebra endowed with  a non-
degenerate symmetric bilinear form $\prs$.

In this context, one considers the bilinear map $J : \g \times \g \to \g$ defined by
\[
\langle J_u(v), w \rangle = \langle u, [v,w] \rangle,
\text{ for all } u,v,w \in \g.
\]
This product has been extensively studied from both geometric and algebraic viewpoints (see, for instance, \cite{ AFMV, CP, Patrick1, Patrick2, JPP,  WG}).
\end{Remark}

\sssbegin{Remark}
Let $(\A, \bullet, \prs)$ be a pseudo-Euclidean Novikov superalgebra. As mentioned previously, $(\A, \bullet)$ is not necessarily nilpotent; indeed, there exist non-nilpotent examples (see Remark~\ref{TT}). Consequently, the associated Lie superalgebra $(\A, \br_\bullet)$ is not necessarily nilpotent.

On the other hand, the Lie superalgebra $(\A, \star)$ is always $2$-step nilpotent. Therefore, $\A$ carries two Lie superalgebra structures: one induced by the commutator $\br_\bullet$ and the other one is given by $\star$.

\end{Remark}

\sssbegin{Proposition}\label{prop-cyclicity}
The relation $\star = -[\cdot, \cdot]_\bullet$ holds between the Lie superalgebra structures $\mathcal{A}^-$ and $(\mathcal{A}, \star)$ if and only if $(\mathcal{A}^-, \prs)$ satisfies the following cyclicity equation:
\begin{equation}\label{cyclicity}
(-1)^{|u||w|}\langle [u,v]_\bullet, w\rangle + (-1)^{|u||v|}\langle [v,w]_\bullet, u\rangle + (-1)^{|v||w|}\langle [w,u]_\bullet, v\rangle = 0
\end{equation}
for all $u, v, w \in \mathcal{A}$. In this case, the following identity holds:
\[
\langle u \bullet v, w\rangle  = -\langle u, [v, w]_\bullet\rangle.
\]
\end{Proposition}

\begin{proof}
Assume that $\star = -[\cdot, \cdot]_\bullet$. Let $u, v, w \in \mathcal{A}$. According to \eqref{produce associe}, we have:
\begin{align*}
\langle u \bullet v, w\rangle &= \langle u, v \star w\rangle 
= -\langle u, [v,w]_\bullet\rangle \\
&= -\langle u, v \bullet w\rangle + (-1)^{|v||w|} \langle u, w \bullet v\rangle \\
&= -(-1)^{|u||v|+|u||w|}\langle v \bullet w, u\rangle - (-1)^{|v||w|+|u||w|} \langle w \bullet u, v\rangle.
\end{align*}
This implies that
\[
(-1)^{|u||w|}\langle u \bullet v, w\rangle + (-1)^{|u||v|}\langle v \bullet w, u\rangle + (-1)^{|v||w|}\langle w \bullet u, v\rangle = 0.
\]
Since $[u,v]_\bullet = u \bullet v - (-1)^{|u||v|} v \bullet u$, a straightforward computation shows that:
\[
(-1)^{|u||w|}\langle [u,v]_\bullet, w\rangle + (-1)^{|u||v|}\langle [v,w]_\bullet, u\rangle + (-1)^{|v||w|}\langle [w, u]_\bullet, v\rangle = 0.
\]
Hence, $(\mathcal{A}^-, \prs)$ satisfies Equation~\eqref{cyclicity}.

\medskip

Conversely, assume that $(\mathcal{A}^-, \prs)$ satisfies Equation~\eqref{cyclicity}. Let $u, v, w \in \mathcal{A}$. We have: 
\begin{align*}
0 &= (-1)^{|u||w|}\langle [u,v]_\bullet, w\rangle + (-1)^{|u||v|}\langle [v,w]_\bullet, u\rangle + (-1)^{|v||w|}\langle [w, u]_\bullet, v\rangle \\
  &= (-1)^{|u||w|} \langle u \bullet v, w \rangle - (-1)^{|u||w|+|u||v|} \langle v \bullet u, w\rangle + (-1)^{|u||v|}\langle [v, w]_\bullet, u\rangle \\
  &\quad + (-1)^{|v||w|} \langle w \bullet u, v\rangle - (-1)^{|v||w|+|u||w|} \langle u \bullet w, v\rangle \\
  &= (-1)^{|u||w|} \langle u \bullet v, w \rangle + (-1)^{|u||w|} \langle u, v \bullet w\rangle + (-1)^{|u||v|}\langle [v, w]_\bullet, u\rangle \\
  &\quad - (-1)^{|v||w|+|u||w|} \langle u, w \bullet v\rangle + (-1)^{|v||w|} \langle w, u \bullet v\rangle \\
  &= 2(-1)^{|u||w|} \langle u \bullet v, w \rangle + (-1)^{|u||v|}2\langle [v, w]_\bullet, u\rangle \\
  &= 2(-1)^{|u||w|} \langle u \bullet v, w \rangle + (-1)^{|u||w|}2\langle u, [v, w]_\bullet\rangle.
\end{align*}
Therefore, 
\[
\langle u \bullet v, w \rangle = -\langle u, [v,w]_\bullet \rangle,
\]
or equivalently, $\langle u \bullet v, w\rangle = -\langle u, [v, w]_\bullet\rangle$, which proves that $\star = -[\cdot, \cdot]_\bullet$.
\end{proof}

Let $(\A, \cdot)$ be a superalgebra. For a linear map $\phi : \A \to \mathrm{End}(V)$, we define a linear map $\phi^* : \A \to \mathrm{End}(V^*)$ by 
\[
 \phi^*(u)(f)(x)  := - (-1)^{|f|(|\phi|+|u|)} f(\phi(u)(x)) 
\quad \text{ for all } u \in \A,\; f \in V^*,\; x\in V.
\]
\sssbegin{Proposition}\label{rep1}
Let $(\A, \bullet,\prs)$ be a pseudo-Euclidean Novikov superalgebra,  then the pair $
\big( -(\Ll^\star)^*, (\Ll^\bullet)^* \big)
$
defines a representation of $\A$ on $\A^*$. 
\end{Proposition}

\begin{proof}
Assume that $(\A, \bullet)$ is a Novikov superalgebra. Let $u,v \in \A$ and $f \in \A^*$. By the Novikov identities  (c), (d), and (e) of Proposition~\ref{Pr associe}, we obtain
\[
(\Ll_u^\bullet)^* \circ (\Ll_v^\bullet)^*(f)
= -(-1)^{|f||v|} (\Ll_u^\bullet)^*(f \circ \Ll_v^\bullet)
= (-1)^{|f||v| + |f||u| + |u||v|} f \circ \Ll_v^\bullet \circ \Ll_u^\bullet.
\]
Using again the left Leibniz $\Ll$-superalgebra identities, this yields
\[
(\Ll_u^\bullet)^* \circ (\Ll_v^\bullet)^*
= (-1)^{|u||v|} (\Ll_v^\bullet)^* \circ (\Ll_u^\bullet)^*.
\]
Moreover,
\[
(\Ll_{u \bullet v}^\bullet)^*(f)
= (-1)^{|f|(|u|+|v|)} f \circ \Ll_{u \bullet v}^\bullet = 0,
\]
by the left Leibniz $\Ll$-superalgebra condition.

Similarly, we have
\[
(\Ll_u^\star)^* \circ (\Ll_v^\star)^*(f)
= (-1)^{|f||v| + |f||u| + |u||v|} f \circ \Ll_v^\star \circ \Ll_u^\star = 0,
\]
by identity $(c)$, and
\[(\Ll_u^\star)^* \circ (\Ll_v^\bullet)^*(f)
= (-1)^{|f||v| + |f||u| + |u||v|} f \circ \Ll_v^\bullet \circ \Ll_u^\star = 0,\]
by identity $(d)$. Moreover, we have
\[
-(\Ll_u^\bullet)^* \circ (\Ll_v^\star)^*(f)
= -(-1)^{|f||v| + |f||u| + |u||v|} f \circ \Ll_v^\star \circ \Ll_u^\bullet
= (-1)^{|f||v| + |f||u|} f \circ \Ll_{u \bullet v}^\star
= -(\Ll_{u \bullet v}^\star)^*(f),
\]
by identity $(e)$.

Therefore, $\big(- (\Ll^\star)^*, (\Ll^\bullet)^* \big)$ defines a representation of $\A$ on $\A^*$.
\end{proof}
\sssbegin{Remark}\label{Re-left}
If $(\A, \bullet,\prs)$ is a pseudo-Euclidean left-symmetric superalgebra, then the pair $
\big( -(\Ll^\star)^*, (\Ll^\bullet)^* \big)
$
defines a representation of $\A$ on $\A^*$. 
\end{Remark}
\sssbegin{Proposition}\label{Prop7.8}
Let $(\A,\bullet,\prs)$ be a pseudo-Euclidean Novikov superalgebra. Let $\Ll^\bullet$ and $\Rr^\bullet$ denote the left and right multiplication operators, respectively, and  let $\Ll^\star$ denote the left multiplication operator associated with the product $\star$ defined in Eq.~\eqref{produce associe}. Then there exists a homogeneous  isomorphism $
\Phi : \A \longrightarrow \A^{*}
$
such that
\[
\Phi(\Ll_u^\bullet(a)) = (-1)^{|\Phi||u |}(\Ll_u^\bullet)^{*}(\Phi(a)) \; \text{and} \; \Phi(\Rr^\bullet_u(a)) = -(-1)^{|\Phi||u |}(\Ll_u^\star)^{*}(\Phi(a)), \text{ for all $u,a \in \A$.}
\]
\end{Proposition}

\begin{proof}
Consider the linear map $\Phi : \A \longrightarrow \A^{*}$ defined by
$
\Phi(a) := \langle a, \cdot \rangle, $ for all $a \in \A.$
Since $\prs$ is non-degenerate, $\Phi$ is an isomorphism of vector superspaces. 

Let $u,a,b \in \A$. Using the $\prs$-antisymmetry of the left multiplication operators, we obtain
\[
\Phi(\Ll^\bullet_u(a))(b) = \langle u \bullet a, b \rangle
= -(-1)^{|u||a|} \langle a, u \bullet b \rangle
= -(-1)^{|u||a|} \Phi(a)(\Ll^\bullet_u(b)).
\]
Thus, 
\[
\Phi(\Ll^\bullet_u(a)) =(-1)^{|\Phi||u |} (\Ll_u^\bullet)^{*}(\Phi(a)).
\]
Similarly, for the right multiplication operator, we have 
\[
\Phi(\Rr^\bullet_u(a))(b) = (-1)^{|u||a|}\langle a \bullet u, b \rangle=(-1)^{|u||a|}\langle a , u\star b \rangle= (-1)^{|u||a|}
\Phi(a)(\Ll_u^\star(b)) 
\]
which implies
\[
\Phi(\Rr^\bullet_u(a))(b)  =- (-1)^{|\Phi||u |}(\Ll_u^\star)^{*}(\Phi(a)).
\]
This completes the proof.
\end{proof}
Now we will prove the converse of Proposition \ref{Prop7.8}. We will divide the proof into two cases depending on the parity of the bilinear form.

\sssbegin{Theorem}\label{Thm7.11}
Let $(\A,\bullet)$ be a Novikov superalgebra such that $\dim \A_{\bar{1}}$ is even,  and let $\Ll^\bullet$ and $\Rr^\bullet$ denote the left and right multiplication operators, respectively. Then there exists an even non-degenerate symmetric bilinear form 
$
\prs : \A \times \A \longrightarrow \mathbb{K}
$
such that $(\A,\bullet, \prs)$ is an even pseudo-Euclidean Novikov superalgebra if and only if there exists an even isomorphism 
$
\Phi : \A \longrightarrow \A^{*}
$
satisfying \[
\Phi(\Ll_u^\bullet(a)) = (\Ll_u^\bullet)^{*}(\Phi(a)) \; \text{ and } \; \Phi(\Rr^\bullet_u(a)) = -(\Ll_u^\star)^{*}(\Phi(a)), 
\]
where the product $\star$  is antisymmetric. 
\end{Theorem}
\begin{proof}
Assume first that $(\A,\bullet, \prs)$ is an even pseudo-Euclidean Novikov superalgebra. Then the result follows directly from Proposition~\ref{Prop7.8}.

Conversely, suppose that there exists an even isomorphism $\Phi : \A \to \A^{*}$ such that
\[
\Phi(\Ll_u^\bullet(a)) = (\Ll_u^\bullet)^{*}(\Phi(a)) \; \text{and} \; \Phi(\Rr^\bullet_u(a)) = -(\Ll_u^\star)^{*}(\Phi(a)), \text{ for all $u,a \in \A$.}
\]
Consider the bilinear form
$
T : \A \times \A \longrightarrow \mathbb{K}$ defined by $T(u,v) := \Phi(u)(v),$
for all $u,v \in \A$. Since $\Phi$ is an even isomorphism, $T$ is an even non-degenerate bilinear form.

Moreover, for all  $u,v,w \in \A$, we have
\[
\begin{aligned}
T(u \bullet v, w) &= \Phi(u \bullet v)(w) 
= \Phi(\Ll_u^\bullet(v))(w) \\&= (\Ll^\bullet_u)^*(\Phi(v))(w)=-\Phi(v)(\Ll^\bullet_u(w))=-(-1)^{|u||v|} T(v, u\bullet w)
\end{aligned}
\]
Thus, we get that 
\begin{equation}\label{s1}
T(u \bullet v, w) =-(-1)^{|u||v|} T(v, u\bullet w).
\end{equation}
Moreover, we have 
\[
\begin{aligned}
T(v \bullet u, w) &= \Phi(v \bullet u)(w) 
=(-1)^{|u||v|} \Phi(\Rr_u^\bullet(v))(w) = -(-1)^{|u||v|}(\Ll^\star_u)^*(\Phi(v))(w)\\& =\Phi(v)(\Ll^\star_u(w))=T(v, u\star w)=-(-1)^{|u||w|}T(v, w\star u)\\
&=-(-1)^{|u||w|} T(v\bullet w, u) \overset{\eqref{s1}}{=}(-1)^{(|u|+|v|)|w|}T(w, v\bullet u)
\end{aligned}
\]
Thus, we obtain  
\begin{equation}\label{s2}
T(u \bullet v, w) = (-1)^{(|u|+|v|)|w|}T(w, u\bullet v).
\end{equation}

Now, we introduce the  even symmetric and anti-symmetric parts of the bilinear form $T$. Define the bilinear forms $T_s, T_a : \A \times \A \to \mathbb{K}$ by
\[
T_s(u,v) := \frac{1}{2}\Big(T( u,v ) + (-1)^{|u||v|}T (v,u)\Big),\esp
T_a(u,v) := \frac{1}{2}\Big(T( u,v ) - (-1)^{|u||v|}T( v,u )\Big),
\]
for all  $u,v \in \A$. Clearly, $T = T_s + T_a$.

Since $T$ satisfies \eqref{s1}, for all  $u,v,w \in \A$, we obtain
\[
\begin{aligned}
T_s(u \bullet v, w) &= \frac{1}{2}\Big(T( u \bullet v, w ) + (-1)^{|w|(|u|+|v|)}T(w, u \bullet v )\Big) \\
&= \frac{1}{2}\Big(-(-1)^{|u||v|}T( v, u \bullet w) - (-1)^{|w||v|}T( u\bullet w ,v )\Big) \\ &= -\frac{1}{2}(-1)^{|u||v|}\Big(T( v, u \bullet w) + (-1)^{|v|(|u|+|w|)}T( u\bullet w ,v )\Big)\\
&= -(-1)^{|u||v|}T_s(v, u \bullet w).
\end{aligned}
\]
It follows that $T_s(u \bullet v, w)=-(-1)^{|u||v|}T_s(v, u \bullet w)$.

Moreover, since $T$   satisfies \eqref{s2}, for all  $u,v,w \in \A$, we obtain
\[
\begin{aligned}
T_a(u \bullet v, w) &= \frac{1}{2}\Big(T( u \bullet v, w ) - (-1)^{|w|(|u|+|v|)}T(w, u \bullet v )\Big) \\
&\overset{\eqref{s2}}{=} \frac{1}{2}\Big(  (-1)^{|w|(|u|+|v|)}T(w, u \bullet v )-T( u \bullet v, w )\Big)
\end{aligned}
\]
Therefore,
\[
T_a(u \bullet v, w) = - T_a(u \bullet v, w), \text{
which implies } 
T_a(u \bullet v, w) = 0.
\]
We define the two superspaces:
\[
N := \{u \in \A \mid T_a(u,v)=0 \ \text{for all } v \in \A\} \text{ and } W := \{u \in \A \mid T_s(u,v)=0 \ \text{for all } v \in \A\}.
\] 
It is clear that $W$ and $N$ are subsuperspaces of $\mathcal{A}$. By the previous argument, we have $\A^2=\A\bullet\A \subseteq N$, and consequently $N$ is a two-sided ideal of $\A$. Since $T$ is non-degenerate, we have $N \cap W = \{0\}$. Hence, there exists a vector subsuperspace $V \subseteq \A$ such that $
\A = W \oplus V,$ and clearly $N \subseteq V$. Moreover, since $\A^2 \subseteq N \subseteq V$, it follows that $V$ is a two-sided ideal of $\A$ and that $T_s|_{V \times V}$ is non-degenerate.  By Lemma~\ref{evendimension}, the dimension of $V_{\bar{1}}$ is even. Since $\dim \mathcal{A}_{\bar{1}}$ is even and the fact that  $\mathcal{A}_{\bar{1}} = W_{\bar{1}} \oplus V_{\bar{1}}$, it follows that $\dim W_{\bar{1}}$ is also even.

Since $\dim W_{\bar{1}}$ is even, there exists an even non-degenerate symmetric bilinear form
$
H : W \times W \longrightarrow \mathbb{K}.
$ We define a bilinear form $\widetilde{H} : \A \times \A \to \mathbb{K}$ by
$
\widetilde{H}|_{W \times W} = H,$ and $ \widetilde{H}(V,\A)=\widetilde{H}(\A,V)=0.
$
Clearly, $\widetilde{H}$ is even and symmetric. Now, define the bilinear form
\[
\prs := T_s + \widetilde{H} : \A \times \A \longrightarrow \mathbb{K}.
\]
Since $T_s$ and $\widetilde{H}$ are both even and symmetric, it follows that $\prs$ is an even symmetric bilinear form. We will show that it is actually non-degenerate.

Let $u = w + v \in \A$, with $w \in W$ and $v \in V$, and assume that
$
\langle u,z\rangle=0,$  for all $z \in \A.$
Then
$
\widetilde{H}(w,z) + T_s(v,z)=0.$

If we choose $z\in W$, then $T_s(v,z)=0$, hence $\widetilde{H}(w,z)=0$, and thus $H(w,z)=0$. Since $H$ is non-degenerate, it follows that $w=0$.

If we choose $z\in V$, then $\widetilde{H}(w,z)=0$, and therefore $T_s(v,z)=0$. Since $T_s|_{V \times V}$ is non-degenerate, it follows that $v=0$.

Thus $u=0$, and consequently $\prs$ is non-degenerate.

Finally, since $\A^2 \subseteq V$, $\widetilde{H}(V,\A)=\widetilde{H}(\A,V)=0$, and $T_s$ satisfies \eqref{s1}, we obtain, for all $u,v,w \in \A$,
\[
\langle u \bullet v, w\rangle  = T_s(u \bullet v, w) = -(-1)^{|u||v|}T_s(v, u \bullet w) = -(-1)^{|u||v|} \langle v, u\bullet w\rangle .
\]
Therefore, $(\A,\bullet, \prs)$ is an even pseudo-Euclidean Novikov superalgebra.
\end{proof}

\sssbegin{Theorem} \label{Thm7.12}
Let $(\A,\bullet)$ be a Novikov superalgebra   and let $\Ll^\bullet$ and $\Rr^\bullet$ denote the left and right multiplication operators, respectively. Then there exists an odd non-degenerate symmetric bilinear form 
$
\prs : \A \times \A \longrightarrow \mathbb{K}
$
such that $(\A,\bullet, \prs)$ is an odd pseudo-Euclidean Novikov superalgebra if and only if there exists an odd isomorphism 
$
\Phi : \A \longrightarrow \A^{*}
$
satisfying \[
\Phi(\Ll_u^\bullet(a)) = (-1)^{|u |}(\Ll_u^\bullet)^{*}(\Phi(a)) \; \text{ and } \; \Phi(\Rr^\bullet_u(a)) = -(-1)^{|u |}(\Ll_u^\star)^{*}(\Phi(a)), 
\]
where the product $\star$  is antisymmetric. 
\end{Theorem}
\begin{proof}
Assume first that $(\A,\bullet, \prs)$ is an odd pseudo-Euclidean Novikov superalgebra. Then the result follows directly from Proposition~\ref{Prop7.8}.

Conversely, suppose that there exists an odd isomorphism $\Phi : \A \to \A^{*}$ such that
\[
\Phi(\Ll_u^\bullet(a)) = (-1)^{|u |}(\Ll_u^\bullet)^{*}(\Phi(a)) \; \text{and} \; \Phi(\Rr^\bullet_u(a)) = -(-1)^{|u |}(\Ll_u^\star)^{*}(\Phi(a)), \text{ for all $u,a \in \A$.}
\]
Consider the bilinear form
$
T : \A \times \A \longrightarrow \mathbb{K}$ defined by $T(u,v) := \Phi(u)(v),$
for all $u,v \in \A$. Since $\Phi$ is an odd isomorphism, $T$ is a non-degenerate odd bilinear form.

Following the same arguments as in the proof of Theorem~\ref{Thm7.11}, one can show that for all  $u,v,w \in \mathcal{A}$:
\begin{equation}\label{s1_odd}
T(u \bullet v, w) = -(-1)^{|u||v|} T(v, u\bullet w),
\end{equation}
and
\begin{equation}\label{s2_odd}
T(u \bullet v, w) = (-1)^{(|u|+|v|)|w|}T(w, u\bullet v).
\end{equation}

Now, we introduce the symmetric and anti-symmetric parts of $T$, denoted respectively by $T_s, T_a : \mathcal{A} \times \mathcal{A} \to \mathbb{K}$ and defined by
\[
T_s(u,v) := \frac{1}{2}\Big(T(u,v) + (-1)^{|u||v|}T(v,u)\Big), \quad
T_a(u,v) := \frac{1}{2}\Big(T(u,v) - (-1)^{|u||v|}T(v,u)\Big),
\]

for all  $u,v \in \mathcal{A}$. Clearly, $T = T_s + T_a$, and both $T_s$ and $T_a$ are odd bilinear forms.

Using \eqref{s1_odd}, a direct calculation shows that for all $u,v,w \in \mathcal{A}$
\[
T_s(u \bullet v, w) = -(-1)^{|u||v|}T_s(v, u \bullet w).
\]
Moreover, applying \eqref{s2_odd} yields
\[
T_a(u \bullet v, w) = \frac{1}{2}\Big(T( u \bullet v, w ) - (-1)^{|w|(|u|+|v|)}T(w, u \bullet v )\Big) \overset{\eqref{s2_odd}}{=} 0.
\]

We define the subsuperspaces:
\[
N := \{u \in \mathcal{A} \mid T_a(u,v)=0 \text{ for all } v \in \mathcal{A}\} \quad \text{and} \quad W := \{u \in \mathcal{A} \mid T_s(u,v)=0 \text{ for all } v \in \mathcal{A}\}.
\]
It is clear that $W$ and $N$ are subsuperspaces of $\mathcal{A}$. By the previous argument, we have $\mathcal{A}^2 = \mathcal{A} \bullet \mathcal{A} \subseteq N$, and thus $N$ is a two-sided ideal of $\mathcal{A}$. Since $T$ is non-degenerate, we have $N \cap W = \{0\}$. Therefore, there exists a subsuperspace $V \subseteq \mathcal{A}$ such that
$
\mathcal{A} = W \oplus V$ and $N \subseteq V.$ Because $\mathcal{A}^2 \subseteq N \subseteq V$, $V$ is a two-sided ideal of $\mathcal{A}$ and $T_s|_{V \times V}$ is a non-degenerate odd bilinear form.

Since $\Phi$ is an odd isomorphism, it immediately implies $\dim \mathcal{A}_{\bar{0}} = \dim \mathcal{A}_{\bar{1}}$. Moreover, since $T_s|_{V \times V}$ is an odd  non-degenerate bilinear form, Lemma~\ref{evendimension} implies that $\dim V_{\bar{0}} = \dim V_{\bar{1}}$. Consequently, using the decomposition $\mathcal{A} = W \oplus V$, we deduce that $\dim W_{\bar{0}} = \dim W_{\bar{1}}$.

Since $\dim W_{\bar{0}} = \dim W_{\bar{1}}$, there exists an odd  non-degenerate symmetric bilinear form 
$
H : W \times W \longrightarrow \mathbb{K}.
$ We define a bilinear form $\widetilde{H} : \A \times \A \to \mathbb{K}$ by
$
\widetilde{H}|_{W \times W} = H,$ and $ \widetilde{H}(V,\A)=\widetilde{H}(\A,V)=0.
$
Clearly, $\widetilde{H}$ is odd and symmetric. Now, define the bilinear form
\[
\prs := T_s + \widetilde{H} : \mathcal{A} \times \mathcal{A} \longrightarrow \mathbb{K}.
\]
As $T_s$ and $\widetilde{H}$ are both odd and symmetric, $\prs$ is an odd  symmetric bilinear form.

Proceeding similarly to the proof of Theorem~\ref{Thm7.11}, one easily verifies that $\prs$ is non-degenerate and satisfies 
\[
\langle u \bullet v, w \rangle = -(-1)^{|u||v|} \langle v, u \bullet w \rangle, \quad \text{for all } u, v, w \in \mathcal{A}.
\]
Therefore, $(\mathcal{A}, \bullet, \prs)$ is an odd pseudo-Euclidean Novikov superalgebra.
\end{proof}

\sssbegin{Remark}
Theorem~\ref{Thm7.11} remains valid when $(\A, \bullet, \prs)$ is a pseudo-Euclidean left-symmetric superalgebra.
\end{Remark}
\sssbegin{Definition}
Let $(\A, \cdot)$ be a non-associative superalgebra, and let $(V_1, r_1, l_1)$ and $(V_2, r_2, l_2)$ be two $\A$-bimodules. A homogeneous linear map $\Phi : V_1 \to V_2$ is called a \emph{morphism of $\A$-bimodules} if, for all $u \in \A$ and $x \in V_1$, the following identities hold:
\[
\Phi\big(r_1(u)(x)\big)
= (-1)^{|\Phi||u|} \, r_2(u)\big(\Phi(x)\big),\esp 
\Phi\big(l_1(u)(x)\big)
= (-1)^{|\Phi||u|} \, l_2(u)\big(\Phi(x)\big).
\]
Moreover, $\Phi$ is called an \emph{isomorphism of $\A$-bimodules} if it is bijective.
\end{Definition}

As a consequence of Theorems~\ref{Thm7.11} and~\ref{Thm7.12}, we obtain the following result for $\mathcal{A}$-bimodules.

\sssbegin{Corollary}
Let $(\mathcal{A}, \bullet)$ be a Novikov superalgebra, and let $\Ll^\bullet$ and $\Rr^\bullet$ denote the operators of left and right multiplication, respectively. 
\begin{enumerate}
    \item There exists an even pseudo-Euclidean structure $\prs$ on $\mathcal{A}$ if and only if $\dim \mathcal{A}_{\bar{1}}$ is even and the $\mathcal{A}$-bimodules $(\A, \Rr^\bullet, \Ll^\bullet)$ and $(\A^*, -(\Ll^\star)^*, (\Ll^\bullet)^*)$ are isomorphic via an even isomorphism.
    
    \item There exists an odd pseudo-Euclidean structure $\prs$ on $\mathcal{A}$ if and only if the $\mathcal{A}$-bimodules $(\A, \Rr^\bullet, \Ll^\bullet)$ and $(\A^*, -(\Ll^\star)^*, (\Ll^\bullet)^*)$ are isomorphic via an odd isomorphism (which directly implies $\dim \mathcal{A}_{\bar{0}} = \dim \mathcal{A}_{\bar{1}}$).
\end{enumerate}

Motivated by the previous results, we now present a new construction of pseudo-Euclidean Novikov superalgebras.
\end{Corollary}

\sssbegin{Theorem}\label{constriction}
Let $(\A,\bullet)$ be a left Leibniz $\Ll$-superalgebra. Assume that a bilinear map $\star : \A \times \A \to \A$ satisfies identities $(a)$, $(c)$, $(d)$, and $(e)$ of Proposition~\ref{Pr associe}. Then the vector superspace $T^*(\A)=\A \oplus \A^*$, endowed with the product
\[
(u+f)\ast (v+g)
= u \bullet v + (\Ll_u^\bullet)^*(g) -(-1)^{|u||v|} (\Ll_v^\star)^*(f),
\]
for all $u+f \in T^*(\A)_{\bar{|u|}}$ and $v+g\in  T^*(\A)_{\bar{|v|}} $, is a left Leibniz $\Ll$-superalgebra.

Moreover, the even bilinear form $\prs_c : T^*(\A) \times T^*(\A) \to \mathbb{K}$ defined by
\[
\langle u+f, v+g \rangle_c = f(v) + (-1)^{|u||v|} g(u),
\]
is non-degenerate and symmetric, and the left multiplication associated with $\ast$ is $\prs_c$-antisymmetric. Moreover, $(T^*(\A), \ast, \prs_c)$ is even pseudo-Euclidean left Leibniz $\Ll$-superalgebra \textup{(}or equivalently, Novikov superalgebra\textup{)}.

In this case, we have
\[
\langle (u+f)\ast (v+g), w+h \rangle_c
= \langle u+f, (v+g)\rhd (w+h) \rangle_c,
\]
where
\[
(v+g)\rhd (w+h)
= v\star w-g\circ \Rr_w^\bullet +(-1)^{|v||w|}h\circ\Rr_v^\bullet. \quad \text{ for all } u, v,w\in \A, \, f,g,h\in \A^*
\]
\end{Theorem}
\begin{proof}
By Proposition~\ref{Pr associe}, the pair 
$
\big(-(\Ll^\star)^*, (\Ll^\bullet)^* \big)
$ defines a representation of $\A$ on $\A^*$. Therefore, by the semidirect product construction, 
$(T^*(\A), \ast)$ is a left Leibniz 
$\Ll$-superalgebra.

By definition, the bilinear form $\prs_c$ is even, non-degenerate and symmetric. Let $u+f \in T^*(\A)_{\bar{|u|}}$, $v+g\in  T^*(\A)_{\bar{|v|}}$ and $w+h\in  T^*(\A)_{\bar{|w|}}$. Then
\begin{align*}
\langle (u+f)\ast (v+g), w+h \rangle_c
&= \left\langle u\bullet v + (\Ll_u^\bullet)^*(g) - (-1)^{|u||v|}(\Ll_v^\star)^*(f),\, w+h \right\rangle_c \\
&= -(-1)^{|u||v|} g(u \bullet w) +  f(v \star w) + (-1)^{|w|(|u|+|v|)} h(u \bullet v) \\
&= -(-1)^{|u||v|} \Big(g(u \bullet w) +(-1)^{(|u|+|w|)|v|}f(w \star v) \\
&\quad - (-1)^{|w|(|u|+|v|)+|u||v|} h(u \bullet v)\Big)  \\
&= -(-1)^{|u||v|} \Big(g(u \bullet w) +(-1)^{(|u|+|w|)|v|}f\circ \Ll^\star_w (v)\\
& \quad - (-1)^{|w|(|u|+|v|)+|u||v|} h\circ \Ll_u^\bullet  ( v)\Big)\\
&= -(-1)^{|u||v|} \Big(g(u \bullet w) -(-1)^{(|u|+|w|)|v|+|u||w|} (\Ll^\star_w)^*(f) (v) \\&\quad+(-1)^{|w|(|u|+|v|)}  (\Ll_u^\bullet)^*(h)  ( v)\Big)\\&=  -(-1)^{|u||v|} \left\langle v+g,  u\bullet w + (\Ll_u^\bullet)^*(g) - (-1)^{|u||w|}(\Ll_w^\star)^*(f),\, w+h \right\rangle_c\\&=  -(-1)^{|u||v|} \left\langle v+g, (u+f)\ast(w+h)  \right\rangle_c.
\end{align*}
Therefore, $(T^*(\A), \ast, \prs_c)$ is an even pseudo-Euclidean left Leibniz $\Ll$-superalgebra. 

On the other hand,
\begin{align*}
\langle u+f, (v+g)\rhd (w+h) \rangle_c
&= \langle (u+f)\ast (v+g), w+h\rangle_c\\&= \left\langle u\bullet v + (\Ll_u^\bullet)^*(g) - (-1)^{|u||v|}(\Ll_v^\star)^*(f),\, w+h \right\rangle_c \\
&= -(-1)^{|u|(|v|+|w|)} g\circ \Rr^\bullet_w(u) +  f(v \star w) \\
& \; + (-1)^{|w|(|u|+|v|)+|u||v|} h\circ \Rr_v^\bullet (u)
\\&= \langle u+f, v\star w-g\circ \Rr_w^\bullet +(-1)^{|v||w|}h\circ\Rr_v^\bullet \rangle_c.
\end{align*}
Since $\prs_c$ is nondegenerate, then 
$$
(v+g)\rhd (w+h)= v\star w-g\circ \Rr_w^\bullet +(-1)^{|v||w|}h\circ\Rr_v^\bullet .$$

This completes the proof.
\end{proof}

\sssbegin{Corollary}\label{PP}
Let $(\A, \bullet, \prs)$ be a   pseudo-Euclidean Novikov superalgebra. Then $(T^*(\A), \prs_c)$ is an even pseudo-Euclidean Novikov superalgebra where the product $\star$ is defined in Eq \eqref{produce associe}.
\end{Corollary}
\begin{proof}
Since  the product $\star$ defined in Eq \eqref{produce associe} satisfies identities $(a)$, $(c)$, $(d)$, and $(e)$ of Proposition~\ref{Pr associe}, then according to Theorem \ref{constriction},  $(T^*(\A), \prs_c)$ is an even pseudo-Euclidean Novikov superalgebra.
\end{proof}

We give an example of the construction of an odd pseudo-Euclidean Novikov superalgebra starting from an ordinary Novikov algebra.
\sssbegin{Example}
Let $(\A, \bullet)$ be a Novikov algebra.  Assume that a bilinear map $\star : \A \times \A \to \A$ satisfies identities $(a)$, $(c)$, $(d)$, and $(e)$ of Proposition~\ref{Pr associe}. Consider the $\mathbb{Z}_2$-graded vector space
$
\h = \h_{\bar{0}} \oplus \h_{\bar{1}}, $ with $ \h_{\bar{0}} := \A, \; \h_{\bar{1}} := \A^*.$ Define a bilinear product $\ast$ on $\h$ by
\[
(u + f) \ast (v + g) := u \bullet v + (\Ll_u^\bullet)^*(g)- (\Ll_v^\star)^*(f)
\quad \forall u,v \in \h_{\bar{0}}, \; f,g \in \h_{\bar{1}},
\]
for all $u, v\in \A$ and $f,g\in \A^* $.

Moreover, the bilinear form $\prs_c : T^*(\A) \times T^*(\A) \to \mathbb{K}$ defined by
\[
\langle u + f, v + g \rangle_c := f(v) + g(u), 
\quad \forall u,v \in \A, \; f,g \in \A^*.
\]
Then $(\h, \ast, \prs_c)$ is an odd pseudo-Euclidean Novikov superalgebra.
\end{Example}

Motivated by the previous construction, we now introduce an analogue of Theorem~\ref{constriction} in the odd pseudo-Euclidean setting.

Let $(\A, \cdot)$ be a superalgebra and let $\phi : \A \to \End(V)$ be a linear map. We now introduce the representation on $\Pi(V^*)$. Define
\[
\phi^\Pi : \A \to \End(\Pi(V^*))
\]
by
\[
(\phi^*)^\Pi(u) = \Pi \circ \phi^*(u) \circ \Pi,
\quad \text{for all } u \in \A.
\]
\sssbegin{Proposition}\label{rep1-b}
Let $(\A, \bullet,\prs)$ be a pseudo-Euclidean Novikov superalgebra, then the pair $
\big( -(\Ll^\star)^*)^\Pi, ((\Ll^\bullet)^*)^\Pi \big)
$
defines a representation of $\A$ on $\Pi(\A^*)$. 
\end{Proposition}
\begin{proof}
The proof is similar to that of Proposition~\ref{rep1}.\end{proof}
\sssbegin{Theorem}\label{constriction1}
Let $(\A,\bullet)$ be a left Leibniz $\Ll$-superalgebra. Assume that a bilinear map $\star : \A \times \A \to \A$ satisfies identities $(a)$, $(c)$, $(d)$, and $(e)$ of Proposition~\ref{Pr associe}. Then the vector superspace $(\Pi T^*)(\A):=\A \oplus \Pi(\A^*)$,  endowed with the product
\[
(u+\Pi(f))\ast (v+\Pi(g))
= u \bullet v + ((\Ll_u^\bullet)^*)^\Pi(\Pi(g)) -(-1)^{|u||v|} ((\Ll_v^\star)^*)^\Pi(\Pi(f)),
\]
for all \( u + \Pi(f) \in (\Pi T^*)(\A)_{|u|} \), \( v + \Pi(g) \in (\Pi T^*)(\A)_{|v|} \), is a left Leibniz $\Ll$-superalgebra.

Moreover, the odd bilinear form $\prs_c : (\Pi T^*)(\A) \times (\Pi T^*)(\A) \to \mathbb{K}$ defined by
\[
\langle u + \Pi(f),\, v + \Pi(g) \rangle_c = f(v) + (-1)^{|u||v|} g(u),
\]
is non-degenerate and symmetric, and the left multiplication associated with $\ast$ is $\prs_c$-antisymmetric. Moreover, $((\Pi T^*)(\A), \ast, \prs_c)$ is an odd pseudo-Euclidean left Leibniz $\Ll$-superalgebra.

In this case, we have the following
\[
\langle (u+\Pi(f))\ast (v+\Pi(g)), w+\Pi(h) \rangle_c
= \langle u+\Pi(f), (v+\Pi(g))\rhd (w+\Pi(h)) \rangle_c,
\]
where
\[
(v+\Pi(g))\rhd (w+\Pi(h))
= v\star w-\Pi(g\circ \Rr^\bullet_w)+(-1)^{|v||w|} \Pi(h\circ \Rr^\bullet_v) \quad \text{ for all } u, v,w\in \A, \, f,g,h\in \A^*
\]
\end{Theorem}
\begin{proof}
The proof is similar to that of Theorem~\ref{constriction}.\end{proof}
\sssbegin{Corollary}
Let $(\A, \bullet, \prs)$ be a  pseudo-Euclidean Novikov superalgebra, then $(\Pi(T^*(\A)), \prs_c)$ is odd pseudo-Euclidean Novikov superalgebra where  the product $\star$ is defined in Eq \eqref{produce associe}.
\end{Corollary}
\begin{proof}
The proof is similar to that of Corollary~\ref{PP}.\end{proof}

\section{Double extensions of pseudo-Euclidean Novikov superalgebras}\label{d-ext}
In this section, we introduce the notion of {\it   double extension} for pseudo-Euclidean Novikov superalgebras. It is inspired by the concept of double extension, originally introduced by Medina and Revoy in \cite{MR} for quadratic Lie algebras.

\subsection{Central extension of left-Leibniz \texorpdfstring{$\Ll$}{Ll}-superalgebras}
Let $(\A, \bullet)$ be a left-Leibniz $\Ll$-superalgebra,  let $\h := \K e$ be a one-dimensional vector superspace, and let $\mu : \A \times \A \rightarrow \K$ be a homogeneous bilinear map of a given parity.  
We define a new product $\diamond$ on the vector space
$
\widetilde{\A} := \A \oplus \K e
$
as follows:
\[
(u+\alpha e)\diamond(v+\beta e) := u \bullet v + \mu(u,v)e,
\qquad \text{ for all } u,v \in \A,\; \alpha, \beta \in \K.
\]
Thus,   $(\widetilde{\A}, \diamond)$ is a left-Leibniz $\Ll$-superalgebra, if and only if the following condition holds
\begin{equation}
\mu(u\bullet v, w)=0,\esp \mu(u,v\bullet w)=(-1)^{|u||v|}\mu(v, u\bullet w), \quad \text{for all $u,v,w\in \A$.} 
\label{con-left}
\end{equation}  

In this case, $(\widetilde{\A}, \diamond)$ is called the \emph{central extension} of $(\A, \bullet)$ by means of  $\mu$. 
\subsection{Almost semi-direct product of left-Leibniz \texorpdfstring{$\Ll$}{L}-superalgebras}

Let $(\A, \bullet)$ be a left-Leibniz $\Ll$-superalgebra. Let $\mathcal{V} := \K d$ be a one-dimensional vector superspace, let $
D : \A \rightarrow \A, \;
\xi : \A \rightarrow \A,
$ be homogeneous endomorphisms of parity $|d|$, and let $b_0 \in \A_{\bar 0}$  and a scalar $\lambda \in \K$.  
On the vector superspace 
$
\overline{\A} := \K d \oplus \A,
$
we define a new product $\bar{\bullet}$ as follows: (here $\lambda=0$ if $d$ is odd)
\[
d \bar{\bullet} d := \lambda d + b_0, \qquad
d \bar{\bullet} u := D(u), \qquad
u \bar{\bullet} d := \xi(u), \qquad
u \bar{\bullet} v := u \bullet v,
\]
for all $u,v \in \A$.  

Then $(\overline{\A}, \bar{\bullet})$ is a left-Leibniz $\Ll$-superalgebra if and only if the following conditions are satisfied (for all $u,v \in \A$): 
\begin{equation}
\begin{aligned}
\xi(u\bullet v)= & 0, & \xi(u)\bullet v=& 0, & D(u)\bullet v=&0, \\[0.2em] 
u\bullet \xi(v) =& (-1)^{|u||v|} v\bullet \xi(u),&  \Ll^\bullet_{b_0}=&0, &
D\circ \xi(u)=&(-1)^{|u||d|} u\bullet b_0,\\[0.2em] 
D(u\bullet v)=&(-1)^{|u||d|} u\bullet D(v), &\lambda=& 0,& 
\xi(b_0)=&0, \\[0.2em]  
(1 - (-1)^{|d|})D(b_0) = & 0,&\xi^2=&0, & \xi\circ D=&0,\\[0.2em] 
(1-(-1)^{|d|})D^2=& 0.
\end{aligned}
\label{produit-semi}
\end{equation}
\ssbegin{Definition}
If $(D, \xi, b_0)$ satisfies the compatibility conditions \eqref{produit-semi},  
then the  left-Leibniz $\Ll$-superalgebra $(\overline{\A}, \bar{\bullet})$ is called the 
\emph{almost semi-direct product} of the  left-Leibniz $\Ll$-superalgebra $(\A,\bullet)$ 
by the one-dimensional vector superspace $\K d$ by means of $(D, \xi, b_0)$.  
In this case, the 5-tuple 
$
(\mathcal{V}, \A, D, \xi, b_0)
$
is referred to as a \emph{context of almost semi-direct product of  left-Leibniz $\Ll$-superalgebra}.

\end{Definition}
\subsection{Double extensions of pseudo-Euclidean Novikov superalgebras}
In this subsection, we introduce the notion of {\it  double extension} of pseudo-Euclidean Novikov superalgebras.  

The term ``even extension" (or ``odd extension") will refer specifically to the parity of the extending element $d$, rather than the parity of the bilinear form itself.
\subsubsection{Even  double extensions of even pseudo-Euclidean Novikov superalgebras.}

Let $(\mathcal{B}, \bullet_\mathcal{B}, \prs_\mathcal{B})$ be an even pseudo-Euclidean Novikov superalgebra.   
Let $V = \K d$ be a one-dimensional vector superspace and $V^* = \K e$ denotes its dual.  

In order to introduce the process of  double extension by $V$, we consider two even linear maps $
\xi, D : \mathcal{B} \rightarrow \mathcal{B},
$  such that $D$ is $\prs_\mathcal{B}$-antisymmetric, and an even element $b_0 \in \mathcal{B}_{\bar 0}$. We extend the linear maps $\xi$  to 
$
\widetilde{D}, \widetilde{\xi} : \mathcal{B} \oplus V^* \rightarrow \mathcal{B} \oplus V^*,
$
as follows:   
\begin{equation}
\label{eq:double-tilde}
\begin{array}{c}
\widetilde{D}(u + \alpha e)
:= D(u) -\langle b_0, u\rangle_\mathcal{B}\, e, \quad   
\widetilde{\xi}(u + \alpha e)
:= \xi(u), \quad \text{for all $u \in \mathcal{B}$ and $\alpha \in \K$.}
\end{array}
\end{equation}

\medskip

\sssbegin{Lemma} \label{Lemma1}
The superspace $\mathcal{H} := \mathcal{B} \oplus V^*$ equipped with the product
\[
(u + \alpha e) \diamond (v + \beta e)
:= u \bullet v -\langle\xi(u), v\rangle_\mathcal{B}\, e,
\qquad \text{for all } \, u,v \in \mathcal{B},\; \alpha, \beta \in \K,
\]
is a left-Leibniz $\Ll$-superalgebra if and only if
\begin{equation}
\xi(u\bullet v)=0\esp u\bullet\xi(v)=(-1)^{|u||v|} v\bullet \xi(u), \quad \text{for all $u,v \in \mathcal{B}$.}
\label{eq:xi-condition0}
\end{equation}
\end{Lemma} 

\begin{proof} 
$(\mathcal{H}, \diamond)$ is left-Leibniz $\Ll$-superalgebra if and only if the bilinear map 
\[
\mu(u,v) := -\langle\xi(u), v\rangle_\mathcal{B}
\]
satisfies Eq. \eqref{con-left}, which is equivalent to Eq.  \eqref{eq:xi-condition0}. 
\end{proof}

\medskip

\sssbegin{Lemma}\label{Lemma2}
$(V, \mathcal{H} := \mathcal{B} \oplus V^*, \widetilde{D}, \widetilde{\xi}, b_0)$ is a context of almost semi-direct product 
of left-Leibniz $\Ll$-superalgebra if and only if,  for all $u,v \in \mathcal{B}$, the following conditions hold:
\begin{equation}\label{eq:claim2}
\begin{cases}
\xi(u\bullet v)=\xi(u)\bullet v=D(u)\bullet v=0, \\
u\bullet \xi(v)=(-1)^{|u||v|} v\bullet \xi(u),\quad D(u\bullet v) =u\bullet D(v), \\[2mm]
D\circ \xi=\Rr_{b_0}^\bullet,\quad \xi^2=\xi\circ D=\Ll^\bullet_{b_0}=0, \quad 
b_0\in \ker{\xi},
\end{cases}
\end{equation}

\end{Lemma}

\begin{proof}
 Let $u, v \in \mathcal{B}$ and $\alpha, \beta \in \mathbb{K}$.
 A straightforward computation shows that the following equivalences hold: \begin{align*}
\widetilde{\xi}\big((u+\alpha e)\diamond (v+\beta e)\big)=0
&\Longleftrightarrow \xi(u\bullet v)=0,\\[6pt]
\widetilde{\xi}(u+\alpha e)\diamond (v+\beta e)=0
&\Longleftrightarrow \xi(u)\bullet v=0\; \text{ and $\xi^2=0$},\\[6pt]
\widetilde{D}(u+\alpha e)\diamond (v+\beta e)=0
    &\Longleftrightarrow D(u)\bullet v=0,\quad \text{$\xi\circ D=0.$}
    \\[6pt]
\widetilde{\xi}^2(u+\alpha e)=0
    &\Longleftrightarrow \text{$\xi^2=0,$}
    \\[6pt]
\widetilde{\xi}\circ \widetilde{D}(u+\alpha e)=0
    &\Longleftrightarrow  \text{$\xi\circ D=0.$}
\end{align*}
Moreover, the equation $$X\diamond \widetilde{\xi}(Y)=(-1)^{|X||Y|} Y\diamond \widetilde{\xi}(X), \quad \text{ for all } X, Y\in \h$$
is equivalent to
\[
u\bullet\xi(v)= (-1)^{|u||v|}v\bullet \xi(u).
\]
Further
\[\widetilde{D}((u+\al e) \diamond (v+\beta e))  =(u+\al e) \diamond \widetilde{D}(v+\beta e)\]
is equivalent to
$$D(u\bullet v)=u\bullet D(v), \esp D\circ \xi=\Rr^\bullet_{b_0}.$$
In addition, the equation $\widetilde{D}\circ \widetilde{\xi}(X)= X\bullet b_0$ for all $X\in \h$
is equivalent to
$$D\circ \xi(u) = u\bullet b_0.$$
Finally, the last identities are equivalent to the conditions
$$\Ll_{b_0}^\bullet=0, \esp \xi(b_0)=0.$$
 This proves the lemma.
\end{proof}

We now introduce the notion of  even double extension of an even pseudo-Euclidean Novikov  superalgebra by a one-dimensional vector space.

\sssbegin{Theorem}\label{double-ex1}

Let $(\mathcal{B}, \bullet_\mathcal{B}, \prs_\mathcal{B})$ be an even pseudo-Euclidean Novikov superalgebra. Let $V = \mathbb{K}d$ be a one-dimensional vector space and $V^* = \mathbb{K}e$ denotes its dual. Assume that there exist two even linear maps
 $\xi, D : \mathcal{B} \to \mathcal{B}$, 
and an element $b_{0} \in \mathcal{B}_{\bar 0}$ such that $D$ is $\prs_\mathcal{B}$-anti-symmetric and satisfying the system \eqref{eq:claim2}. Define 
$
\A := \mathbb{K}d \oplus \mathcal{B} \oplus \mathbb{K}e,
$
equipped with the product $\bullet$
\begin{equation}\label{Produit1}
\begin{array}{lll} \displaystyle 
e \bullet u = u \bullet e =0,& e \bullet d = d \bullet e = 0,& \displaystyle  
d \bullet u = D(u) -\langle b_0, u\rangle_\mathcal{B} e, \\[0.3em]
\displaystyle  u \bullet d = \xi(u), &d \bullet d =  b_0, & 
u \bullet v = u \bullet_{\mathcal{B}} v - \langle \xi(u), v\rangle_\mathcal{B} e, 
\end{array}
\end{equation}
for all $ u,v \in \mathcal{B}$, and an even bilinear form $\prs$ given by 
$$
\prs|_{\mathcal{B} \times \mathcal{B}} = \prs_{\mathcal{B}},
\qquad \langle e,d\rangle = \langle d,e\rangle  = 1,
\qquad \langle d,\mathcal{B}\rangle = \langle e,\mathcal{B}\rangle= \{0\}.
$$
Then, $(\A, \bullet, \prs)$ is an even pseudo-Euclidean Novikov superalgebra. 
\end{Theorem}
The pseudo-Euclidean Novikov superalgebra $(\A, \bullet, \prs)$ is called the even  double extension of $(\mathcal{B}, \bullet_\mathcal{B}, \prs_{\mathcal{B}})$ by means of $(D, \xi, b_0).$
\begin{proof}
According to Proposition \ref{Pr1}, $(\A, \bullet, \prs)$ is a pseudo-Euclidean Novikov superalgebra if and only if $(\A, \bullet, \prs)$ is a pseudo-Euclidean left Leibniz $\Ll$-superalgebra.

It is easy to check that the left multiplication operators 
$\Ll^\bullet$ are $\prs$-antisymmetric.  Now, we show that $(\A, \bullet)$ is a left-Leibniz $\Ll$-superalgebra.  We have 
\[
(u + \alpha e) \bullet (v + \beta e) = u \bullet_{\mathcal{B}} v -\langle \xi(u), v\rangle_\mathcal{B} e, \quad \text{for all $u,v \in \mathcal{B}$ and $\alpha, \beta \in \mathbb{K}$.}
\]
From system~\eqref{eq:claim2}, for all $u,v\in \mathcal{B}$, we have
\[
\xi(u\bullet v)=0, \qquad 
u\bullet\xi(v)=(-1)^{|u||v|} v\bullet \xi(u).
\]
Then, Lemma~\ref{Lemma1} implies that the above product defines on the vector superspace
$\mathcal{H}=\mathcal{B}\oplus \mathbb{K}e$ a left Leibniz $\Ll$-superalgebra structure.
Moreover, it is a central extension of $\mathcal{B}$ by $V$  by means of  $\mu$, where
$
\mu(u,v)=-\langle \xi(u),v\rangle_{\mathcal{B}},$
 for all  $u,v\in\mathcal{B}$.  

Now,  since $(D, \xi, b_0)$ satisfies the system \eqref{eq:claim2}, it follows from Lemma \ref{Lemma2}, that 
$
(V, \mathcal{H} := \mathcal{B} \oplus V^*, \widetilde{D}, \widetilde{\xi}, b_0)
$
is a context of an almost semi-direct product of left-Leibniz  $\Ll$-superalgebras. 
Hence, $(\A, \bullet)$ is a left-Leibniz  $\Ll$-superalgebra, 
which implies that $(\A, \bullet, \prs)$ is an even  pseudo-Euclidean Novikov superalgebra.
\end{proof}

\sssbegin{Remark}\label{TT}
Let $(\mathcal{B}, \bullet_{\mathcal B}, \prs_{\mathcal B})$ be an even  pseudo-Euclidean Novikov superalgebra such that $\prs_{\mathcal B}$ is even, and let
$(\mathcal A, \bullet, \prs)$ be the even double extension of
$(\mathcal{B}, \bullet_{\mathcal B}, \prs_{\mathcal B})$
by means of $(D,\xi,b_0)$.
From the definition of the product, it is clear that if $(\mathcal A,\bullet)$ is nilpotent, then
$(\mathcal{B}, \bullet_{\mathcal B})$ is nilpotent as well.
However, the converse does not hold in general.

Indeed, let $(\mathcal{B}, \bullet_{\mathcal B}, \prs_{\mathcal B})$
be an even pseudo-Euclidean trivial superalgebra and let $V = \mathbb{K}d$ be a one-dimensional vector space and $V^* = \mathbb{K}e$ be its dual.
Assume that $\xi = 0$ and $b_0 = 0$, and choose an arbitrary endomorphism
$D : \mathcal{B} \to \mathcal{B}$ that is not nilpotent and $\prs_{\mathcal B}$-antisymmetric.
Then $(D,0,0)$ satisfies the system~\eqref{eq:claim2}, and hence
$(\mathcal A:= \mathbb{K} e \oplus \mathcal{B}\oplus \mathbb{K} d, \bullet, \prs)$ is an even double extension of pseudo-Euclidean Novikov superalgebra $(\mathcal{B}, \bullet_\mathcal{B}, \prs_\mathcal{B})$.
Moreover, from the product~\eqref{Produit1}, we have
\[
d \bullet u = D(u), \quad \text{for all } u \in \mathcal{B},
\]
and consequently
\[
(\Ll_d^\bullet)^n(u) = D^n(u), \quad \text{for all } u \in \mathcal{B} , n\in \mathbb{N}.
\]
Since $D$ is not nilpotent, it follows that $(\mathcal A,\bullet)$ is not nilpotent and consequently 
$\A^-$ is also not nilpotent.
\end{Remark}

\subsubsection{Odd double extensions of  even pseudo-Euclidean Novikov superalgebras}
We now introduce the notion of  {\it odd double extension} of even pseudo-Euclidean Novikov superalgebras by a purely odd one-dimensional vector superspace.

Let $(\mathcal{B}, \bullet_\mathcal{B}, \prs_\mathcal{B})$ be an even  pseudo-Euclidean Novikov superalgebra. Let $V = \K d$ be a purely odd one-dimensional vector superspace and let $V^* = \K e$ be its dual.  

In order to introduce the process of the odd  double extension, we consider two odd linear maps $
\xi, D : \mathcal{B} \rightarrow \mathcal{B},
$  such that $D$ is $\prs_\mathcal{B}$-anti-symmetric, and two even elements $b_0, c_0 \in \mathcal{B}_{\bar 0}$. We extend the linear maps $D, \xi$  to 
$
\widetilde{D}, \widetilde{\xi} : \mathcal{B} \oplus V^* \rightarrow \mathcal{B} \oplus V^*,
$
as follows:   
\begin{equation}
\label{eq:double-tilde1}
\begin{array}{c}
\widetilde{D}(u + \alpha e)
:= D(u) -\langle b_0, u\rangle_\mathcal{B}\, e, \quad   
\widetilde{\xi}(u + \alpha e)
:= \xi(u)+ \langle c_0, u\rangle_\mathcal{B}\, e, 
\end{array}
\end{equation}
for all $u \in \mathcal{B}$ and $\alpha \in \K$.

\medskip

\sssbegin{Lemma}\label{Lemma3} 
The superspace $\mathcal{H} := \mathcal{B} \oplus V^*$ equipped with the product
\[
(u + \alpha e) \diamond (v + \beta e)
:= u \bullet v -(-1)^{|v|} \langle\xi(u), v\rangle\, e,
\qquad \text{for all }\, u,v \in \mathcal{B},\; \alpha, \beta \in \K,
\]
is a left-Leibniz $\Ll$-superalgebra if and only if
\begin{equation}
\xi(u\bullet v)=0\esp u\bullet\xi(v)=(-1)^{|u||v|} v\bullet \xi(u), \text{ for all $u,v \in \mathcal{B}$.}
\label{eq:xi-condition}
\end{equation}

\end{Lemma}

\begin{proof}
$(\mathcal{H}, \diamond)$ is left-Leibniz $\Ll$-superalgebra if and only if the bilinear map 
\[
\mu(u,v) := -(-1)^{|v|}\langle\xi(u), v\rangle_\mathcal{B}
\]
satisfies condition \eqref{con-left}, which is equivalent to equations \eqref{eq:xi-condition}.\end{proof}

\medskip

\sssbegin{Lemma}\label{Lemma4}
$(V, \mathcal{H} := \mathcal{B} \oplus V^*, \widetilde{D}, \widetilde{\xi}, \widetilde{b_0}, 0)$ is a context of almost semi-direct product 
of left-Leibniz $\Ll$-superalgebras if and only if 
\begin{equation}\label{eq:claim4}
\begin{cases}
\xi(u\bullet v)=\xi(u)\bullet v=D(u)\bullet v=0, \\
u\bullet \xi(v)=(-1)^{|u||v|} v\bullet \xi(u),\quad D(u\bullet v) =(-1)^{|u|}u\bullet D(v), \\[2mm]
D\circ \xi(u)=(-1)^{|u|}\Rr_{b_0}^\bullet(u),\quad \xi^2=\xi^*\circ\xi=\xi\circ D=D^2=\Ll^\bullet_{b_0}=\Rr_{c_0}^{\bullet}=0, \\  
b_0\in \ker{\xi}\cap \ker{D},\,  c_0\in \ker\xi^*\cap \ker D,\\
\langle b_0, c_0\rangle_\mathcal{B}=\langle b_0, b_0\rangle_\mathcal{B}=0, 
\end{cases}
\end{equation}
for all $u, v \in \mathcal{B}$. 
\end{Lemma}
\begin{proof}
 Let $u, v \in \mathcal{B}$ and $\alpha, \beta \in \mathbb{K}$. 
A straightforward computation shows the following equivalences:
\begin{align*}
\widetilde{\xi}\big((u+\alpha e)\diamond (v+\beta e)\big)=0
&\Longleftrightarrow
\xi(u\bullet v)=0,
\quad \text{ and } \quad \Rr_{c_0}^{\bullet}=0,\\[6pt]
\widetilde{\xi}(u+\alpha e)\diamond (v+\beta e)=0
&\Longleftrightarrow
\xi(u)\bullet v=0,
\quad \text{ and } \quad\xi^2=0,\\[6pt]
\widetilde{D}(u+\alpha e)\diamond (v+\beta e)=0
&\Longleftrightarrow
D(u)\bullet v=0,
\quad \text{ and  } \quad \xi\circ D=0,
\\[6pt]
\widetilde{\xi}^2(u+\alpha e)=0
    &\Longleftrightarrow \text{$\xi^2=0,$} \quad\quad\quad\quad \text{ and  } \quad \xi^*(c_0)=0,
    \\[6pt]
\widetilde{\xi}\circ \widetilde{D}(u+\alpha e)=0
    &\Longleftrightarrow \text{$\xi\circ D=0,$} \quad\quad\quad \text{ and  } \quad D(c_0)=0.
\end{align*}
Moreover, the equation $$X\diamond \widetilde{\xi}(Y)=(-1)^{|X||Y|} Y\diamond \widetilde{\xi}(X), \quad \text{ for all }  X, Y\in\h, $$
is equivalent to
\[
u\bullet\xi(v)= (-1)^{|u||v|}v\bullet \xi(u), \esp (-1)^{|v|+|d|}\langle \xi(u), \xi(v)\rangle_\mathcal{B}= (-1)^{|u|+|d|} \langle \xi(v), \xi(u)\rangle_\mathcal{B}.
\]
The latter equality is equivalent to
\[
\langle \xi^{*}\circ \xi(u), v\rangle_{\mathcal B}
= - \langle \xi^{*}\circ \xi(u), v\rangle_{\mathcal B},
\]
 which is equivalent to
\[
\xi^{*}\circ \xi = 0.
\]
Further 
\[\widetilde{D}(X \diamond Y ) =(-1)^{|X|}X \diamond \widetilde{D}(Y),\quad \text{ for all } X,Y\in \h, \]
is equivalent to
$$D(u\bullet v)=(-1)^{|u|}u\bullet D(v), \esp \langle u\bullet b_0, v\rangle_\mathcal{B}=(-1)^{|v|}\langle D\circ \xi(u), v\rangle_\mathcal{B}. $$
The latter equality is equivalent to
\[
\langle u\bullet b_0, v\rangle_\mathcal{B}=(-1)^{|v|+|u|+|u|}\langle D\circ \xi(u), v\rangle_\mathcal{B}.
\]
Since both $D \circ \xi$ and $b_0$ are even, and $\prs$ is even, it follows that when $|u| + |v| = \bar{1}$, the equality reduces to the trivial identity $0 = 0$. Consequently, the only nontrivial case occurs when $|u| + |v| = \bar{0}$, which yields the relation:
\[
D\circ \xi(u) = (-1)^{|u|}u\bullet b_0.
\]

The equation $\widetilde{D}\circ \widetilde{\xi}(X)=(-1)^{|X|} X\diamond b_0$, for all $X\in \h$, 
is equivalent to
$$D\circ \xi(u) = (-1)^{|u|}u\bullet b_0, \esp \langle b_0, \xi(u)\rangle_\mathcal{B}=-(-1)^{|u|} \langle \xi(u), b_0\rangle_\mathcal{B}.$$
Since both $b_0$ and $\prs$ are even, and $\xi$ is odd,  we deduce that $  \langle b_0, \xi(u)\rangle_\mathcal{B}=-(-1)^{|u|} \langle \xi(u), b_0\rangle_\mathcal{B}$.

Finally, the last identities are equivalent to the conditions
$$\Ll_{b_0}^\bullet=D^2=0, \esp \xi(b_0)=D(b_0)=0,  \esp \langle b_0, c_0\rangle_\mathcal{B}=\langle b_0, b_0\rangle_\mathcal{B}=0 .$$
 This proves the lemma.
\end{proof}

 Now, we are in a position to introduce the odd double extensions of  even pseudo-Euclidean Novikov superalgebras.

\sssbegin{Theorem}\label{double-ex2}
Let $(\mathcal{B}, \bullet_\mathcal{B}, \prs_\mathcal{B})$ be an even  pseudo-Euclidean Novikov superalgebra.
Let $V = \mathbb{K}d$ be a purely odd one-dimensional
vector superspace and let $V^* = \mathbb{K}e$ be its dual. Assume that there exist two odd homogeneous linear maps 
$D, \xi :\mathcal{B} \to \mathcal{B}$, 
and two even elements $b_0 , c_0\in \mathcal{B}_{\bar{0}}$, such that $D$ is $\prs_\mathcal{B}$-anti-symmetric and  satisfy the system \eqref{eq:claim4}. Define the vector superspace
$
\A := \mathbb{K}d \oplus \mathcal{B} \oplus \mathbb{K}e,
$
equipped with the product $\bullet$ given by: 
\begin{equation}\label{Produit2}
\begin{array}{lll} \displaystyle 
e \bullet u = u \bullet e =0,& e \bullet d = d \bullet e = 0,& \displaystyle  
d \bullet u = D(u) -\langle b_0, u\rangle_\mathcal{B} e, \\[0.3em]
\displaystyle  u \bullet d = \xi(u)+\langle c_0, u\rangle_\mathcal{B}e, &d \bullet d =  b_0, & 
u \bullet v = u \bullet_{\mathcal{B}} v -(-1)^{|v|} \langle \xi(u), v\rangle_\mathcal{B} e, 
\end{array}
\end{equation}
for all $ u,v \in \mathcal{B}$,  and an even bilinear form $\prs$ given by 
$$
\prs|_{\mathcal{B} \times \mathcal{B}} = \prs_{\mathcal{B}},
\qquad \langle e,d\rangle =- \langle d,e\rangle  = 1,
\qquad \langle d,\mathcal{B}\rangle = \langle e,\mathcal{B}\rangle= \{0\}.
$$
Then, $(\A, \bullet, \prs)$ is an even pseudo-Euclidean Novikov superalgebra.

\end{Theorem}
The even pseudo-Euclidean Novikov superalgebra $(\A, \bullet, \prs)$ is called the odd  double extension of $(\mathcal{B}, \bullet_\mathcal{B}, \prs_\mathcal{B})$ by means of $(D,\xi, b_0, c_0)$.
\begin{proof}
The proof is analogous to that of Theorem~\ref{double-ex1}, and makes use of Lemmas~\ref{Lemma3} and~\ref{Lemma4}.
\end{proof}

\sssbegin{Proposition}\label{nil-ex}
Let $(\mathcal{B}, \bullet_{\mathcal B}, \prs_{\mathcal B})$
be an even pseudo-Euclidean Novikov superalgebra, and let
$(\mathcal A, \bullet, \prs)$ be the odd double extension of
$(\mathcal{B}, \bullet_{\mathcal B}, \prs_{\mathcal B})$
by means of $(D,\xi,b_0,c_0)$.
Then $(\mathcal A,\bullet)$ is nilpotent if and only if
$(\mathcal{B}, \bullet_{\mathcal B})$ is nilpotent.
\end{Proposition}

\begin{proof}
Assume first that $(\mathcal A,\bullet)$ is nilpotent.
According to the product~\eqref{Produit2}, we have
\[
u \bullet v
= u \bullet_{\mathcal B} v
- (-1)^{|v|}\langle \xi(u), v\rangle_{\mathcal B}\, e, \quad   \text{for all $u,v\in\mathcal B$}.
\]
By induction, for any $n \ge 1$
\[
(\Ll_u^\bullet)^n(v)
= (\Ll_u^{\bullet_{\mathcal B}})^n(v)
- (-1)^{|v| + (n-1)|u|}\,
\langle \xi(u), (\Ll_u^{\bullet_{\mathcal B}})^{n-1}(v)\rangle_{\mathcal B}\, e.
\]
Since $(\mathcal A,\bullet)$ is nilpotent, there exists $n_0\in\mathbb N$ such that
$(\Ll_u^\bullet)^{n_0}=0$ for all $u\in\mathcal B$.
Hence $(\Ll_u^{\bullet_{\mathcal B}})^{n_0}=0$ for all $u\in\mathcal B$, and by
Proposition~\ref{Nil}, $(\mathcal B,\bullet_{\mathcal B})$ is nilpotent.

Conversely, assume that $(\mathcal B,\bullet_{\mathcal B})$ is nilpotent.  Let $u, v \in \mathcal{B}$ and $\lambda, \al \in \mathbb{K}$. Using the product formula \eqref{Produit2} together with the identities in \eqref{eq:claim4}, we compute by induction $\Ll^{\bullet n}_{\lambda d+u+\al e}(v)$ for $n \ge 1$:

For $n = 1$
    \[
    \Ll^{\bullet}_{\lambda d+u+\al e}(v) = (\lambda d + u+\al e) \bullet v = \Ll^{\bullet_\mathcal{B}}_u(v) + \lambda D(v) - \Big( \lambda \langle b_0, v \rangle_{\mathcal{B}} + (-1)^{|v|} \langle \xi(u), v \rangle_{\mathcal{B}} \Big) e.
    \]
For $n = 2$
    Using $D^2 = 0$, $D(u \bullet_{\mathcal{B}} v) = (-1)^{|u|} u \bullet_{\mathcal{B}} D(v)$, $D(u) \bullet_{\mathcal{B}} v = 0$, and $\A \bullet e = \{0\}$, we obtain
    \[
    \begin{aligned}
    \Ll^{\bullet 2}_{\lambda d+u+\al e}(v) &= (\lambda d + u) \bullet \Big( \Ll^{\bullet_\mathcal{B}}_u(v) + \lambda D(v) \Big) \\
    &= \Ll^{\bullet_\mathcal{B} 2}_u(v) + \lambda (1 + (-1)^{|u|}) \Ll^{\bullet_\mathcal{B}}_u(D(v)) \\
    &\quad - \Big( \lambda \langle b_0, x_1 \rangle_{\mathcal{B}} + (-1)^{|u|+|v|} \langle \xi(u), y_1 \rangle_{\mathcal{B}} \Big) e,
    \end{aligned}
    \]
    where $x_1 = \Ll^{\bullet_\mathcal{B}}_u(v) + \lambda D(v)$, $y_1 = \Ll^{\bullet_\mathcal{B}}_u(v) - \lambda(-1)^{|u|} D(v)$ . 
    \[
\begin{aligned}
\Ll_{\lambda d+u+\al e}^{\bullet\,2}(v)
={}&
\Ll_u^{\bullet_{\mathcal B}2}(v)
+
\lambda\bigl(1+(-1)^{|u|}\bigr)
\Ll_u^{\bullet_{\mathcal B}}\bigl(D(v)\bigr)
\\
&-
\Bigl(
\lambda\langle b_0,x_1\rangle_{\mathcal B}
+
(-1)^{|u|+|v|}
\langle\xi(u),x_1\rangle_{\mathcal B}
\Bigr)e,
\end{aligned}
\]

By an immediate induction on $n \ge 1$, we obtain 
\[
\Ll^{\bullet n}_{\lambda d+u+\al e}(v) = \Ll^{\bullet_\mathcal{B} n}_u(v) + C_n \Ll^{\bullet_\mathcal{B} n-1}_u(D(v)) - \Big( \lambda \langle b_0, x_{n-1} \rangle_{\mathcal{B}} + (-1)^{(n-1)|u|+|v|} \langle \xi(u), y_{n-1} \rangle_{\mathcal{B}} \Big) e,
\]
where $C_n = \lambda \sum_{i=0}^{n-1} (-1)^{i|u|} $,  $x_{n-1} = \Ll^{\bullet_\mathcal{B} n-1}_u(v) + C_{n-1} \Ll^{\bullet_\mathcal{B} n-2}_u(D(v))$, and $y_{n-1} = \Ll^{\bullet_\mathcal{B} n-1}_u(v) -(-1)^{|u|} C_{n-1} \Ll^{\bullet_\mathcal{B} n-2}_u(D(v))$. 

Since $(\mathcal{B}, \bullet_{\mathcal{B}})$ is nilpotent of index $n_0$ (that is, $\Ll^{\bullet_\mathcal{B} n_0}_{w} = 0$ for all $w \in \mathcal{B}$) and $e \bullet \A = \A \bullet e = \{0\}$, then we have $\Ll^{\bullet n_0+2}_{\la d+u+\al e}(v) = 0$.

Next, using the product formula \eqref{Produit2}, we have
\[
\Ll^{\bullet}_{\lambda d+u+\al e}(d) = (\lambda d + u) \bullet d = \lambda b_0 + \xi(u) + \langle c_0, u \rangle_{\mathcal{B}} e.
\]
Since $\A \bullet e = \{0\}$ and $\lambda b_0 + \xi(u) \in \mathcal{B}$, using the result established for elements in $\mathcal{B}$, we deduce that
\[
\Ll^{\bullet n_0+3}_{\lambda d+u+\al e}(d) = \Ll^{\bullet n_0+2}_{\lambda d+u+\al e}(\lambda b_0 + \xi(u)) = 0.
\]
Thus, $\Ll_x^\bullet$ is nilpotent for all $x \in \mathcal{A}$. By Proposition~\ref{Nil},  $(\mathcal{A}, \bullet)$ is nilpotent.
\end{proof}

\subsubsection{Even  double extension of odd pseudo-Euclidean Novikov superalgebras}
Following the approach of Theorem~\ref{double-ex1}, we now define the even  double extension for  odd pseudo-Euclidean Novikov superalgebras. Combining Lemmas~\ref{Lemma1} and~\ref{Lemma2}, we obtain the following result.
\sssbegin{Theorem}\label{double-ex3}
Let $(\mathcal{B}, \bullet_\mathcal{B}, \prs_\mathcal{B})$ be an odd pseudo-Euclidean Novikov superalgebra.
Let $V = \mathbb{K}d$ be a one-dimensional vector space and $\Pi(V^*) = (\mathbb{K} e)_{\bar{1}}$, where $V^*$ denotes its dual. Assume that there exist two even linear maps
$D, \xi : \mathcal{B} \to \mathcal{B}$, 
and an even element $b_0 \in \mathcal{B}_{\bar 0}$, such that $D$ is $\prs_\mathcal{B}$-anti-symmetric and  and satisfies the system \eqref{eq:claim2}.
 Define the vector superspace
$
\A := \mathbb{K}d \oplus \mathcal{B} \oplus \mathbb{K}e,
$
equipped with the product $\bullet$ given by: 
\begin{equation}\label{Produit3}
\begin{array}{lll} \displaystyle 
e \bullet u = u \bullet e =0,& e \bullet d = d \bullet e = 0,& \displaystyle  
d \bullet u = D(u) -\langle b_0, u\rangle_\mathcal{B} e, \\[0.3em]
\displaystyle  u \bullet d = \xi(u), &d \bullet d =  b_0, & 
u \bullet v = u \bullet_{\mathcal{B}} v - \langle \xi(u), v\rangle_\mathcal{B} e, 
\end{array}
\end{equation}
for all $ u,v \in \mathcal{B}$, and an even bilinear form $\prs$ given by 
$$
\prs|_{\mathcal{B} \times \mathcal{B}} = \prs_{\mathcal{B}},
\qquad \langle e,d\rangle = \langle d,e\rangle  = 1,
\qquad \langle d,\mathcal{B}\rangle = \langle e,\mathcal{B}\rangle= \{0\}.
$$
Then, $(\A, \bullet, \prs)$ is an odd pseudo-Euclidean Novikov superalgebra.

\end{Theorem}
The  odd pseudo-Euclidean Novikov superalgebra $(\A, \bullet, \prs)$ is called the even  double extension of $(\mathcal{B}, \bullet_\mathcal{B}, \prs_\mathcal{B})$ by means of $(D,\xi, b_0)$.

\begin{proof} The proof is similar to that of Theorem~\ref{double-ex1}, using Lemmas~\ref{Lemma1} and~\ref{Lemma2}.
\end{proof}
\subsubsection{Odd  double extension of odd pseudo-Euclidean Novikov superalgebras}
We now introduce the notion of  odd  double extension of odd pseudo-Euclidean Novikov superalgebras by a purely odd one-dimensional vector superspace.

Let $(\mathcal{B}, \bullet_\mathcal{B}, \prs_\mathcal{B})$ be an odd pseudo-Euclidean Novikov superalgebra. Let $V = \K d$ be a purely odd one-dimensional vector superspace and $\Pi(V^*) = \K e$, where $V^*$ denotes its dual.  

To introduce the process of odd double extension by $V$, we consider two  odd linear maps $
D, \xi : \mathcal{B} \rightarrow \mathcal{B},$ and even elements $b_0, c_0\in \mathcal{B}_{\bar 0}$. 
These two maps can be extended as 
$
\widetilde{D}, \widetilde{\xi} : \mathcal{B} \oplus \Pi(V^*) \rightarrow \mathcal{B} \oplus \Pi(V^*),
$ as follows:  
\begin{equation}
\widetilde{D}(u + \alpha e)
= D(u)+ \, \langle b_0, u\rangle_\mathcal{B} e,\quad 
\widetilde{\xi}(u + \alpha e)
= \xi(u)+\langle c_0, u\rangle_\mathcal{B} e, 
\label{eq:double-tilde2}
\end{equation}
for all $u \in \mathcal{B}$ and $\alpha \in \K$. Put $\widetilde{b_0}
=  b_0+\la e$. 
\sssbegin{Lemma}\label{Lemma5}
The superspace $\mathcal{H} := \mathcal{B} \oplus V^*$ equipped with the product
\[
(u + \alpha e) \diamond (v + \beta e)
:= u \bullet v -(-1)^{|v|} \langle\xi(u), v\rangle\, e,
\qquad \text{ for all } u,v \in \mathcal{B},\; \alpha, \beta \in \K,
\]
is a left-Leibniz $\Ll$-superalgebra if and only if
\begin{equation}
\xi(u\bullet v)=0\esp u\bullet\xi(v)=(-1)^{|u||v|} v\bullet \xi(u), \quad \text{for all $u,v \in \mathcal{B}$.}
\label{eq:xi-condition3}
\end{equation}

\end{Lemma}
\begin{proof}
The proof is similar to that of Lemma \ref{Lemma1}. \end{proof}

\sssbegin{Lemma}\label{Lemma6}
$(V, \mathcal{H} := \mathcal{B} \oplus \Pi(V^*), \widetilde{D}, \widetilde{\xi}, \widetilde{b_0})$ is a context of almost semi-direct product 
of left-Leibniz $\Ll$-superalgebras if and only if 
\begin{equation}\label{eq:claim6}
\begin{cases}
\xi(u\bullet v)=\xi(u)\bullet v=D(u)\bullet v=0, \\
u\bullet \xi(v)=(-1)^{|u||v|} v\bullet \xi(u),\quad D(u\bullet v) =(-1)^{|u|}u\bullet D(v),\quad  D\circ \xi(u) = (-1)^{|u|}u\bullet b_0 \\[2mm]
\xi^2=D^2=\xi^*\circ\xi=\xi\circ D=\Ll^\bullet_{b_0}=\Rr_{c_0}^\bullet=0, \\
b_0\in \ker{\xi}\cap \ker\xi^*\cap \ker{D},\quad c_0\in \ker\xi^*\cap \ker{D},
\end{cases}
\end{equation}
for all $u, v \in \mathcal{B}$.
\end{Lemma}

\begin{proof}
 Let $u, v \in \mathcal{B}$ and $\alpha, \beta \in \mathbb{K}$.  A straightforward computation shows the following equivalences:
\begin{align*}
\widetilde{\xi}((u+\al e)\diamond (v+\beta e))=0& \Longleftrightarrow  \xi(u\bullet v)=0,\esp \Rr_{c_0}^\bullet=0,\\
\widetilde{\xi}(u+\al e) \diamond (v+\beta e)  =0& \Longleftrightarrow \xi(u)\bullet v=0,\esp \xi^2=0,\\
\widetilde{D}(u+\al e) \diamond (v+\beta e)  =0& \Longleftrightarrow D(u)\bullet v=0,\esp \xi\circ D=0, \\[6pt]
\widetilde{\xi}^2(u+\alpha e)=0
    &\Longleftrightarrow \text{$\xi^2=0,$} \quad\quad\quad\quad \text{ and  } \quad \xi^*(c_0)=0,
    \\[6pt]
\widetilde{\xi}\circ \widetilde{D}(u+\alpha e)=0
    &\Longleftrightarrow \text{$\xi\circ D=0,$} \quad\quad\quad \text{ and  } \quad D(c_0)=0.
\end{align*}
Moreover, the equation 
$$X\diamond \widetilde{\xi}(Y)=(-1)^{|X||Y|} Y\diamond \widetilde{\xi}(X), \quad \text{ for all } X, Y\in \h,$$
is equivalent to
\[
u\bullet\xi(v)= (-1)^{|u||v|}v\bullet \xi(u), \esp (-1)^{|v|+|d|}\langle \xi(u), \xi(v)\rangle_\mathcal{B}= (-1)^{|u|+|d|} \langle \xi(v), \xi(u)\rangle_\mathcal{B}.
\]
The latter equality is equivalent to
\[
\langle \xi^{*}\circ \xi(u), v\rangle_{\mathcal B}
= - \langle \xi^{*}\circ \xi(u), v\rangle_{\mathcal B},
\text{ 
 which is equivalent to }
\xi^{*}\circ \xi = 0.
\]
Further 
\[\widetilde{D}(X \diamond Y)  =(-1)^{|X|}X \diamond \widetilde{D}(Y),\quad \text{ for all } X, Y\in \h, \]
is equivalent to
$$D(u\bullet v)=(-1)^{|u|}u\bullet D(v), \esp \langle u\bullet b_0, v\rangle_\mathcal{B}=-(-1)^{|v|}\langle D\circ \xi(u), v\rangle_\mathcal{B}.$$
The latter equality is equivalent to
\[
\langle u\bullet b_0, v\rangle_\mathcal{B}=-(-1)^{|v|+|u|+|u|}\langle D\circ \xi(u), v\rangle_\mathcal{B}.
\]
Since both $D \circ \xi$ and $b_0$ are even, and $\prs$ is odd, it follows that when $|u| + |v| = \bar{0}$, the equality reduces to the trivial identity $0 = 0$. Consequently, the only nontrivial case occurs when $|u| + |v| = \bar{1}$, which yields the relation
\[
D\circ \xi(u) = (-1)^{|u|}u\bullet b_0.
\]
The equation $\widetilde{D}\circ \widetilde{\xi}(X)=(-1)^{|X|} X\diamond b_0$, for all $X\in \h$, 
is equivalent to
$$D\circ \xi(u) = (-1)^{|u|}u\bullet b_0, \esp \langle b_0, \xi(u)\rangle_\mathcal{B}=-(-1)^{|u|} \langle \xi(u), b_0\rangle_\mathcal{B}.$$
The latter equality is equivalent to
$$\langle \xi^*(b_0), u\rangle_\mathcal{B}=-(-1)^{|u|} \langle  \xi^*(b_0), u\rangle_\mathcal{B}.$$
Since both  $\prs$ and $\xi$ are odd, and $b_0$ is even, it follows that when $|u|= \bar{1}$, the equality reduces to the trivial identity $0 = 0$. Consequently, the only nontrivial case occurs when $|u|  = \bar{0}$, which yields the relation 
$$\xi^*(b_0)=0.$$
Finally, the last identities are equivalent to the conditions
$$\Ll_{b_0}^\bullet=D^2=0, \esp \xi(b_0)=D(b_0)=0.$$
 This proves the lemma.
\end{proof}

\sssbegin{Theorem}\label{double-ex4}
Let $(\mathcal{B}, \bullet_\mathcal{B}, \prs_\mathcal{B})$ be an  odd pseudo-Euclidean Novikov superalgebra. Let $V=V_{\bar{1}} = (\mathbb{K} d)_{\bar{1}}$ be a purely odd one-dimensional vector superspace  and $\Pi(V^*) = \mathbb{K} e$ where $V^*$ denotes its dual.  Assume that there exist two odd linear maps 
$D, \xi :\mathcal{B} \to \mathcal{B}$, 
and two even elements  $b_0 , c_0\in \mathcal{B}_{\bar{0}}$, such that $D$ is $\prs_\mathcal{B}$-anti-symmetric and  satisfying the system \eqref{eq:claim6}. Define the vector superspace
$
\A := \mathbb{K}d \oplus \mathcal{B} \oplus \mathbb{K}e,
$
equipped with the product $\bullet$
\begin{equation}\label{Produit4}
\begin{array}{lll} \displaystyle 
e \bullet u = u \bullet e =0,& e \bullet d = d \bullet e = 0,& \displaystyle  
d \bullet u = D(u) +\langle b_0, u\rangle_\mathcal{B} e, \\[0.3em]
\displaystyle  u \bullet d = \xi(u)+\langle c_0, u\rangle_\mathcal{B}e, &d \bullet d =  b_0, & 
u \bullet v = u \bullet_{\mathcal{B}} v -(-1)^{|v|} \langle \xi(u), v\rangle_\mathcal{B} e, 
\end{array}
\end{equation}
for all $ u,v \in \mathcal{B}$,  and an odd bilinear form $\prs$ given by
$$
\prs|_{\mathcal{B} \times \mathcal{B}} = \prs_{\mathcal{B}},
\qquad \langle e,d\rangle = \langle d,e\rangle  = 1,
\qquad \langle d,\mathcal{B}\rangle = \langle e,\mathcal{B}\rangle= \{0\}.
$$
Then, $(\A, \bullet, \prs)$ is an odd pseudo-Euclidean Novikov superalgebra.

\end{Theorem}
The odd pseudo-Euclidean Novikov superalgebra $(\A, \bullet, \prs)$ is called the odd  double extension of $(\mathcal{B}, \bullet_\mathcal{B}, \prs_\mathcal{B})$ by means of $(D,\xi, b_0, c_0)$.

\begin{proof}
The structure of this proof is analogous to that of Theorem~\ref{double-ex1}; the key difference is the use of Lemmas~\ref{Lemma5} and~\ref{Lemma6} instead of the lemmas used there.
\end{proof}

\sssbegin{Proposition}\label{nil-ex2}
Let $(\mathcal{B}, \bullet_{\mathcal B}, \prs_{\mathcal B})$
be an odd pseudo-Euclidean Novikov superalgebra, and let
$(\mathcal A, \bullet, \prs)$ be the odd double extension of
$(\mathcal{B}, \bullet_{\mathcal B}, \prs_{\mathcal B})$
by means of $(D,\xi,b_0,c_0)$.
Then $(\mathcal A,\bullet)$ is nilpotent if and only if
$(\mathcal{B}, \bullet_{\mathcal B})$ is nilpotent.
\end{Proposition}

\begin{proof}
The proof is similar to that of Proposition \ref{nil-ex}.
\end{proof}

\section{Structure of Pseudo-Euclidean Novikov Superalgebras}\label{con-ex}

\ssbegin{Proposition}\label{reduit}  
Let $(\A, \bullet, \prs)$ be a pseudo-Euclidean Novikov superalgebra.  
Choose  a non-zero homogeneous element $e\in  \A\bullet\A \cap (\A\bullet\A)^{\perp}$,  and denote  
$
I := \mathbb{K} e.$ Consider    $\mathcal{B} := I^{\perp} / I, 
$
and let $\pi_{\mathcal{B}} : I^{\perp} \to \mathcal{B}$ be the canonical projection. We have
\begin{enumerate}
    \item[$(i)$] $I$  is a totally isotropic graded two-sided ideal of $(\A, \bullet)$ and  $I^\perp$ is a graded two-sided ideal of $(\A, \bullet)$. 
    \item[$(ii)$] The quotient $\mathcal{B} = I^{\perp} / I$ admits a canonical structure of a   pseudo-Euclidean Novikov superalgebra, defined by
    \[
        \pi_{\mathcal{B}}(u)\bullet_\mathcal{B} \pi_{\mathcal{B}}(v) = \pi_{\mathcal{B}}(u\bullet v) \quad \text{and} \quad 
     \langle \pi_{\mathcal{B}}(u), \pi_{\mathcal{B}}(v)\rangle_\mathcal{B} = \langle u, v\rangle, \quad \text{for all $u,v \in I^{\perp}$.}
    \] 

\end{enumerate}
\end{Proposition}
\begin{proof}
Let us prove Part ($i$). Since  $e\in \A\bullet \A\cap(\A\bullet \A)^\bot$, according to Proposition \ref{Pr1} and Eq. \eqref{ortogonal},  we have $\Ll^\bullet_e=\Rr^\bullet_e=0$, which implies that $I$ is a totally isotropic graded two-sided ideal of $(\A, \bullet)$.

Now, let $v \in I^{\perp}$ and $w \in \A$. We have 
    \[
       \langle w \bullet v, e\rangle = -(-1)^{|w||v|} \langle v, w \bullet e\rangle = 0, 
        \qquad 
        \langle v \bullet w, e\rangle = -(-1)^{|w||v|} \langle w, v \bullet e\rangle = 0,
    \]
    hence $w \bullet v,\, v \bullet w \in I^{\perp}$.
    
In the case where the form $\prs$ is even, then $I^{\perp}$ is a graded two-sided ideal of $(\A, \bullet)$. 

In the case where the form  $\prs$ is odd, we distinguish two cases.  
If $e$ is even, and since $\prs$ is odd, then $\A_{\bar{0}} \subseteq I^{\perp}$, 
and therefore $(I^{\perp})_{\bar{0}} = \A_{\bar{0}}$. 
If $e$ is odd, then $\A_{\bar{1}} \subseteq  I^{\perp}$, and hence 
$(I^{\perp})_{\bar{1}} = \A_{\bar{1}}$. 
In both cases, we conclude that $I^{\perp}$ is also a graded two-sided ideal of $(\A, \bullet)$.
 
Let us prove Part $(ii)$. Since $I^{\perp}$ is a graded two-sided ideal of $(\A, \bullet)$, the quotient space $\mathcal{B} = I^{\perp} / I$ naturally inherits a structure of a pseudo-Euclidean Novikov superalgebra. Indeed, for $u, v \in I^{\perp}$, we define
    \[
        \pi_{\mathcal{B}}(u)\bullet \pi_{\mathcal{B}}(v) := \pi_{\mathcal{B}}(u\bullet v), 
        \qquad 
        \langle\pi_{\mathcal{B}}(u), \pi_{\mathcal{B}}(v)\rangle_\mathcal{B} := \langle u, v\rangle.
    \]
    These operations are well-defined since $I$ is totally isotropic, i.e., $\langle u, e\rangle=0$ for all $u \in I^{\perp}$ and $e \in I$.

    Moreover, since $(\A, \bullet, \prs)$ is a pseudo-Euclidean Novikov superalgebra, so is $(\mathcal{B}, \bullet_\mathcal{B}, \prs_\mathcal{B})$. 
\end{proof}

The pseudo-Euclidean Novikov superalgebra $\mathcal{B}$ obtained in Prop.~\ref{reduit} is referred to as the  \emph{pseudo-Euclidean Novikov superalgebra deduced from $I^{\perp}$ by means of $I$}.\\

We now prove the converse of Theorems \ref{double-ex1} and \ref{double-ex2}. 

\ssbegin{Theorem}\label{central}
Let $(\A, \bullet, \prs)$ be an even   pseudo-Euclidean Novikov superalgebra such that $\A\bullet \A$ is degenerate.
Then $(\A, \bullet, \prs)$ is either:
\begin{enumerate}
    \item[$(i)$] an even double extension of an even pseudo-Euclidean Novikov superalgebra $(\mathcal{B}, \bullet_\mathcal{B}, \prs_\mathcal{B})$ by means of $(D, \xi, b_0)$, or
    \item[$(ii)$] an odd double extension of an  even pseudo-Euclidean Novikov superalgebra $(\mathcal{B}, \bullet_\mathcal{B}, \prs_\mathcal{B})$ by means of $(D, \xi, b_0, c_0)$.
\end{enumerate}
\end{Theorem}
\begin{proof}
Since $\A\bullet\A\cap (\A\bullet\A)^\bot\neq 0$, 
by Prop. \ref{reduit}, there exists a totally isotropic two-sided graded ideal \( I = \K e \) of \( (\A, \bullet) \), and its orthogonal \( I^\perp \) that is also a two-sided graded ideal of \( (\A, \bullet) \).  We now distinguish two cases.

\underline{The case where  \( e \) is even.} Since \( \prs \) is even and non-degenerate, there exists an element \( d \in \A_{\bar{0}} \setminus \{0\} \) such that
$
\langle e, d\rangle  = \langle d, e\rangle = 1.
$

Since $I\subseteq I^\perp$, there exists a subsuperspace ${\mathfrak B}$ such that  $I^\perp = I \oplus {\mathfrak B}$ and  the restriction  $\prs_{{\mathfrak B}}=\prs|_{{\mathfrak B}\times {\mathfrak B}}$  
 is non-degenerate.  Let us write
$
\A = \mathbb{K} e \oplus {\mathfrak B} \oplus \mathbb{K} d,
$
and since \( I^\perp = I \oplus {\mathfrak B} \) is a two-sided graded ideal of \( (\A, \bullet) \), for any  \( u, v \in {\mathfrak B} \), we have
\[
u \bullet v = u \bullet_{\mathfrak B} v + \mu(u, v)e,
\]
where \( \mu : {\mathfrak B} \times {\mathfrak B} \to \K \) is an even bilinear form, and 
\( \bullet_{\mathfrak B} : {\mathfrak B} \times {\mathfrak B} \to {\mathfrak B} \) is an even bilinear map.  

It follows that \( ({\mathfrak B}, \bullet_{\mathfrak B}) \) is a left-Leibniz  $\Ll$-superalgebra,
\( \mu \) satisfies the compatibility condition \eqref{con-left},
and the restriction \( \prs_{\mathfrak B} = \prs|_{{\mathfrak B} \times {\mathfrak B}} \) is an  even pseudo-Euclidean structure on \( ({\mathfrak B}, \bullet_{\mathfrak B}) \). Moreover, the canonical projection
$
\pi : {\mathfrak B} \rightarrow I^\perp / I=\mathcal{B}$,  defined by $\pi(u) = \pi_{\mathcal{B}}(u)$,  for all $  u \in {\mathfrak B},$ is an isomorphism of left-symmetric superalgebras. Thus, we can identify  \( {\mathfrak B} \cong I^\perp / I=\mathcal{B} \), and hence
$
{\cal A} = \mathbb{K} e \oplus \mathcal{B} \oplus \mathbb{K} d,
$
where \( (\mathcal{B}, \bullet_{\mathcal{B}}, \prs_{\mathcal{B}}) = (I^\perp / I, \bullet_{\mathcal{B}}, \prs_\mathcal{B}) \)
is a  pseudo-Euclidean Novikov superalgebra.

Since \( I^\perp =  \mathbb{K} e \oplus \mathcal{B} \) is a two-sided graded ideal of  \( (\A, \bullet) \), then, for all \( u,v \in \mathcal{B} \), the product $\bullet$ on \( \A \) is  given by
\[
\begin{aligned}
u \bullet v &= u \bullet_\mathcal{B} v + \mu(u,v)e, &
u \bullet d &= \xi(u) + f(u)e,& 
d \bullet u &= D(u) + g(u)e, 
\\ 
d \bullet d & =  b_0+\al e+\la d, & 
e\bullet d&=d\bullet e=0, &  e\bullet u &=u\bullet e=0,
\end{aligned}
\]
where \( \xi, D \in \mathrm{End}(\mathcal{B})_{\bar{0}} \),  \( f, g : \mathcal{B} \to \K \) are even linear maps, and $\al, \la\in \K$.

Since the left multiplication operators $\Ll^\bullet$ are
$\prs$-antisymmetric, for any
$u,v \in \mathcal{B}$, we have
\[
\mu(u,v) = \langle u \bullet v, d\rangle = -(-1)^{|u||v|}\langle v, u \bullet d\rangle
= -\langle\xi(u), v\rangle_\mathcal{B}.
\]
Similarly,
\[
g(u) = \langle d \bullet u, d\rangle = -\langle  u, d \bullet d\rangle=-\langle b_0, u\rangle_\mathcal{B}.
\]
Moreover,
\[
\langle D(u), v\rangle_\mathcal{B}=\langle d \bullet u, v\rangle = -\langle  u, d \bullet v\rangle=-\langle u, D(v)\rangle_\mathcal{B}.
\]
It follows that $D$ is antisymmetric with respect to $\prs_\mathcal{B}$.

Moreover, we have 
\[
f(u)=\langle u\bullet d, d\rangle=-\langle d, u\bullet d\rangle =-f(u),
\]
hence \(f=0 \).

Furthermore, we have
$$\al=\langle d\bullet d,  d\rangle=-\langle d, d\bullet d\rangle=-\al, \quad \la=\langle d\bullet d, e\rangle=-\langle d, d\bullet e\rangle =0,$$
hence $\al=\la=0$.

Therefore, the product $\bullet$ takes the following simplified form:
\[
\begin{aligned}
u \bullet v &= u \bullet_\mathcal{B} v -\langle\xi(u), v\rangle_\mathcal{B}e, & 
u \bullet d &= \xi(u) , & 
d \bullet u &= D(u) -\langle b_0, u\rangle_\mathcal{B}e, \\
d \bullet d & =  b_0,&
e\bullet d&=d\bullet e=0,&  e\bullet u & =u\bullet e=0.
\end{aligned}
\]
Now, it is easy to show that \( (\A, \bullet) \) is a Novikov superalgebra if and only if the triple \( (D, \xi, b_0) \) satisfies the system \eqref{eq:claim2}. Hence, \( (\A, \bullet, \prs) \) is the even  double extension of the  pseudo-Euclidean Novikov superalgebra \( (\mathcal{B}, \bullet_\mathcal{B}, \prs_\mathcal{B}) \) by means of \( (D, \xi, b_0) \).

\underline{ The case, where \( e \) is odd.}  Since \( \prs \) is even and non-degenerate, there exists an element \( d \in \A_{\bar{1}} \setminus \{0\} \) such that
$
\langle e, d\rangle =- \langle d, e\rangle  = 1.
$

Since $I\subseteq I^\perp$, there exists a subsuperspace $\mathfrak{B}$ such that  $I^\perp = I \oplus \mathfrak{B}$ and  the restriction  $\prs_\mathfrak{B}=\omega|_{\mathcal{B}\times \mathcal{B}}$  
 is non-degenerate.  Let us write
$
\A = \mathbb{K} e \oplus {\mathfrak B} \oplus \mathbb{K} d,
$
and since \( I^\perp = I \oplus {\mathfrak B} \) is a two-sided graded ideal of \( (\A, \bullet) \), for any \( u, v \in {\mathfrak B} \), we have
\[
u \bullet v = u \bullet_{\mathfrak B} v + \mu(u, v)e,
\]
where \( \mu : {\mathfrak B} \times {\mathfrak B} \to \K \) is an even bilinear form, and 
\( \bullet_{\mathfrak B} : {\mathfrak B} \times {\mathfrak B} \to {\mathfrak B} \) is an even bilinear map.  

It follows that \( ({\mathfrak B}, \bullet_{\mathfrak B}) \) is a left-Leibniz $\Ll$-superalgebra,
\( \mu \) satisfies the compatibility condition \eqref{con-left},
and the restriction \( \prs_{\mathfrak B} = \prs|_{{\mathfrak B} \times {\mathfrak B}} \) is an even pseudo-Euclidean structure on\( ({\mathfrak B}, \bullet_{\mathfrak B}) \). Moreover, the canonical projection
$
\pi : {\mathfrak B} \rightarrow I^\perp / I=\mathcal{B}$, defined by $ \pi(u) = \pi_{\mathcal{B}}(u)$,  for all $  u \in {\mathfrak B},$ is an isomorphism of   left-Leibniz $\Ll$-superalgebras. Thus, we can identify  \( {\mathfrak B} \cong I^\perp / I=\mathcal{B} \), and hence
$
\A = \mathbb{K} e \oplus \mathcal{B} \oplus \mathbb{K} d,
$
where \( (\mathcal{B}, \bullet_{\mathcal{B}}, \prs_{\mathcal{B}}) = (I^\perp / I, \bullet_{\mathcal{B}}, \prs_\mathcal{B}) \)
is a pseudo-Euclidean Novikov superalgebra.

Since \( I^\perp =  \mathbb{K} e \oplus \mathcal{B} \) is a two-sided graded ideal of  \( (\A, \bullet) \).
Then, for all \( u,v \in \mathcal{B} \), the products $\bullet$ on \( \A \) are given by
\[
\begin{aligned}
u \bullet v &= u \bullet_\mathcal{B} v + \mu(u,v)e, & 
u \bullet d & = \xi(u) + f(u)e, & 
d \bullet u &= D(u) + g(u)e, \\ 
d \bullet d &=  b_0, & 
e\bullet d&=d\bullet e=0, &  e\bullet u &=u\bullet e=0,
\end{aligned}
\]
where \( \xi, D \in \mathrm{End}(\mathcal{B})_{\bar{1}} \), and \( f, g : \mathcal{B} \to \K \) are even linear maps.

Since the left multiplication operators $\Ll^\bullet$ are
$\prs$-antisymmetric, for any
$u,v \in \mathcal{B}$ we have
\[
\mu(u,v) = \langle u \bullet v, d\rangle = -(-1)^{|u||v|}\langle v, u \bullet d\rangle
= -(-1)^{|v|}\langle\xi(u), v\rangle_\mathcal{B}.
\]
Similarly,
\[
g(u) = \langle d \bullet u, d\rangle = -(-1)^{|u |}\langle  u, d \bullet d\rangle=-(-1)^{|u|}\langle b_0, u\rangle_\mathcal{B},
\]
since $\prs$ is even and $b_0$ is even. It follows that when $|u|=\bar{1}$, the
equality reduces to the trivial identity   $0=0$. Consequently, the nontrivial case occurs when $|u|=\bar{0}$, leading to the following relation
$$g(u) =-\langle b_0, u\rangle_\mathcal{B}.$$
Moreover,
\[
\langle D(u), v\rangle_\mathcal{B}=\langle d \bullet u, v\rangle = -(-1)^{|u|}\langle  u, d \bullet v\rangle=-(-1)^{|u|}\langle u, D(v)\rangle_\mathcal{B}.
\]
It follows that $D$ is $\prs_\mathcal{B}$-antisymmetric.

Furthermore, we have 
\[
f(u)=\langle u\bullet d, d\rangle=-(-1)^{|u|}\langle d, u\bullet d\rangle =(-1)^{|u|}f(u).
\]
Since $f$ is an even linear map and the bilinear form
$\prs_{\mathcal{B}}$ is non-degenerate, there exists an
element $c_0 \in \mathcal{B}_0$ such that
\[
f(u)=\langle c_0, u\rangle_{\mathcal{B}}, \quad \text{for all } u\in\mathcal{B}.
\]
Therefore, the product $\bullet$ takes the simplified form:
\[
\begin{aligned}
u \bullet v &= u \bullet_\mathcal{B} v -(-1)^{|v|}\langle\xi(u),  v\rangle_\mathcal{B}e, &
u \bullet d & = \xi(u)+\langle c_0, u\rangle_\mathcal{B} e, & 
d \bullet u &= D(u) -\langle b_0, u\rangle_\mathcal{B}e,\\ 
d \bullet d & =  b_0, & 
e\bullet d&=d\bullet e=0,&  e\bullet u &=u\bullet e=0.
\end{aligned}
\]
Now, it is easy to show that \( (\A, \bullet) \) is a Novikov superalgebra if and only if the  quadruple \( (D, \xi, b_0, c_0) \) satisfies the system \eqref{eq:claim4}. Hence, \( (\A, \bullet, \prs) \) is the odd  double extension of the  even pseudo-Euclidean Novikov superalgebra \( (\mathcal{B}, \bullet_\mathcal{B}, \prs_\mathcal{B}) \) by means of \( (D, \xi, b_0, c_0) \).
\end{proof}

\ssbegin{Corollary}\label{start-abel}

Every even pseudo-Euclidean non-trivial nilpotent Novikov superalgebra  $(\A,\bullet, \prs)$ can be obtained by a finite sequence of even or odd double extensions of even pseudo-Euclidean nilpotent Novikov superalgebras, starting from the trivial superalgebra. 
\end{Corollary}
\begin{proof}
If $(\A,\bullet)$ is a non-trivial nilpotent Novikov superalgebra, then $\A^-$ is also nilpotent. According to Corollary \ref{nilpotent}, $\A\bullet\A$ is degenerate. Hence, by Theorem \ref{central}, 
$(\A, \bullet, \prs)$ is an even or an odd  double extension of an even pseudo-Euclidean Novikov superalgebra
$(\mathcal{B}_1,\bullet_{\mathcal{B}_1},  \prs_{\mathcal{B}_1})$ by means of  $(D_1, \xi_1, b_1)$ or $(D_1, \xi_1, b_1,  c_1)$. If $(\mathcal{B}_1,\bullet_{\mathcal{B}_1})$ is a non-trivial superalgebra, then it is itself an even or an odd  double extension of an even pseudo-Euclidean Novikov superalgebra
$(\mathcal{B}_2,\bullet_{\mathcal{B}_2},  \prs_{\mathcal{B}_2})$ by means of  $(D_2, \xi_2, b_2)$  or $(D_2, \xi_2, b_2,  c_2)$. This process has to stop since $\A$ is finite-dimensional. 
\end{proof}

\ssbegin{Theorem} \label{Milnor1}
Any even pseudo-Euclidean Novikov superalgebra is either a Milnor  superalgebra or can be
obtained by a sequence of even or odd double extension from a Milnor superalgebra
\end{Theorem}
\begin{proof}
If $\A\bullet\A$ is non-degenerate, then according to Theorem \ref{ ideal non degenere}, $(\A, \bullet)$ is Milnor superalgebra. If $\A\bullet\A$ is degenerate, then according to Theorem \ref{central}, 
 $(\A,\bullet,\prs)$ is an even or odd double extension
of an even pseudo-Euclidean Novikov superalgebra
$(\mathcal{B}_1,\bullet_{\mathcal{B}_1},\prs_{\mathcal{B}_1})$
by means of $(D_1,\xi_1,b_1)$  or $(D_1, \xi_1, b_1,  c_1)$.
If 
$\mathcal{B}_1\bullet_{\mathcal{B}_1}\mathcal{B}_1$ is non-degenerate, then
$(\mathcal{B}_1,\bullet_{\mathcal{B}_1},
\prs_{\mathcal{B}_1})$
is a Milnor superalgebra.
If it is degenerate, then
$(\mathcal{B}_1,\bullet_{\mathcal{B}_1},
\prs_{\mathcal{B}_1})$
is itself an even or odd double extension of an even pseudo-Euclidean Novikov
superalgebra
$(\mathcal{B}_2,\bullet_{\mathcal{B}_2},
\prs_{\mathcal{B}_2})$
by means of $(D_2,\xi_2,b_2)$  or $(D_2, \xi_2, b_2,  c_2)$.
Since $\A$ is finite-dimensional, this process must terminate after finitely
many steps.
\end{proof}

Now, we give the converse of Theorem \ref{double-ex3} and Theorem \ref{double-ex4}.
\ssbegin{Theorem}\label{central1}
Let $(\A, \bullet, \prs)$ be an odd   pseudo-Euclidean Novikov superalgebra such that $\A\bullet \A$ is degenerate.
Then $(\A, \bullet, \prs)$ is either:
\begin{enumerate}
    \item[$(i)$] an even double extension of an odd pseudo-Euclidean Novikov superalgebra $(\mathcal{B}, \bullet_\mathcal{B}, \prs_\mathcal{B})$ by means of $(D, \xi, b_0)$, or
    \item[$(ii)$] an odd double extension of an  odd pseudo-Euclidean Novikov superalgebra $(\mathcal{B}, \bullet_\mathcal{B}, \prs_\mathcal{B})$ by means of $(D, \xi, b_0, c_0)$.
\end{enumerate}
\end{Theorem}
\begin{proof}
The proof is similar to that of Theorem \ref{central}.
\end{proof}
\ssbegin{Corollary}\label{start-abel2}
Every non-trivial nilpotent odd pseudo-Euclidean Novikov superalgebra \allowbreak $(\A,\bullet, \prs)$ can be obtained by a finite sequence of even or odd double extensions of odd pseudo-Euclidean nilpotent Novikov superalgebras, starting from $\{0\}$.
\end{Corollary}
\begin{proof}
If $(\A, \bullet)$ is non-trivial nilpotent superalgebra,  then $\A^-$ is also nilpotent. According to Corollary \ref{nilpotent}, $\A\bullet\A$ is degenerate. Hence,   by Theorem \ref{central1}, 
$(\A, \bullet, \prs)$ is an even or an odd  double extension of an odd pseudo-Euclidean Novikov superalgebra
$(\mathcal{B}_1,\bullet_{\mathcal{B}_1},  \prs_{\mathcal{B}_1})$ by means of  $(D_1, \xi_1, b_1)$ or $(D_1, \xi_1, b_1,  c_1)$. If $(\mathcal{B}_1,\bullet_{\mathcal{B}_1})$ is a non-trivial superalgebra, then it is itself an even or an odd  double extension of an odd pseudo-Euclidean Novikov superalgebra
$(\mathcal{B}_2,\bullet_{\mathcal{B}_2},  \prs_{\mathcal{B}_2})$ by means of  $(D_2, \xi_2, b_2)$  or $(D_2, \xi_2, b_2,  c_2)$.  Iterating this process, and since the superdimension of $\mathcal A$ is
finite, there exists $k\in\mathbb N$ such that $\mathcal B_k$ is a trivial
superalgebra. Consequently,
$(\mathcal B_k,\prs_{\mathcal B_k})$
is an odd pseudo-Euclidean trivial superalgebra.

By Lemma~\ref{evendimension}, the superdimension of $\mathcal B_k$ is even and satisfies
\(
\dim (\mathcal B_k)_{\bar 0}=\dim (\mathcal B_k)_{\bar 1}.
\)
Therefore, there exists a totally isotropic supersubspace of dimension $1$.
This implies that
$(\mathcal B_k,\prs_{\mathcal B_k})$
is itself a double extension of an odd pseudo-Euclidean nilpotent Novikov superalgebra
\(
(\mathcal B_{k+1},\prs_{\mathcal B_{k+1}}).
\)
Proceeding inductively, we obtain a sequence of double extensions until there
exists $n\in\mathbb N$ such that $\mathcal B_n=\{0\}$.
\end{proof}

\ssbegin{Theorem}
Any odd pseudo-Euclidean Novikov superalgebra is either a Milnor  superalgebra or can be
obtained by a sequence of even or odd double extension from a Milnor superalgebra
\end{Theorem}
\begin{proof}The proof is similar to that of Theorem \ref{Milnor1}.\end{proof}
\section{Classification of pseudo-Euclidean Novikov superalgebras of low
dimensions}\label{classi}
In this section, we provide a classification of pseudo-Euclidean Novikov superalgebras of dimension at most 
4, over the field of real numbers.  The classification follows as a consequence of Theorems~\ref{ ideal non degenere}, \ref{central}, and~\ref{central1}.

In order to complete the classification, we analyze separately the cases determined by the degeneracy of the product space ${\cal A} \bullet {\cal A}$. 

In the case of Novikov algebras,  a classification up to dimension 3 was offered in \cite{BEL} in the case where ${\cal A} \bullet {\cal A}$ is non-degenerate. Here we complete the classification in the case where ${\cal A} \bullet {\cal A}$ is degenerate. We cover the super case as well. 

The classification of \emph{nilpotent} Novikov algebras up to dimension 4 has been carried out in \cite{KKK} using different techniques.
\subsection{The case where \texorpdfstring{$\A\bullet\A$}{A\textperiodcentered{}A} is non-degenerate}
We have proved in Theorem \ref{ ideal non degenere} that any pseudo-Euclidean Novikov superalgebra $(\A, \bullet, \prs)$ for which $\A\bullet\A$ is non-degenerate is of Milnor type. We have also shown in Corollary \ref{nilpotent} that such $\A^-$ is nilpotent if and only if $\bullet$ is trivial. 
\sssbegin{Proposition}\label{de-23}
Let $(\A, \bullet, \prs)$ be an even  pseudo-Euclidean Novikov superalgebra of dimension
at most $3$ such that $\A\bullet\A$ is non-degenerate. We have:
\begin{enumerate}
    \item[$(1)$] If $\dim \A=2$, then $(\A,\bullet,\prs)$ is isomorphic to
    $(\A^2_1, \bullet, \prs)$, the even pseudo-Euclidean trivial superalgebra.
    \item[$(2)$] If $\dim \A=3$,  then $(\A,\bullet,\prs)$ is isomorphic to one of the following even pseudo-Euclidean Novikov superalgebras: 
    \begin{enumerate}
        \item[$(i)$] $(\A^3_1,\bullet, \prs)$, the  trivial superalgebra, i.e,. $u\bullet v=0$, for $u,v\in \A^3_1.$
        \item[$(ii)$] $(\A^3_2, \bullet, \prs)$, where the product and the form are given in the basis $\{ e_1, e_2, f_1\}$ by 
        \[
f_1\bullet e_1=\la e_2,\qquad
f_1\bullet e_2=-\epsilon\la e_1,\quad \prs_\epsilon=e_1^*\otimes e_1^*+\epsilon e_2^*\otimes e_2^*+ f_1^*\otimes f_1^*,
\]  
where $\la\neq 0$ and $\epsilon=\mp 1$.
        \item[$(iii)$]  $(\A^3_3, \bullet, \prs)$, where the product and the form are given in the basis $\{ f_1\mid e_1, e_2\}$ by 
        \[
f_1\bullet e_1=\la e_1,\qquad
f_1\bullet e_2=-\la e_2,  \quad \prs=-e_1^*\odot  e_2^*+ f_1^*\otimes f_1^*.
\]
        \item[$(iv)$] $(\A^3_4, \bullet, \prs)$, where the product and the form are given in the basis $\{ f_1\mid e_1, e_2\}$ by    \[
f_1\bullet e_1=\lambda e_2,\qquad
f_1\bullet e_2=e_1,\quad \quad \prs=- e_1^*\odot  e_2^*+ f_1^*\otimes f_1^*.
\]
where $\lambda\neq0$. 
    \end{enumerate}
\end{enumerate}
\end{Proposition}

\begin{proof}
Assume that $\A\bullet\A$ is non-degenerate. Then
\(
\A = \A\bullet\A \oplus (\A\bullet\A)^{\perp}.
\)
First, observe
that  $\A\neq \A\bullet\A$  by
Proposition~\ref{Pr1}. Moreover, if $\A=(\A\bullet\A)^{\perp}$, then, by
Eq.~\ref{ortogonal}, the Novikov superalgebra
$(\A,\bullet)$ is trivial. 

\medskip
\noindent
\textbf{The case where $\dim( \A)=2$.} Assume that $\A\bullet\A\neq 0$, in this case
\(
\dim(\A\bullet\A)=\dim((\A\bullet\A)^{\perp})=1.
\)
By Lemma~\ref{evendimension}, both subspaces  $\A\bullet\A$ and $(\A\bullet\A)^\bot$ are even.
Let $\A\bullet\A=\mathbb{R} e_1$ and $(\A\bullet\A)^{\perp}=\mathbb{R} f_1$ such that $\langle e_1, e_1\rangle=\epsilon$ where $\epsilon=\mp 1$.
By Theorem~\ref{milnor}, there exists $\alpha\in\mathbb{R}$ such that
$f_1\bullet e_1=\alpha e_1$.
However,
\[
\epsilon\alpha=\langle f_1\bullet e_1,e_1\rangle
=-\langle e_1,f_1\bullet e_1\rangle=-\epsilon \al=0,
\]
which contradicts $\A\bullet\A\neq\{0\}$. Hence $(\A,\bullet)$ is trivial.

\medskip
\noindent
\textbf{The case where $\dim (\A)=3$.}
We distinguish two subcases: either 
\(
\dim(\A\bullet\A)=2\) and \(\dim((\A\bullet\A)^{\perp})=1,\)
 or \(
\dim(\A\bullet\A)=1\) and \( \dim((\A\bullet\A)^{\perp})=2.
\)

\medskip
\noindent
\textbf{Subcase 1:} Let us assume that $\dim(\A\bullet\A)=2$ and $\dim((\A\bullet\A)^{\perp})=1$.
By Lemma~\ref{evendimension}, either $\A\bullet\A$ is purely even or purely odd, and $(\A\bullet\A)^{\perp}$ is purely even.

\smallskip
\noindent
(i) Suppose  $\A\bullet\A\subseteq\A_{\bar 0}$.
Let $\{e_1,e_2\mid 0\}$ be a basis of $\A\bullet\A$, with $\langle e_1, e_1\rangle=1$   and $\langle e_2, e_2\rangle= \epsilon$, where $\epsilon=\mp 1$  and let
\(
(\A\bullet\A)^{\perp}=\mathbb{R} f_1\)  with \( \langle f_1,f_1\rangle=1.
\)
According to Theorem~\ref{milnor}, we have 
\[
f_1\bullet e_1=\alpha_1 e_1+\beta_1 e_2,\qquad
f_1\bullet e_2=\alpha_2 e_1+\beta_2 e_2,
\]
for some $\alpha_i,\beta_i\in\mathbb{R}$.  
Since the left multiplication operators are $\prs$-antisymmetric, we obtain
\[
\alpha_1=\beta_2=0,\qquad \alpha_2=-\beta_1.
\]
Consequently, 
\[
f_1\bullet e_1=\beta_1 e_2,\qquad
f_1\bullet e_2=-\epsilon\beta_1 e_1.
\]
Since $\A\bullet\A\neq\{0\}$, we have $\beta_1\neq0$. Hence $(\A,\bullet,\prs)$ is isomorphic to
$(\A^3_2, \bullet, \prs)$.

\smallskip
\noindent
(ii) Suppose $\A\bullet\A\subseteq\A_{\bar 1}$. Let $\{0\mid e_1,e_2\}$ be a basis of $\A\bullet\A$ such that
$\langle e_1,e_2\rangle=1$, and let
\(
(\A\bullet\A)^{\perp}=\mathbb{R} f_1\)  with  \(\langle f_1,f_1\rangle=1.
\)
Then
\[
f_1\bullet e_1=\alpha_1 e_1+\beta_1 e_2,\qquad
f_1\bullet e_2=\alpha_2 e_1+\beta_2 e_2,
\]
with $\alpha_1\beta_2-\beta_1\alpha_2\neq0$. Using the antisymmetry of the left multiplication, we obtain $\beta_2=-\alpha_1$.
Thus,
\[
f_1\bullet e_1=\alpha_1 e_1+\beta_1 e_2,\qquad
f_1\bullet e_2=\alpha_2 e_1-\alpha_1 e_2,
\]
where $\alpha_1^2+\beta_1\alpha_2\neq0$.

If $\alpha_2=\beta_1=0$, then $\alpha_1\neq0$ and
\[
f_1\bullet e_1=\alpha_1 e_1,\qquad
f_1\bullet e_2=-\alpha_1 e_2.
\]
Hence $(\A,\bullet,\prs)$ is isomorphic to
$(\A^3_3, \bullet, \prs)$.

If $\alpha_2\neq0$, set
\[
x_1:=\frac{1}{\sqrt{|\alpha_2|}}(\alpha_2 e_1-\alpha_1 e_2),\quad
x_2:=\frac{1}{\sqrt{|\alpha_2|}}e_2,\quad
y_1:=f_1.
\]
Then
\[
y_1\bullet x_1=\lambda x_2,\qquad
y_1\bullet x_2=x_1,
\]
where $\lambda=\alpha_1^2+\alpha_2\beta_1\neq0$.

If \(\alpha_2>0\), then \((\mathcal A,\bullet,\prs)\) is isomorphic to \((\mathcal A^3_4,\bullet,\prs)\). If \(\alpha_2<0\), we set \(y'_1:=-y_1\) and \(x'_1:=-x_1\). In this case, \((\mathcal A,\bullet,\prs)\) is also isomorphic to \((\mathcal A^3_4,\bullet,\prs)\).

If $\beta_1\neq0$, a similar argument shows that
$(\A,\bullet,\prs)$ is isomorphic to
$(\A^3_4, \bullet, \prs)$.

\medskip
\noindent
\textbf{Subcase 2:} Let us assume that $\dim(\A \bullet \A) = 1$, $\dim((\A \bullet \A)^\bot) = 2$. By Lemma~\ref{evendimension}, we deduce that  $\A\bullet\A$ is even. Set
\(
\A\bullet\A=\mathbb{R} e_1,\) and \( (\A\bullet\A)^{\perp}=\mathbb{R} f_1\oplus \mathbb{R} f_2.
\)
According to Theorem~\ref{milnor}, there exist $\alpha, \beta\in\mathbb{R}$ such that
\[
f_1\bullet e_1=\alpha e_1, \esp f_2\bullet e_1=\beta e_1.
\]
However,
\[
\al=\langle f_1\bullet e_1,e_1\rangle
= -\langle e_1,f_1\bullet e_1\rangle=0,\esp \beta=\langle f_2\bullet e_1,e_1\rangle
= -\langle e_1,f_2\bullet e_1\rangle=0,
\]
which contradicts the assumption that $\A\bullet\A\neq\{0\}$.
 Hence $(\A,\bullet)$ is a trivial superalgebra.

\end{proof}

\sssbegin{Theorem}
Let $(\A, \bullet, \prs)$ be an even pseudo-Euclidean Novikov superalgebra of dimension
 $4$ such that $\A\bullet\A$ is non-degenerate.  Then
$(\A,\bullet,\prs)$ is isomorphic to one of the following pseudo-Euclidean Novikov superalgebras:
\begin{enumerate}
\item[$(i)$] $(\A^4_1, \bullet, \prs)$, the trivial superalgebra, i.e,. $u\bullet v=0$, for all $u,v\in \A^4_1.$
        \item[$(ii)$] $(\A^4_2, \bullet, \prs_\varepsilon^1 \text{ or } \prs_\varepsilon^2  )$, where the product is given in the basis $\{e_1, e_2,  f_1, f_2\}$ by 
        \begin{align*}
&f_1\bullet e_1=\la e_2,\qquad
f_1\bullet e_2=-\varepsilon\la e_1,
\end{align*}
and the bilinear forms are 
\begin{align*}
&\prs^1_{\varepsilon, \alpha} =e_1^*\otimes e_1^*+\varepsilon e_2^*\otimes e_2^*+ f_1^*\otimes f_1^* -\al f_1^*\odot f_2^*+ (\al^2-1)f_2^*\otimes f_2^*, \\ &\prs^2_\varepsilon =e_1^*\otimes e_1^*+\varepsilon e_2^*\otimes e_2^*- f_1^*\otimes f_1^* + f_2^*\otimes f_2^*,
\end{align*} 
where $0 \neq \la$, $\al\in \mathbb{R}$,  and $\varepsilon=\mp 1$. 

\item[$(iii)$] $(\A^4_3, \bullet, \prs)$, where the product is given in the basis $\{f_1, f_2\mid e_1, e_2\}$ by 
        \begin{align*}
&f_1\bullet e_1=\la e_2,\qquad
f_1\bullet e_2=-\varepsilon\la e_1,
\end{align*}
and the bilinear forms are 
\begin{align*}
\prs_\alpha& =e_1^*\odot e_2^*+ f_1^*\otimes f_1^* -\al f_1^*\odot f_2^*+ (\al^2-1)f_2^*\otimes f_2^*,
\end{align*} 
where $0 \neq \la$, $\al\in \mathbb{R}$,  and $\varepsilon=\mp 1$.

        \item[$(iv)$] $(\A^4_4, \bullet, \prs)$, where the product is given in the basis $\{ f_1, f_2\mid e_1, e_2\}$ by 
\[
f_1\bullet e_1=\lambda e_1,\qquad
f_1\bullet e_2=-\la e_2,\] 
and the bilinear forms are given by 
\begin{align*}
\prs_\epsilon^1=&-e_1^*\odot  e_2^*+ \epsilon f_1^*\otimes f_1^*+ f_2^*\otimes f_2^*, \\ \prs_\epsilon^2= &-e_1^*\odot  e_2^*+  f_1^*\otimes f_1^*+ \epsilon f_2^*\otimes f_2^*,\\ \prs^3_\varepsilon =&e_1^*\otimes e_1^*+\varepsilon e_2^*\otimes e_2^*+ f_1^*\otimes f_1^*+ f_2^*\otimes f_2^*,
\end{align*}
where $\lambda\neq0$,  and $\epsilon=\mp 1$.
 \item[$(v)$] $(\A^4_5, \bullet, \prs)$, where the product is given in the basis $\{ f_1, f_2\mid e_1, e_2\}$ by 
  \begin{align*}
&f_1\bullet e_1=\lambda e_2,\;
f_1\bullet e_2=e_1,
\end{align*}
and the bilinear forms are given by 
\begin{align*}
\prs_{\epsilon}^1 & =-e_1^*\odot  e_2^*+\epsilon t f_1^*\otimes f_1^*+ f_2^*\otimes f_2^*,\\\
\prs_{\epsilon}^2&=-e_1^*\odot  e_2^*+  f_1^*\otimes f_1^*+ \epsilon f_2^*\otimes f_2^*, \\
\prs_\ga & =-e_1^*\odot e_2^*+  f_1^*\otimes f_1^*-\ga f_1^*\odot f_2^*+(\ga^2+1)f_2^*\otimes f_2^*,
\end{align*}
where $\lambda\neq0$, $t>0$, $\ga\in \mathbb{R}$,  and $\epsilon=\mp 1$.
\end{enumerate}
\end{Theorem}
\begin{proof} Assume that $\A\bullet\A$ is non-degenerate. Then
\(
\A = \A\bullet\A \oplus (\A\bullet\A)^{\perp}.
\)
First, observe
that  $\A\neq \A\bullet\A$  by
Proposition~\ref{Pr1}. Moreover, if $\A=(\A\bullet\A)^{\perp}$, then, by
Eq.~\ref{ortogonal}, the Novikov superalgebra
$(\A,\bullet)$ is trivial. 

We distinguish several cases according to the dimension and the parity of
$\A\bullet\A$ and $(\A\bullet\A)^{\perp}$.

\medskip
\noindent\textbf{Case 1.} Assume that $\dim(\A\bullet\A) = 1$ and $\dim(\A\bullet\A)^\perp = 3$. By Lemma~\ref{evendimension}, the superspace $\A\bullet\A$ must be purely even. Let us set $\A\bullet\A := \mathbb{R} e$ with $\langle e, e \rangle = 1$. Then, for any $x \in \A$, we have
$
x \bullet e = \alpha(x)\, e,
$
where $\alpha \in \A^*$. Since the left multiplication operators are $\prs$-antisymmetric, it follows that
\[
\alpha(x) = \langle x \bullet e, e \rangle = - \langle e, x \bullet e \rangle = -\alpha(x),
\]
which implies $\alpha(x) = 0$ for all $x \in \A$. Consequently, $\A \bullet \A = \{0\}$, contradicting the assumption. Hence, this case is impossible.

\medskip
\noindent\textbf{Case 2.} Assume $\dim(\A\bullet\A)=2$ and $\dim(\A\bullet\A)^{\perp}=2$. By Lemma~\ref{evendimension}, these two spaces  are either purely even or purely odd. 

\medskip
\noindent\textbf{Subcase 2.1.} Suppose that both $\A\bullet\A$ and $(\A\bullet\A)^{\perp}$ are purely  even. Let $\{e_1,e_2\}$ be a basis of $\A\bullet\A$ such that
\(
\langle e_1,e_1\rangle=1,\) and \( \langle e_2,e_2\rangle=\varepsilon,\) where
\( \varepsilon=\pm 1\) 
and let $\{f_1,f_2\}$ be a basis of $(\A\bullet\A)^{\perp}$ satisfying
\(
\langle f_1,f_1\rangle=1,\) and \( \langle f_2,f_2\rangle=\epsilon,\) where
\( \epsilon= \pm 1\). 
According to Theorem~\ref{milnor}, and using the $\prs$-antisymmetry of the left
multiplication, we obtain
\begin{align*}
f_1\bullet e_1&=\alpha e_2, &f_1\bullet e_2&=-\varepsilon\alpha e_1,& 
f_2\bullet e_1&=\beta e_2, &f_2\bullet e_2&=-\varepsilon\beta e_1,
\end{align*}
with $(\alpha,\beta)\neq(0,0)$.

If $\epsilon=1$, set
\[
(x_1,x_2,y_1,y_2):=
\Bigl(e_1,e_2,
\frac{\alpha f_1+\beta f_2}{\sqrt{\alpha^2+\beta^2}},
\frac{\beta f_1-\alpha f_2}{\sqrt{\alpha^2+\beta^2}}\Bigr).
\]
Then
\[
y_1\bullet x_1=\lambda x_2,\qquad
y_1\bullet x_2=-\varepsilon\lambda x_1,
\]
 and
\[
\langle x_1,x_1\rangle=1,\; \langle x_2,x_2\rangle=\varepsilon,\; \langle y_1,y_1\rangle= \langle y_2,y_2\rangle=1.
\]
where $\lambda=\sqrt{\alpha^2+\beta^2}$. Hence $(\A,\bullet,\prs)$ is isomorphic to $(\A^4_4, \bullet,  \prs^3_\varepsilon)$.

If $\epsilon=-1$ and $\alpha\neq0$, set
\[
(x_1,x_2,y_1,y_2):=
\bigl(e_1,e_2,f_1,f_2-\tfrac{\beta}{\alpha}f_1\bigr).
\]
Therefore,
\[
y_1\bullet x_1=\al x_2,\qquad
y_1\bullet x_2=-\varepsilon\al x_1,
\]
 and
\[
\langle x_1,x_1\rangle=1,\; \langle x_2,x_2\rangle=\varepsilon,\;\langle y_1, y_1\rangle=1,\; \langle y_2,y_2\rangle=\Bigl(\frac{\beta}{\alpha}\Bigr)^2-1,
\quad
\langle y_1,y_2\rangle=-\frac{\beta}{\alpha}.
\]
It follows that   $(\A,\bullet,\prs)$ is isomorphic to $(\A^4_2, \bullet  \prs^1_\varepsilon)$ .

Now, if $\alpha=0$,  then $(\A,\bullet,\prs)$ is isomorphic to $(\A^4_2, \bullet, \prs^2_\varepsilon)$.

\medskip
\noindent\textbf{Subcase 2.2.} Assume that $\A\bullet\A$ is even and $(\A\bullet\A)^{\perp}$ is odd.  Since $\Ll_u=0$ for all $u\in \A\bullet \A$, we have
$
(\A\bullet \A)^\perp \bullet (\A\bullet \A)
\subseteq
(\A\bullet \A)\cap (\A\bullet \A)^\perp.
$
Hence $(\A\bullet \A)^\perp \bullet (\A\bullet \A)=\{0\}$, which is impossible. Therefore, this case cannot occur.

\medskip
\noindent\textbf{Subcase 2.3.} Assume that $\A\bullet\A$ is odd and $(\A\bullet\A)^{\perp}$ is even. Let $\{f_1,f_2\}$ be a basis of $(\A\bullet\A)^\bot$ such that
\(
\langle f_1,f_1\rangle=1,\) and \(\langle f_2,f_2\rangle=\varepsilon,\)
where $\varepsilon= \pm 1.$ We first analyze the possible values of $\dim(\operatorname{Im}\Ll_{f_1})$.
Assume that $\dim(\operatorname{Im}\Ll_{f_1})=1$.
According to Theorem~\ref{milnor}, there exists $e_1\in \A\bullet\A$
such that
$
f_1\bullet x=\alpha(x)e_1,$ for all $x\in\A\bullet\A.$ Since the restriction of $\prs$ to $\A\bullet\A$ is non-degenerate,
there exists $e_2\in \A\bullet\A$ such that
\(
\langle e_1,e_2\rangle=1.
\)
Thus, we may write
\[
f_1\bullet e_1=\alpha e_1,\qquad
f_1\bullet e_2=\beta e_1,
\]
with $(\alpha,\beta)\neq(0,0)$.

Since left multiplications are $\prs$-antisymmetric, we obtain
\[
\alpha=\langle f_1\bullet e_1,e_2\rangle
=-\langle e_1,f_1\bullet e_2\rangle=0.
\]
Hence $\alpha=0$.

Now let
\[
f_2\bullet e_1=\alpha_1 e_1+\beta_1 e_2,\qquad
f_2\bullet e_2=\alpha_2 e_1-\alpha_1 e_2.
\]
Since $\Ll_{f_1}\circ\Ll_{f_2}=\Ll_{f_2}\circ\Ll_{f_1}$, we deduce that
$
\alpha_1=\beta_1=0,
$
and therefore
\[
f_2\bullet e_2=\alpha_2 e_1.
\]
This implies that $\dim(\A\bullet\A)=1$, which contradicts the assumption
$\dim(\A\bullet\A)=2$.
Consequently,
\[
\dim(\operatorname{Im}\Ll_{f_i})=0 \text{ or } 2.
\]

\medskip
Assume now that $\Ll_{f_1}=0$.
Let $\{e_1,e_2\}$ be a basis of $\A\bullet\A$ such that
$
\langle e_1,e_2\rangle=1$. 
According to Theorem~\ref{milnor}, we have
\[
f_2\bullet e_1=\alpha_1 e_1+\beta_1 e_2,\qquad
f_2\bullet e_2=\alpha_2 e_1+\beta_2 e_2,
\]
with $\alpha_1\beta_2-\beta_1\alpha_2\neq0$.
Using the antisymmetry of the left multiplication, we obtain $\beta_2=-\alpha_1.$ Hence
\[
f_2\bullet e_1=\alpha_1 e_1+\beta_1 e_2,\qquad
f_2\bullet e_2=\alpha_2 e_1-\alpha_1 e_2,
\]
where $\alpha_1^2+\beta_1\alpha_2\neq0$.

If $\alpha_2=\beta_1=0$, then $\alpha_1\neq0$ and
\[
f_2\bullet e_1=\alpha_1 e_1,\qquad
f_2\bullet e_2=-\alpha_1 e_2.
\]
In this case, $(\A,\bullet,\prs)$ is isomorphic to $(\A^4_4, \bullet, \prs_\epsilon^1)$.

If $\alpha_2\neq0$, set
\[
x_1:=\frac{1}{\sqrt{|\alpha_2|}}(\alpha_2 e_1-\alpha_1 e_2),\qquad
x_2:=\frac{1}{\sqrt{|\alpha_2|}}e_2.
\]
Then
\[
f_2\bullet x_1=\lambda x_2,\qquad
f_2\bullet x_2=x_1,
\]
where $\lambda=\alpha_1^2+\alpha_2\beta_1\neq0$.

If \(\alpha_2>0\), then \((\mathcal A,\bullet,\prs)\) is isomorphic to $(\A^4_5,\bullet, \prs_{\epsilon}^1)$. If \(\alpha_2<0\), we set \(f'_2:=-f_2\) and \(x'_1:=-x_1\). In this case, \((\mathcal A,\bullet,\prs)\) is also isomorphic to $(\A^4_5,\bullet, \prs_{\epsilon}^1)$.

Thus, $(\A,\bullet,\prs)$ is isomorphic to $(\A^4_5,\bullet, \prs_{\epsilon}^1)$.
If $\beta_1\neq0$, a similar argument leads to the same conclusion.

\medskip
If $\Ll_{f_2}=0$, the same reasoning shows that
$(\A,\bullet,\prs)$ is isomorphic to $(\A^4_4, \bullet, \prs_\epsilon^2)$ or to $(\A^4_5, \bullet, \prs_{\epsilon}^2)$.

\medskip
 Finally, assume that $\Ll_{f_1}\neq0$. Then
$\dim(\operatorname{Im}\Ll_{f_1})=2$.
Let $\{e_1,e_2\}$ be a basis of $\A\bullet\A$ with
$\langle e_1,e_2\rangle=1$.
By Theorem~\ref{milnor}, we have
\[
f_1\bullet e_1=\alpha_1 e_1+\beta_1 e_2,\qquad
f_1\bullet e_2=\alpha_2 e_1+\beta_2 e_2,
\]
with $\alpha_1\beta_2-\beta_1\alpha_2\neq0$.
The antisymmetry condition yields $\beta_2=-\alpha_1$, and hence
\[
f_1\bullet e_1=\alpha_1 e_1+\beta_1 e_2,\qquad
f_1\bullet e_2=\alpha_2 e_1-\alpha_1 e_2,
\]
where $\alpha_1^2+\beta_1\alpha_2\neq0$.

If $\alpha_2=\beta_1=0$, then $\alpha_1\neq0$ and
\[
f_1\bullet e_1=\alpha_1 e_1,\qquad
f_1\bullet e_2=-\alpha_1 e_2.
\]
Since $\Ll_{f_1}\circ\Ll_{f_2}=\Ll_{f_2}\circ\Ll_{f_1}$, and $\Ll_{f_2}$ is $\prs$-antisymmetric, we deduce that 
\[
f_2\bullet e_1=\ga e_1,\qquad
f_2\bullet e_2=-\ga e_2.
\]

If $\epsilon=1$, set
\[
(y_1,y_2, x_1, x_2):=
\Bigl(
\frac{\alpha_1 f_1+\ga f_2}{\sqrt{\alpha_1^2+\ga^2}},
\frac{\ga f_1-\alpha_1 f_2}{\sqrt{\alpha_1^2+\ga^2}}, e_1,e_2\Bigr).
\]
Then
\[
y_1\bullet x_1=\lambda x_1,\qquad
y_1\bullet x_2=-\lambda x_2,
\]
 and
\[
\langle x_1,x_2\rangle=1,\quad  \langle y_1,y_1\rangle= \langle y_2,y_2\rangle=1.
\]
where $\lambda=\sqrt{\alpha_1^2+\ga^2}$. Hence $(\A,\bullet,\prs)$ is isomorphic to $(\A^4_4, \bullet, \prs^1_{\epsilon=1})$.

If $\epsilon=-1$, set
\[
(y_1,y_2, x_1,x_2):=
\bigl(f_1,f_2-\tfrac{\ga}{\alpha_1}f_1, e_1,e_2\bigr).
\]
Then
\[
y_1\bullet x_1=\al_1 x_2,\qquad
y_1\bullet x_2=-\al_1 x_1,
\]
 and
\[
\langle x_1,x_2\rangle=1,\;\langle y_1, y_1\rangle=1,\; \langle y_2,y_2\rangle=\Bigl(\frac{\ga}{\alpha_1}\Bigr)^2-1,
\quad
\langle y_1,y_2\rangle=-\frac{\ga}{\alpha_1},
\]
Thus, we deduce that   $(\A,\bullet,\prs)$ is isomorphic to $(\A^4_3,\bullet,  \prs)$ .

 If $\alpha_2\neq0$, set
\[
x_1:=\frac{1}{\sqrt{|\alpha_2|}}(\alpha_2 e_1-\alpha_1 e_2),\quad
x_2:=\frac{1}{\sqrt{|\alpha_2|}}e_2.
\]
Then
\[
f_1\bullet x_1=\lambda x_2,\qquad
f_1\bullet x_2=x_1, \text{ and } \langle x_1, x_2\rangle =s\]
where $\lambda=\alpha_1^2+\alpha_2\beta_1\neq0$ and $s\in \{-1, 1\}$.

 Since $\Ll_{f_1}\circ\Ll_{f_2}=\Ll_{f_2}\circ\Ll_{f_1}$, and $\Ll_{f_2}$ is $\prs$-antisymmetric, we deduce that 
\[
f_2\bullet x_1=\la\ga x_2,\qquad
f_2\bullet x_2=\ga x_1.
\]

If $\epsilon=1$, set
\[
(y_1,y_2, x_1',x_2'):=
\Bigl(
 f_1+\ga f_2,
\frac{ f_2-\ga f_1}{\sqrt{1+\ga^2}}, x_1,x_2\Bigr).
\]
Then
\[
y_1\bullet x_1'=\ \lambda x_2',\qquad
y_1\bullet x_2'= x_1',
\]
 and
\[
\langle x_1',x_2'\rangle=s,\quad  \langle y_1,y_1\rangle=\lambda_1,\; \langle y_2,y_2\rangle=1.
\]
where $\lambda_1=\sqrt{1+\ga^2}$ and $s\in\{-1, 1\}$. 

If $s=1$, then  $(\A,\bullet,\prs)$ is isomorphic to $(\A^4_5, \bullet,  \prs^1)$. If $s=-1$, we set \(y_1^{'}:=-y_1\) and \(x^{''}_1:=-x_1\). In this case, \((\mathcal A,\bullet,\prs)\) is also isomorphic to $(\A^4_5,\bullet, \prs_{\epsilon}^1)$.

If $\epsilon=-1$, set
\[
(y_1,y_2, x_1',x_2')=
\bigl(f_1,f_2-\ga f_1, x_1,x_2\bigr).
\]
Then
\[
y_1\bullet x_1'=\la x_2',\qquad
y_1\bullet x_2'= x_1',
\]
 and 
\[
\langle x_1',x_2'\rangle=s,\quad \langle y_1, y_1\rangle=1,\quad  \langle y_2,y_2\rangle=\ga^2-1,
\quad
\langle y_1,y_2\rangle=-\ga,
\]
If $s=1$, then  $(\A,\bullet,\prs)$ is isomorphic to $(\A^4_5,\bullet,  \prs_\ga)$. If \(s=-1\), we set \(y_1^{'}:=-y_1\) and \(x^{''}_1:=-x_1\). In this case, \((\mathcal A,\bullet,\prs)\) is also isomorphic to $(\A^4_5,\bullet,  \prs_\ga)$. 

If $\beta_1\neq0$, a similar argument leads to the same conclusion.

\medskip
\noindent\textbf{Case 3.}  Assume that $\dim(\mathcal{A}\bullet\mathcal{A}) = 3$ and $\dim(\mathcal{A}\bullet\mathcal{A})^{\perp} = 1$. By Lemma~\ref{evendimension}, the superspace $(\mathcal{A}\bullet\mathcal{A})^{\perp}$ is purely even. Consequently, we distinguish the following cases: either $\mathcal{A}\bullet\mathcal{A}$ is purely even, or $\dim(\mathcal{A}_0 \cap (\mathcal{A}\bullet\mathcal{A})) = 1$ and $\dim(\mathcal{A}_1 \cap (\mathcal{A}\bullet\mathcal{A})) = 2$.

We first consider the case where $\mathcal{A}\bullet\mathcal{A}$ is purely even. Choose a basis $\{e_1,e_2,e_3\}$ of $\mathcal{A}\bullet\mathcal{A}$ such that
\[
\langle e_1,e_1\rangle = \langle e_2,e_2\rangle = 1, \qquad \langle e_3,e_3\rangle = \epsilon,
\]
where $\epsilon = \pm 1$, and let $(\mathcal{A}\bullet\mathcal{A})^{\perp} = \mathbb{R} f_1$ with $\langle f_1,f_1\rangle = 1$.

Since the left multiplications are $\prs$-antisymmetric, we obtain
\[
f_1\bullet e_1 = \beta_1 e_2 + \lambda_1 e_3, \quad
f_1\bullet e_2 = -\beta_1 e_1 + \lambda_3 e_3, \quad
f_1\bullet e_3 = -\epsilon \lambda_1 e_1 - \epsilon \lambda_3 e_2.
\]
Hence, the family
\[
\{\beta_1 e_2 + \lambda_1 e_3,\; -\beta_1 e_1 + \lambda_3 e_3,\; -\epsilon \lambda_1 e_1 - \epsilon \lambda_3 e_2\}
\]
is linearly dependent, and thus cannot form a basis of $\mathcal{A}\bullet\mathcal{A}$, which contradicts the fact that $\dim(\mathcal{A}\bullet\mathcal{A})=3$. Therefore, this case is impossible.

Now, assume that $\dim(\mathcal{A}_0 \cap (\mathcal{A}\bullet\mathcal{A})) = 1$ and $\dim(\mathcal{A}_1 \cap (\mathcal{A}\bullet\mathcal{A})) = 2$. Proceeding in a similar way, we obtain again a contradiction. Hence, this case is also impossible.
\end{proof}

We now study the case where the form $\prs$ is odd. 
\sssbegin{Theorem}\label{de-24}
Let $(\A,\bullet,\prs)$ be an odd pseudo-Euclidean Novikov superalgebra of
dimension at most $4$ such that $\A\bullet\A$ is non-degenerate.
We have:
\begin{enumerate}
\item[$(1)$] If $\dim\A=2$, then $(\A,\bullet,\prs)$ is isomorphic to
$(\A_1^2, \bullet, \prs)$, the odd pseudo-Euclidean trivial superalgebra.

\item[$(2)$] If $\dim\A=4$, then $(\A,\bullet,\prs)$ is isomorphic to one of the
following pseudo-Euclidean Novikov superalgebras:
\begin{enumerate}
\item[$(i)$] the trivial superalgebra;
\item[$(ii)$] $(\A_6^4,\bullet,\prs)$ with homogeneous basis
$\{e_1,f_1\mid e_2,f_2\}$, whose non-zero product and bilinear form are given by
\[
f_1\bullet e_1=\lambda e_1,\qquad
f_1\bullet e_2=-\lambda e_2,
\qquad
\prs=e_1^*\odot e_2^*+f_1^*\odot f_2^*,
\]
where $\lambda\neq0$. 
\end{enumerate}
\end{enumerate}
\end{Theorem}
\begin{proof}
Assume that $\A\bullet\A$ is non-degenerate. It follows that 
\(
\A=\A\bullet\A\oplus(\A\bullet\A)^{\perp}.
\)
By Lemma~\ref{evendimension}, the superalgebra $\A$ is even-dimensional, and
both $\A\bullet\A$ and $(\A\bullet\A)^{\perp}$ are even-dimensional superspaces.

If $\dim\A=2$, then necessarily $\A=\A\bullet\A$ or $\A=(\A\bullet\A)^\bot$, and it is easy to see that
$(\A,\bullet)$ is the trivial superalgebra. 

Assume now that $\dim\A=4$.
Let $\{e_1\mid e_2\}$ be a homogeneous basis of $\A\bullet\A$ such that
\(
\langle e_1,e_2\rangle=1,
\)
and let $\{f_1\mid f_2\}$ be a homogeneous basis of $(\A\bullet\A)^{\perp}$
such that 
\(
\langle f_1,f_2\rangle=1.
\)
According to Theorem~\ref{milnor}, the multiplication is given by
\[
f_1\bullet e_1=\alpha e_1,\qquad
f_1\bullet e_2=\beta e_2,\qquad
f_2\bullet e_1=\alpha_1 e_2,\qquad
f_2\bullet e_2=\beta_1 e_1.
\]
Since the left multiplication operators are $\prs$-antisymmetric, we obtain
\[
\alpha=-\beta,\qquad \alpha_1=\beta_1=0.
\]
Moreover, since $\A\bullet\A\neq\{0\}$, we have $\alpha\neq0$.
Therefore, $(\A,\bullet,\prs)$ is isomorphic to $(\A_6^4,\bullet,\prs)$.
This completes the proof.
\end{proof}

\subsection{The case where \texorpdfstring{$\A\bullet\A$}{bullet} is degenerate}

Here, we classify pseudo-Euclidean Novikov superalgebras of total dimension at most $4$ such that $\A\bullet\A$ is degenerate.  
Theorems~\ref{central} and Theorem ~\ref{central1} imply that every four-dimensional pseudo-Euclidean Novikov superalgebra $(\A,\bullet, \prs)$ is an odd or even  double extension of a two-dimensional pseudo-Euclidean Novikov superalgebra $(\mathcal{B},\bullet_\mathcal{B}, \prs_\mathcal{B})$, defined by  an admissible triple $(D,\xi,b_0)$ or an admissible quadruple
$(D,\xi,b_0,c_0)$.

\sssbegin{Proposition}\label{Dim2}
\begin{enumerate}
    \item[$(i)$] There is no two-dimensional even pseudo-Euclidean Novikov superalgebra $(\A,\bullet,\prs)$ such that $\A\bullet \A$ is degenerate.
    \item[$(ii)$] Every two-dimensional odd pseudo-Euclidean Novikov superalgebra $(\A,\bullet,\prs)$ such that $\A\bullet \A$ is degenerate  is isomorphic to the superalgebra $(\A_2^2,\bullet,\prs)$ with homogeneous basis
$\{e_1\mid f_1\}$, whose non-zero product and bilinear form are given by
\[
f_1\bullet f_1=\al e_1,
\qquad
\prs=e_1^*\odot f_1^*,
\]
where $\al\neq0$.
    \end{enumerate}
\end{Proposition}

\begin{proof}
Part $(i)$. Assume that $\A \bullet \A$ is degenerate. Then there exists a nonzero element
$e \in \A \bullet \A \cap (\A \bullet \A)^{\perp}$. By Proposition~\ref{Pr1} and
Eq.~\eqref{ortogonal}, we have
\(
\Ll_e = \Rr_e = 0.
\) Since $\prs$ is even, Lemma~\ref{evendimension}
implies that $\A$ is either purely even or purely odd.  If $\A$ is purely odd, then $(\A,\bullet)$ is clearly trivial.  
If $\A$ is purely even, there exists $d \in \A$ such that
$\langle e,d\rangle = 1$. Since $\dim \A = 2$, we have  $\A\bullet\A=\mathbb{R} e$, hence we may write
\(
d \bullet d = \alpha e,\) where $\alpha \in \mathbb{R}.$
However,
\[
\langle d \bullet d, d\rangle = - \langle d, d \bullet d\rangle =\al= 0
\quad \text{and} \quad
\langle d \bullet d, e\rangle = 0,
\]
which forces $\alpha = 0$. This contradicts the fact that  $\A\bullet\A$ is degenerate. Hence $(\A,\bullet)$ is trivial.

Part $(ii)$. Assume that $\A \bullet \A$ is degenerate. Then there exists a nonzero element
$e \in \A \bullet \A \cap (\A \bullet \A)^{\perp}$. By Proposition~\ref{Pr1} and
Eq.~\eqref{ortogonal}, we have
\(
\Ll_e = \Rr_e = 0.
\) Since $\prs$ is odd, Lemma~\ref{evendimension} implies that 
$\A_{\bar{0}}\neq 0$ and $\A_{\bar{1}}\neq 0$. Since $\A\neq \A\bullet\A$, we have
$
\A\bullet\A=\mathbb{R} e .
$
Hence $e$ is even. Since $\prs$ is odd and nondegenerate, we may choose 
$d\in \A_{\bar{1}}$ such that
$
\langle e,d\rangle =1 .
$
Consequently, we may write
\[
d\bullet d=\alpha e, \quad \alpha\in\mathbb{R} .
\]
Since $\A\bullet\A\neq 0$, then $\al\neq 0$.
\end{proof}

\sssbegin{Definition}
Let $(\mathcal{B}, \bullet_{\mathcal{B}}, \prs_{\mathcal{B}})$ be a pseudo-Euclidean Novikov superalgebra.

\begin{itemize}
\item If the form $\prs_{\mathcal{B}}$ is even, a triple
\(
(D,\xi,b_0)\in \End(\mathcal{B})\times \End(\mathcal{B}) \times \mathcal{B}
\)
is called \emph{even admissible} if it satisfies System~\eqref{eq:claim2}.  
A quadruple
\(
(D,\xi,b_0,c_0)\in \End(\mathcal{B})\times \End(\mathcal{B}) \times \mathcal{B}\times \mathcal{B}
\)
is called \emph{odd admissible} if it satisfies System~\eqref{eq:claim4}.

\item If the form $\prs_{\mathcal{B}}$ is odd, a triple
\(
(D,\xi,b_0)\in \End(\mathcal{B})\times \End(\mathcal{B}) \times \mathcal{B}
\)
is called \emph{even admissible} if it satisfies System~\eqref{eq:claim2}.  
A quadruple
\(
(D,\xi,b_0,c_0)\in \End(\mathcal{B})\times \End(\mathcal{B}) \times \mathcal{B}\times \mathcal{B}
\)
is called \emph{odd admissible} if it satisfies System~\eqref{eq:claim6}.
\end{itemize}
\end{Definition}
\sssbegin{Proposition}\label{solutioneven}
Let $(\mathcal{B},\prs)$ be a two-dimensional even pseudo-Euclidean trivial Novikov superalgebra  \textup{(}i.e., with product $\mathcal{B} \bullet \mathcal{B} = \{0\}$\textup{)}.
\begin{itemize}
    \item[$(i)$] The triple $(D, \xi, b_0)$ is \emph{even admissible} if and only if  there exists a basis $\mathbb{B} = \{e_1, e_2\}$ of $\mathcal{B}$ in which both elements are even \textup{(}resp. both odd\textup{)}, such that  the following conditions are satisfied:
    \begin{itemize}
        \item If $\{e_1, e_2\}$ are both even, then we either have \textup{(}where $\alpha, \beta, a \in \mathbb{R}$ and  $\epsilon=\mp 1$\textup{)}
\[
\xi = 0, \; D = 
        \begin{pmatrix}
        0 & a \\
       - \epsilon a & 0
        \end{pmatrix}, \;
        b_0 = \alpha e_1+\beta e_2, \;
        \prs = e_1^* \otimes e_1^*+\epsilon e_2^* \otimes e_2^*,
\]
or, 
  \begin{align*}
       & \xi = 
        \begin{pmatrix}
        0 & a \\
        0 & 0
        \end{pmatrix}, \quad D=0, \quad 
        b_0 = \alpha e_1, \; \text{ with one of the forms}
        \end{align*}
        \begin{align*}
        \prs &=e_1^* \otimes e_1^*+ e_2^* \otimes e_2^*,& \prs&=-e_{1}^{*}\otimes e_{1}^{*}+e_{2}^{*}\otimes e_{2}^{*},\\
 		\prs&=e_{1}^{*}\otimes e_{1}^{*}-e_{2}^{*}\otimes e_{2}^{*}, &  \prs&=e_{1}^{*}\odot e_{2}^{*}.
        \end{align*}
       where $\alpha \in \mathbb{R}$ and $a \neq 0$.\\
 \item If $\{e_1, e_2\}$ are both odd, then we either have
\[
\xi = 0, \quad  D = 
        \begin{pmatrix}
        a & b \\
        c & -a
        \end{pmatrix},\quad b_0 =0, \quad \prs = -e_1^*\odot e_2^*, \text{ where } a, b, c\in \mathbb{R},
\]
or
        \[
        \xi = 
        \begin{pmatrix}
        0 & a \\
        0 & 0
        \end{pmatrix}, \quad D = 
        \begin{pmatrix}
        0 & b \\
        0 & 0
        \end{pmatrix},\quad
        b_0 = 0, \quad
        \prs = -e_1^* \odot e_2^*, \text{ where  $a \neq 0$, $b\in \mathbb{R}$.}
        \]

    \end{itemize}

    \item[$(ii)$] The quadruple  $(D,\xi, b_0, c_0)$ is \emph{odd admissible} if and only if there exists a basis $\mathbb{B} = \{e_1, e_2\}$ of $\mathcal{B}$ in which both elements are even \textup{(}resp. both odd\textup{)}, such that \\
\begin{itemize} \item If $\{e_1, e_2\}$ are both even, we have
\[
\xi =D= 0, \quad b_0=0,\quad c_0 = \alpha e_1 + \beta e_2, \quad \prs = e_1^* \otimes e_1^*+\epsilon e_2^*\otimes e_2^*, \text{where } \epsilon=\mp1 \quad \alpha, \beta \in \mathbb{R}.
\]
Or,
\[
\xi =D= 0, \quad b_0=\la e_1,\quad c_0 = \alpha e_1, \quad \prs = e_1^* \odot e_2^*, \quad \la\neq 0 \text{ and } \alpha \in \mathbb{R}.
\]
\item If $\{e_1, e_2\}$ are both odd, we have
\[
\xi=D = 0, \quad b_0=c_0 =0, \quad \prs = -e_1^*\odot e_2^*.
\]
\end{itemize}       
\end{itemize}
\end{Proposition} 

\begin{proof}
Let us first prove Part (i). If the triple $(D, \xi, b_0)$ is even admissible,  then 
\[
\xi^2= D\circ \xi = \xi \circ D =0,  \qquad b_0 \in \ker( \xi) \esp D \text{ is $\prs$-antisymmetric}.
\]
Hence, $\xi$ is nilpotent, and therefore $\xi = 0$ or $\xi^2 = 0$ with $\xi \neq 0$.  
Since $\prs$ is even and $\mathcal{B}$ is two-dimensional, Lemma \ref{evendimension} implies that we must have $\mathcal{B} = \mathcal{B}_{\bar 0}$ or $\mathcal{B} =\mathcal{B}_{\bar 1}$.

\medskip
\noindent
\underline{The case where $\mathcal{B} = \mathcal{B}_{\bar 0}$, and $\xi=0$.}  Any $(D,b_0)$ with $D$ antisymmetric gives a solution. This gives the first case when the form $\prs$ is Euclidean or Lorentzian.

\medskip
\noindent
\underline{The case where $\mathcal{B} = \mathcal{B}_{\bar 0}$, and $\xi\neq 0$.} 
	Since $\xi^2=0$ then  $\ker \xi=\mathrm{ Im }\, \xi$ and $\dim\ker \xi=1$. Moreover, since $\xi\circ D=D\circ \xi=0$ then $D$ leaves $\ker \xi$ invariant. 
	If the form $\prs$ is Euclidean, then $D=0$. If form $\prs$ is Lorentzian, we distinguish two cases:
	
 $\bullet$ $\ker \xi=\mathbb{R} e_1$ where $\langle e_1,e_1\rangle=\pm 1$  and in this case $D=0$.
	We complete by a vector $e_2$ such that $\langle e_2,e_1\rangle=0$, $\langle e_2,e_2\rangle=\mp 1$; hence, $\xi(e_2)=be_1$ where $b\neq0$. Since $\xi(b_0)=0$, it follows that $b_0=\al e_1$, where $\al\in \mathbb{R}$.
	
	$\bullet$ $\ker \xi=\mathbb{R} e_1$ where $\langle e_1,e_1\rangle=0$ and $D(e_1)=ae_1$. Choose an isotropic vector $e_2$ such that $\langle e_1,e_2\rangle=1$. Then $D(e_2)=-ae_2$ and $\xi(e_2)=be_1$ with $b\not=0$. We have
	\[ D\circ \xi(e_2)=ab e_1=0, \] and therefore $D=0$. Since $\xi(b_0)=0$, it follows that $b_0=\al e_1$, where $\al\in \mathbb{R}.$

\medskip
\noindent
\underline{The case where $\mathcal{B} = \mathcal{B}_{\bar 1}$ and $\xi=0$.}  Any $(D,b_0)$ with $D$ antisymmetric gives a solution. We can choose a basis $\{e_1, e_2\}$ of $\mathcal{B}$ such that $\prs=-e_1^*\odot e_2^*$, then  we have $b_0=0$ and $D(e_1)=a e_1+b e_2$ and $D(e_2)= c e_1-ae_2$ where $a, b,c \in \mathbb{R}$.

\medskip
\noindent
\underline{The case where $\mathcal{B} = \mathcal{B}_{\bar 1}$, and $\xi\neq 0$.} Since $\xi^2=0$ then  $\ker \xi=\mathrm{ Im}\, \xi$ and $\dim\ker \xi=1$. We can choose a basis $\{e_1, e_2\}$ of $\mathcal{B}$ such that $\prs= -e_1\odot e_2$, and 
$$\xi = 
        \begin{pmatrix}
        0 & a \\
        0 & 0
        \end{pmatrix},$$
where $a\neq 0$. Since $D$ is $\prs$-antisymmetric, we have 
\[
D=\begin{pmatrix} b&c\\ f&-b\end{pmatrix},
\qquad b,c,f\in\mathbb{R}.
\] 
The condition  $\xi\circ D=D\circ \xi=0$ yields $b=f=0$. Finally, since $b_0$ is even, then $b_0 = 0$.

Let us prove Part (ii). If the quadruple $(D, \xi, b_0, c_0)$ is odd admissible,  then 
\begin{align*}
&D\circ \xi= \xi^2=\xi^*\circ\xi=\xi\circ D=D^2=0, \;
b_0\in \ker{\xi}\cap \ker{D}, \,  c_0\in \ker\xi^*\cap \ker D,\\&
\langle b_0, c_0\rangle=\langle b_0, b_0\rangle=0  \esp D \text{ is $\prs$-antisymmetric}.
\end{align*} 
Since $\prs$ is even and $\mathcal{B}$ is two-dimensional, Lemma \ref{evendimension} implies that we must have $\mathcal{B} = \mathcal{B}_{\bar 0}$ or $\mathcal{B} =\mathcal{B}_{\bar 1}$. Moreover, because $D$ and $\xi$ are odd operators, they map homogeneous components to their opposite parity, which forces $\xi = D = 0$.

\underline{The case where $\mathcal{B} = \mathcal{B}_{\bar 0}$.} If $b_0 = 0$, we can choose an orthonormal  basis $\{e_1, e_2\}$ of $\mathcal{B}$ such that 
    $
    \prs = e_1^* \otimes e_1^* + \epsilon e_2^* \otimes e_2^*, $ where  $\epsilon = \pm 1,$ and $c_0 = \alpha e_1 + \beta e_2$ for some $\alpha, \beta \in \mathbb{R}$.
    
If $b_0 \neq 0$, the condition $\langle b_0, b_0 \rangle = 0$ implies that $b_0$ is isotropic (hence $\prs$ must be of Lorentzian signature). We set $e_1 = b_0$ and complete it to a hyperbolic basis $\{e_1, e_2\}$ satisfying $\langle e_1, e_2 \rangle = 1$ and $\langle e_1, e_1 \rangle_ = \langle e_2, e_2 \rangle = 0$. The orthogonality relation $\langle b_0, c_0 \rangle = 0$ then implies $c_0 = \alpha e_1$ for some $\alpha \in \mathbb{R}$.

    \underline{The case where $\mathcal{B} = \mathcal{B}_{\bar{1}}$.}
    Since $b_0$ and $c_0$ are even elements while $\mathcal{B}$ is purely odd, we immediately obtain $b_0 = c_0 = 0$.

This completes the proof.
\end{proof}

\sssbegin{Proposition} \label{solutioneven1}
Let $(\mathcal{B},\prs)$ be a two-dimensional odd pseudo-Euclidean  trivial Novikov superalgebra  \textup{(}i.e., with product $\mathcal{B} \bullet \mathcal{B} = \{0\}$\textup{)}.
\begin{itemize}
    \item[$(i)$] The triple $(D, \xi, b_0)$ is \emph{even admissible} if and only if there exists a basis $\mathbb{B} = \{e_1\mid  f_1\}$ of $\mathcal{B}$ such that the following conditions are satisfied:
 $$
         \xi=0, \quad D = 
        \begin{pmatrix}
        a & 0 \\
        0 & -a
        \end{pmatrix},\; b_0 = \alpha e_1 ,\quad \text{and}\quad \prs = e_1^*\odot f_1^*,\quad \text{where}\quad a, \alpha \in \mathbb{R}.$$

   \item[$(ii)$]  
The quadruple $(D,\xi, b_0, c_0)$ is \emph{odd admissible} if and only if there exists 
a homogeneous basis $\mathbb{B}=\{e_1 \mid f_1\}$ of $\mathcal{B}$ such that  
$\prs = e_1^{*} \odot f_1^{*}$, and one of the following holds:

\[
\begin{array}{lll}
\text{$(iia)$} & \xi=0, \quad \quad\quad D = \begin{pmatrix}
        0 & b \\
        0 & 0
        \end{pmatrix}, 
& b_0 = \alpha\, e_1, \; c_0=\beta e_1, \quad b, \al, \beta \in \mathbb{R},
\\[4pt]
\text{$(iib)$} & 
\xi = 
\begin{pmatrix}
0 & 0 \\[2pt]
a & 0
\end{pmatrix},\; D=0,\;
& b_0 = 0, \; c_0=0,\qquad  0\neq a \in \mathbb{R},
\\[10pt]
\text{$(iic)$} &
\xi = 
\begin{pmatrix}
0 & a \\[2pt]
0 & 0
\end{pmatrix},\; D = 
        \begin{pmatrix}
        0 & b \\
        0 & 0
        \end{pmatrix},
& b_0 = \alpha\, e_1,\; c_0=\beta e_1, \qquad 0 \neq a, b,  \al, \beta \in \mathbb{R}.
\end{array}
\]

\end{itemize}
\end{Proposition}

\begin{proof}
Let us prove Part (i). If the triple $(D, \xi, b_0)$ is even admissible,  then 
\[
\xi^2= D\circ \xi = \xi \circ D =0,  \qquad b_0 \in \ker( \xi) \esp D \text{ is $\prs$-antisymmetric}.
\]
Hence, $\xi$ is nilpotent.   Since $\prs$ is odd and $\xi$ and $D$ are even, there exists a homogeneous basis $\{e_1\mid f_1\}$ of $\mathcal{B}$ such that
\[
\prs = e_1^*\odot f_1^*,
\; 
\xi(e_1) = a e_1, \; 
\xi(f_1) = b f_1,\; D(e_1) = a_1 e_1, \; 
D(f_1) = b_1 f_1,\; \text{and}\;  b_0=\al e_1, 
\]
for some scalars $a, b, a_1, b_1, \al \in \mathbb{R}$.  Since $\xi$ is nilpotent, we must have $a = b = 0$, and therefore $\xi = 0$. Since $D$ is $\prs$-antisymmetric, we deduce that $b_1=-a_1.$

Let us prove Part (ii). If the quadruple $(D, \xi, b_0, c_0)$ is odd admissible,  then 
\begin{align*}
&D\circ \xi= \xi^2=\xi^*\circ\xi=\xi\circ D=D^2=0, \;
b_0\in \ker{\xi}\cap\ker{\xi^*}\cap \ker{D}, \,  c_0\in \ker\xi^*\cap \ker D,\\&
\langle b_0, c_0\rangle=\langle b_0, b_0\rangle=0  \esp D \text{ is $\prs$-antisymmetric}.
\end{align*}  
Since $\prs$ is odd and $\xi$ and $D$ are odd, there exists a homogeneous basis $\{e_1\mid f_1\}$ of $\mathcal{B}$ such that
\[
\prs = e_1^*\odot f_1^*,
\; 
\xi(e_1) = a f_1, \; 
\xi(f_1) = b\, e_1,\; D(e_1) = a_1\, f_1, \; 
D(f_1) = b_1\, e_1,\; b_0=\al e_1, \text{and}\;  c_0=\beta e_1, 
\]
for some $a, b, a_1, b_1, \alpha, \beta \in \mathbb{R}$. The $\prs$-antisymmetry of $D$ implies $a_1 = 0$, while $\xi^2 = 0$ requires $ab = 0$. We distinguish three cases:
\begin{itemize}
   \item[\emph{Case 1.}] Suppose \(a=b=0\). Then \(\xi=0\), while \(b_1,\alpha,\beta\in\mathbb{R}\) are arbitrary. Hence, we obtain~(iia).

   \item[\emph{Case 2.}] Suppose \(a\neq 0\) and \(b=0\). The condition \(\xi\circ D=0\) yields \(ab_1=0\), and hence \(b_1=0\), so \(D=0\). Since \(b_0,c_0\in\ker\xi\), we have \(\alpha a=\beta a=0\), which implies \(\alpha=\beta=0\). Hence, we obtain~(iib).

   \item[\emph{Case 3.}] Suppose \(a=0\) and \(b\neq 0\). All admissibility conditions are satisfied for arbitrary \(b_1,\alpha,\beta\in\mathbb{R}\). Hence, we obtain~(iic).\qedhere

\end{itemize}

\end{proof}

\sssbegin{Proposition}\label{A22}
Let $(\A_2^2, \bullet, \prs)$ be  the odd pseudo-Euclidean Novikov superalgebra. Then 
\begin{itemize}
    \item[$(i)$] The triple $(D, \xi, b_0)$ is \emph{even admissible} if and only if there exists a basis $\mathbb{B} = \{e_1\mid  f_1\}$ of $\mathcal{B}$ such that the following conditions are satisfied:
 $$
         \xi=D=0,\; b_0 = \beta e_1 ,\quad \text{where}\quad  \beta\in \mathbb{R}.$$
   \item[$(ii)$] The quadruple  $(D, \xi, b_0, c_0)$ is \emph{odd admissible} if and only if there exists a basis $\mathbb{B} = \{e_1\mid  f_1\}$ of $\mathcal{B}$ such that the following conditions are satisfied:     
    $$
         \xi=\begin{pmatrix}
        0 & a \\
        0 & 0
        \end{pmatrix}, \quad D = 
        \begin{pmatrix}
        0 & b \\
        0 & 0
        \end{pmatrix},\; b_0 = \beta e_1 ,\quad  \text{and}\quad c_0=\la e_1,\quad \text{where}\quad a, b, \beta,\la \in \mathbb{R}.$$
        \end{itemize}

\end{Proposition}
\begin{proof}
Let us only prove Part (i). According to Proposition~\ref{Dim2}, there exists a basis 
$\{e_1 \mid f_1\}$ of $\A_2^2$ 
whose only non-zero product and non-degenerate bilinear form are
\[
f_1 \bullet f_1 = \alpha e_1, 
\qquad 
\prs = e_1^* \odot f_1^*,
\]
where $\alpha \neq 0$. If the triple $(D, \xi, b_0)$ is even admissible, then
    \begin{equation*}
    \begin{cases}
    \xi(u\bullet v)=\xi(u)\bullet v=D(u)\bullet v=0, \\
    u\bullet \xi(v)=(-1)^{|u||v|} v\bullet \xi(u),\quad D(u\bullet v) =u\bullet D(v), \\[1mm]
    D\circ \xi=\Rr_{b_0}^\bullet,\quad \xi^2=\xi\circ D=\Ll^\bullet_{b_0}=0, \quad 
    b_0\in \ker{\xi}.
    \end{cases}
    \end{equation*}
    Since $0=\xi(f_1\bullet f_1)=\alpha \,\xi(e_1)$, it follows that $\xi(e_1)=0$. 
    Since $\xi$ is even and $\xi^2=0$, we conclude that $\xi=0$. 
    Since $D$ is even, we have $D(e_1)=\beta e_1$ and $D(f_1)=\lambda f_1$ 
    with $\beta, \lambda \in \mathbb{R}$. 
    From $0=D(f_1)\bullet f_1=\lambda \alpha e_1$, we deduce $\lambda=0$. 
    Moreover, $$\alpha D(e_1) = D(f_1\bullet f_1) = f_1 \bullet D(f_1) = 0,$$ 
    which implies $D(e_1)=0$, hence $D=0$. 
    Since $b_0$ is even, we have $b_0=\gamma e_1$ for some $\gamma\in \mathbb{R}$.
\end{proof}

\sssbegin{Proposition}
Let $(\A, \bullet, \prs)$ be an even pseudo-Euclidean Novikov superalgebra of dimension $3$ such that $\A\bullet\A$ is degenerate,  then $(\A,\bullet,\prs)$ is isomorphic to one of the following pseudo-Euclidean superalgebras
    \begin{enumerate}
        \item[$(i)$] $(\A^3_5,\bullet, \prs)$, given in the basis $\{e_1, e_2, e_3\}$ by 
        \[
e_1\bullet e_2=-\la e_3,\qquad
e_1\bullet e_1=\la e_2,\quad \prs=e_1\odot e_3+ e_2\otimes e_2, \quad \text{where $\la\neq 0$.}
\] 
        \item[$(ii)$]  $(\A^3_6,\bullet, \prs)$, given in the basis $\{e_1\mid e_2, e_3\}$ by  
        \[
e_1\bullet e_2=\la e_3,  \quad \prs=e_1\otimes e_1+e_2\odot e_3, \quad \text{where $\lambda\neq0$.}\]
       
    \end{enumerate}

\end{Proposition}
\begin{proof}
Let $(\A, \bullet, \prs)$ be a pseudo-Euclidean Novikov superalgebra such that
$\A \bullet \A$ is degenerate. 
 According to Theorem~\ref{central}, such a superalgebra is either an even or an odd double extension
of the one-dimensional trivial superalgebra
$(\mathcal{B}, \bullet_{\mathcal{B}}, \prs_{\mathcal{B}})$
by means of $(D, \xi, b_0)$ or $(D, \xi, b_0, c_0)$, respectively.
Moreover, by Lemma~\ref{evendimension}, we have $\mathcal{B} = \mathcal{B}_{\bar{0}}$.

\medskip
\noindent
\underline{\textbf{The even admissible case.}}
Let $\mathcal{B} = \mathbb{R} e_1$ with
$\langle e_1, e_1 \rangle_{\mathcal{B}} = 1$ and $b_0 = \lambda e_1$, where $\lambda \in \mathbb{R}$.
Since $\xi^2 = 0$ and $D$ is $\prs_{\mathcal{B}}$-antisymmetric, it follows that $D = \xi = 0$.
Eq.~\eqref{Produit1}  implies that there exists a basis $\{e, d, e_1\}$ of $\A$
such that the Novikov product and the bilinear form are given by
\[
d \bullet e_1 = -\lambda e, \qquad
d \bullet d = \lambda e_1, \qquad
\prs = e^* \odot d^* + e_1^*\otimes e_1^*.
\]
Therefore, $(\A, \bullet, \prs)$ is isomorphic to
$(\A^3_5, \bullet,  \prs)$.

\medskip
\noindent
\underline{\textbf{The odd admissible case.}}
Let $\mathcal{B} = \mathbb{R} e_1$ with
$\langle e_1, e_1 \rangle_{\mathcal{B}} = 1$,
$b_0 = \lambda e_1$, and $c_0 = \alpha e_1$, where $\lambda, \alpha \in \mathbb{R}$.
Since $\xi^2 = 0$, $D$ is $\prs_{\mathcal{B}}$-antisymmetric, and
$\langle b_0, b_0 \rangle_{\mathcal{B}} = 0$, we obtain $D = \xi = 0$ and $b_0=0$.
Eq.~\eqref{Produit2} implies that there exists a homogeneous basis
$\{e_1 \mid e, d\}$ of $\A$ such that
\[
e_1 \bullet d = \alpha e, \qquad
\prs = e^* \odot d^* + e_1^*\otimes e_1^*.
\]
Hence, $(\A, \bullet, \prs)$ is isomorphic to
$(\A^3_6, \bullet,  \prs)$.
\end{proof}

\sssbegin{Theorem}
Every 4-dimensional even pseudo-Euclidean Novikov superalgebra $(\mathcal{A}, \bullet, \prs)$ such that $\mathcal{A} \bullet \mathcal{A}$ is degenerate  is isomorphic to one of the following superalgebras:
		\begin{enumerate}
			\item[$(i)$] $(\A^4_7,\bullet,\prs)$ given in the basis $\{e,e_1,e_2,d\}$ by 
			\[\begin{cases} d\bullet e_1=-a\epsilon e_2-\al e,\quad  d\bullet e_2=ae_1-\beta\epsilon e\esp d\bullet d=\al e_1+\beta e_2, \; \text{ where }a, \al, \beta\in\mathbb{R},\\
				\prs=e^*\odot d^*+ e_1^*\otimes e_1^*+\epsilon e_2^*\otimes e_2^*,\quad\epsilon=\pm 1.			\end{cases}	\] 

        \item[$(ii)$] $(\A_8^4,\bullet,\prs)$ given in the basis $\{e,e_1,e_2,d\}$ by  
		\[\begin{cases}
			e_2\bullet e_1=-a\epsilon e,\quad  e_2\bullet d=ae_1,\quad  d\bullet e_1=-\al\epsilon e, \esp d\bullet d=\al e_1, \text{ where } 0 \not =a, \al\in \mathbb{R},\\
			\prs=e^*\odot d^*+\epsilon e_1^*\otimes e_1^*+\rho e_2^*\otimes e_2^*, \text{ where }\epsilon,\rho=\pm 1.
		\end{cases}
		\]
		\item[$(iii)$] $(\A_9^4,\bullet,\prs)$ given in the basis $\{e,e_1,e_2,d\}$ by  
		\[\begin{cases}
			e_2\bullet e_2=-a e,\quad  e_2\bullet d=ae_1, \quad d\bullet e_2=-\al e, \esp d\bullet d=\al e_1, \text{ where } \; 0 \not =a,  \al\in \mathbb{R},\\
			\prs=e^*\odot d^*+e_1^*\odot e_2^*. 

		\end{cases}
		\]

        \item[$(iv)$] $(\A_{10}^4,\bullet,\prs)$ given in the basis $\{e,d\mid e_1,e_2\}$ by  
		\[\begin{cases}
			d\bullet e_1=a e_1+c e_2, \esp d\bullet e_2=be_1-ae_2, \text{ where }  a,b, c\in \mathbb{R},\\
			\prs=e^*\odot d^*-e_1^*\odot e_2^*. 

		\end{cases}
		\]

        \item[$(v)$] $(\A_{11}^4,\bullet,\prs)$ given in the basis $\{e,d\mid e_1,e_2\}$ by 
		\[\begin{cases}
			e_2\bullet e_2=-a e,\; e_2\bullet d=ae_1,\esp d\bullet e_2=be_1, \text{ where } 0 \not =a, b\in \mathbb{R}, \\
			\prs=e^*\odot d^*-e_1^*\odot e_2^*.

		\end{cases}
		\]
         \item[$(vi)$] $(\A_{12}^4,\bullet,\prs)$ given in the basis $\{e_1,e_2\mid e, d\}$ by  
		\[\begin{cases}
			e_1\bullet d=\al e, \esp e_2\bullet d=\epsilon\beta e,\text{ where }   \al, \beta\in \mathbb{R},\\
			\prs=e_1^*\otimes e_1^*+\epsilon e_2^*\otimes e_2^*-e^*\odot  d^*. \end{cases}
		\]
        
             \item[$(vii)$] $(\A_{13}^4,\bullet,\prs)$ given in the basis $\{e_1,e_2\mid e, d\}$ by
		\[\begin{cases}
			e_2\bullet d=\al e,  \, d\bullet d=\la e_1\esp d\bullet e_2=-\la e \text{ where }   \la\neq 0 \text{ and }\al\in \mathbb{R},\\
			\prs=e_1^*\odot e_2^*-e^*\odot d^*. 

		\end{cases}
		\]
       
		\end{enumerate}
		\end{Theorem}

\begin{proof}
According to Theorem~\ref{central}, such a superalgebra is either an even or an odd double extension
of the two-dimensional trivial superalgebra
$(\mathcal{B}, \bullet_{\mathcal{B}}, \prs_{\mathcal{B}})$
by means of $(D, \xi, b_0)$ or $(D, \xi, b_0, c_0)$, respectively.

By Proposition~\ref{solutioneven}, all admissible solutions of the data
$(D, \xi, b_0)$ and $(D, \xi, b_0, c_0)$ are explicitly determined.
Applying the Novikov product associated with the corresponding double extension,
we obtain the above  superalgebras.
\end{proof}

\sssbegin{Theorem}
Every 4-dimensional odd pseudo-Euclidean Novikov superalgebra $(\mathcal{A}, \bullet, \prs)$ such that $\mathcal{A} \bullet \mathcal{A}$ is degenerate is isomorphic to one of the following superalgebras:
		\begin{enumerate}
			\item[$(i)$] $(\A^4_{14},\bullet,\prs)$ given in the basis $\{d,e_1\mid f_1,e\}$ by 
			\[\begin{cases} d\bullet e_1=a e_1,\; d\bullet f_1=-af_1-\al e,  \esp d\bullet d=\al e_1, \quad \al \in\mathbb{R},\\
				\prs=e^*\odot d^*+ e_1^*\odot f_1^*.			\end{cases}	\]

        \item[$(ii)$] $(\A_{15}^4,\bullet,\prs)$ given in the basis $\{e,e_1\mid f_1,d\}$ by 
		\[\begin{cases}
			d\bullet f_1=be_1+\al e,\;  f_1\bullet d=\beta e\esp d\bullet d=\al e_1, \text{ where } b, \al, \beta\in \mathbb{R},\\
			\prs=e^*\odot d^*+ e_1^*\odot f_1^*. 
		\end{cases}
		\]
       
		\item[$(iii)$] $(\A_{16}^4,\bullet,\prs)$ given in the basis $\{e,e_1\mid f_1,d\}$ by 
		\[\begin{cases}
			e_1\bullet d=a f_1, \esp e_1\bullet e_1=-a e, \text{ where } 0 \not =a\in \mathbb{R},\\
			\prs=e^*\odot d^*+ e_1^*\odot f_1^*.
		\end{cases}
		\]
        \item[$(iv)$] $(\A_{17}^4,\bullet,\prs)$ given in the basis $\{e,e_1\mid f_1,d\}$ by 
		\[\begin{cases}
			d\bullet f_1=b e_1+\al e,\; f_1\bullet d=ae_1+\beta e, \quad  f_1\bullet f_1=a  e, \esp d\bullet d=\al e_1, \text{ where } \; 
			0\not =a, b, \al, \beta\in \mathbb{R},\\
			\prs=d^*\odot e^*+ e_1^*\odot f_1^*.
		\end{cases}
		\]
         \item[$(v)$] $(\A_{18}^4,\bullet,\prs)$ given in the basis $\{d,e_1\mid f_1,e\}$ by 
		\[\begin{cases}
			d\bullet f_1=-\beta e,\; f_1\bullet f_1=\al e_1\quad   \esp d\bullet d=\beta e_1, \text{ where } \; 
			0\not =\al, \beta\in \mathbb{R},\\
			\prs=d^*\odot e^*+ e_1^*\odot f_1^*.
		\end{cases}
		\]
       
         \item[$(vi)$] $(\A_{19}^4,\bullet,\prs)$ given in the basis $\{e,e_1\mid f_1, d\}$ by 
		\[\begin{cases}
			d\bullet f_1=b e_1+\beta e, \quad   f_1\bullet d=ae_1+\la  e, \quad  f_1\bullet f_1=\al e_1+a e, \\  d\bullet d=\beta e_1, \text{ where } 0\not=\alpha, \la, \beta, a,b\in \mathbb{R},\\
			\prs=e^*\odot d^*+ e_1^*\odot f_1^*.
		\end{cases}
		\]
        \end{enumerate}

\end{Theorem}

\begin{proof}
According to Theorem~\ref{central1}, such a superalgebra is either an even or an odd double extension
of the two-dimensional trivial or non-trivial Novikov superalgebra
$(\mathcal{B}, \bullet_{\mathcal{B}}, \prs_{\mathcal{B}})$
by means of $(D, \xi, b_0)$ or $(D, \xi, b_0, c_0)$, respectively.

By Proposition~\ref{solutioneven1} and Proposition \ref{A22}, all admissible solutions of the data
$(D, \xi, b_0)$ and $(D, \xi, b_0, c_0)$ are explicitly determined.
Applying the Novikov product associated with the corresponding double extension,
we obtain the above superalgebras.
\end{proof}
\sssbegin{Remark}
The list of pseudo-Euclidean Novikov superalgebras we have provided is exhaustive, but may contain redundancies. It would be interesting to check for which values of the parameters these algebras are isomorphic to each other. Carrying out this classification does not seem straightforward and requires a number of computational techniques that are beyond the scope of this paper.
\end{Remark}

{\bf Acknowledgment.} We are very grateful to the anonymous referee for the careful reading of the manuscript, which considerably improved the proofs of some statements, as well as the presentation.



\begin{thebibliography}{}
\bibitem{ABL}
   Ait Ben Haddou M.,  Boucetta M., and  Lebzioui H., \emph{Left-invariant Lorentzian Flat Metrics on Lie Groups.} Journal of Lie Theory 22.1 (2012): 269-289.
   \bibitem{AAO}
Albeverio S.,  Ayupov Sh A., and  Omirov B. A., On nilpotent and simple Leibniz algebras. Communications in Algebra 33.1 (2005): 159-172.
\bibitem{Medina2}
 Aubert A.,   Medina A., \emph{Groupes de Lie pseudo-riemanniens plats.} Tohoku Mathematical Journal, Second Series 55.4 (2003): 487-506.
  \bibitem{AB}
  Ayadi I.,   Benayadi S., \emph{Symmetric Novikov superalgebras.} Journal of mathematical physics 51.2 (2010).
\bibitem{AFMV}
Autenried C., Furutani K., Markina I., Vasiľev, A., Pseudo-metric 2-step nilpotent Lie algebras. Advances in Geometry, 18(2),(2018)  237-263.




\bibitem{Bajo}
Bajo, I,  Benayadi S.  Lebzioui  H., Classification of flat Lorentzian nilpotent Lie algebras. Bulletin of the London Mathematical Society 56.6 (2024): 2132-2149.








 



 
 





\bibitem{BBE}
Benayadi S.,  Bouarroudj S.,  El Ouali H., Flat pseudo-Euclidean Leibniz superalgebras. arXiv preprint arXiv:2510.15328 (2025).
\bibitem{BO}
Benayadi S.,  Oubba H., Complete descriptions of pseudo-Euclidean left-symmetric $\Ll$-algebras and their pseudo-Euclidean modules. Journal of Algebra 657 (2024): 600-637.
\bibitem{BO1}
Benayadi S.,  Oubba H., Nonassociative algebras of biderivation-type. Linear Algebra and its Applications 701 (2024): 22-60.
\bibitem{Bordemann}
 Bordemann M., \emph{Nondegenerate associative bilinear forms on nonassociative algebras}. Acta Math.
Univ. Com. LXVI(2), (1997) :151-201.
\bibitem{BoO} Bouarroudj S.,  El Ouali H., Flat quasi-Frobenius Lie superalgebras,  European Journal of Mathematics, Vol. {\bf 12}, No. 4, 2026. \url{https://doi.org/10.1007/s40879-026-00889-2} 

\bibitem{BoO2} Bouarroudj S.,  El Ouali H., Structure of flat quadratic quasi-Frobenius Lie superalgebras via double extensions, J. Algebra Appl. (2026), 2850020.

\bibitem{BR}  
Bouarroudj S.,  Radu A-M.,
Lagrangian extensions and left-symmetric structures on the four-dimensional real Lie superalgebras, Preprint.
\url{arXiv:2410.01520} To appear in the Ukrainian Mathematical Journal.
\bibitem{BEL}
Boucetta M.,  El Ouali H.,   Lebzioui H., Pseudo-Euclidean Novikov algebras of arbitrary signature. Journal of Geometry and Physics 206 (2024): 105334.
\bibitem{BL1}
 Boucetta M.,  Lebzioui H., \emph{Flat nonunimodular Lorentzian Lie algebras.} Communications in Algebra 44.10 (2016): 4185-4195.
\bibitem{BL2}
 Boucetta M., Lebzioui H., \emph{On flat pseudo-Euclidean nilpotent Lie algebras.} Journal of Algebra 537 (2019): 459-477.
 \bibitem{BDD}Burde D.,  Deschamps S., and Dekimpe K., LR-algebras. Contemporary Mathematics (CONM) 491 (2009): 125--140.

\bibitem{CP}
Cordero, L. A.,  Parker P. E., Pseudoriemannian 2-step nilpotent Lie groups. arXiv preprint math/9905188 (1999).


\bibitem{DI}
Duncan D. C., Ihrig E. C., Homogeneous spacetimes of zero curvature.
Proceedings of the American Mathematical Society (1989) 107(3):785-795.


\bibitem{Patrick1}
Eberlein P., Geometry of $2 $-step nilpotent groups with a left invariant metric. Annales scientifiques de l'Ecole normale supérieure. Vol. 27. No. 5. 1994.

\bibitem{Patrick2}
Eberlein, P., Geometry of 2-step nilpotent groups with a left invariant metric. II. Transactions of the American Mathematical Society 343.2 (1994): 805-828.



\bibitem{G}
 Guédiri M., Novikov algebras carrying an invariant Lorentzian symmetric bilinear form, J. Geom.
Phys. 82 (2014) 132–144.
\bibitem{guediri}
	 Guédiri M., \emph{On the nonexistence of closed timelike geodesics in flat Lorentz
	2-step nilmanifolds}.  Transactions of the American Mathematical Society 355.2
	(2003), pp. 775–786.

    \bibitem{JPP}
    Jang C.,   Parker P. E., Conjugate loci of pseudo-Riemannian 2-step nilpotent Lie groups with nondegenerate center. Annals of Global Analysis and Geometry 28.1 (2005): 1-18.
  \bibitem{KKK} Karimjanov  I., Kaygorodov I., and Khudoyberdiyev A., The algebraic and geometric classification of nilpotent Novikov algebras. J. Geom. Phys. 143 (2019), 11--21.
    \bibitem{Ko}
    Kostant B., The Weyl algebra and the structure of all Lie superalgebras of Riemannian type. Transformation Groups 6.3 (2001): 215-226.
\bibitem{LH1}
 Lebzioui H., On a remarkable class of flat pseudo-Riemannian Lie groups, J. Geom. Phys. 186 (2023)
104774.

\bibitem{LH}
Lebzioui H., Flat left-invariant pseudo-Riemannian metrics on unimodular Lie groups. Proceedings of the American Mathematical Society 148.4 (2020): 1723-1730.
\bibitem{LH2}
Lebzioui H., On pseudo-Euclidean Novikov algebras. Journal of Algebra 564 (2020): 300-316.

\bibitem{L}
Leites D. (ed.) \textit{Seminar on supersymmetry v. $1$. Algebra and
Calculus: Main chapters}, (J.~Bernstein, D.~Leites, V.~Molotkov,
V.~Shander), MCCME, Moscow, 2012, 410 pp (in Russian; a~version in
English is in preparation but available for perusal).



 
\bibitem{Youri}
 Manin  Y.I., \emph{ Gauge field theory and complex geometry.} Vol. 289. Springer Science and Business Media, 1997.

 \bibitem{MR}
	 Medina A.,  Revoy P., \emph{ Alg\`ebres de Lie et produit scalaire invariant.} Annales scientifiques de l'\'{E}cole Normale Sup\'erieure {\bf 18} No. 3 (1986) 553-561.
    
\bibitem{M}  Milnor J., \emph{Curvatures of left invariant metrics on Lie groups}, Adv. Math. 21(3) (1976) 293-329.

\bibitem{JZ}
 Junna NI.,  Zhiqi C., \emph{Novikov super-algebras with associative non-degenerate super-symmetric bilinear forms.} Journal of Nonlinear Mathematical Physics 17.2 (2010): 159-166.

\bibitem{WG}
Walschap G., Cut and conjugate loci in two-step nilpotent Lie groups. The Journal of Geometric Analysis 7.2 (1997): 343-355.
 \bibitem{Xu}
   Xu X., Variational calculus of supervariables and related algebraic structures, J. Algebra 223 (2000) 396-437.

\bibitem{Z}
  Zelmanov E.I., On a class of local translation invariant Lie algebras, Sov. Math. Dokl. 35 (1987)
216–218.

\end{thebibliography}
\end{document}